\documentclass{article}

\usepackage{enumerate}
\date{}
\usepackage{amsfonts,amsmath,amsthm}
\usepackage{amssymb,latexsym}
\usepackage[dvips]{epsfig}
\usepackage{amscd}  
\usepackage[all]{xy}
\usepackage{mathrsfs,dsfont}
\title{\bf A categorification of finite-dimensional irreducible representations
of quantum $sl(2)$ and their tensor products}
\vspace{5cm}
\author{Igor Frenkel, Mikhail Khovanov and Catharina Stroppel} 

\newtheorem{theorem}{Theorem}[section]
\newtheorem{prop}[theorem]{Proposition}
\newtheorem{lemma}[theorem]{Lemma}
\newtheorem{exs}[theorem]{Example}
\newtheorem{corollary}[theorem]{Corollary}
\newtheorem{definition}[theorem]{Definition}
\newtheorem{conjecture}[theorem]{Conjecture} 
\newtheorem{remark}[theorem]{Remark}

\newcommand{\ov}{\overline}
\def\eg{\emph{e.g. }}
\def\ie{\emph{i.e. }}
\newcommand{\mf}{\mathfrak}
\newcommand{\mc}{\mathcal}
\newcommand{\ID}{\operatorname{ID}}
\newcommand{\cG}{\operatorname{G}}
\newcommand{\cL}{\mathcal{L}}
\newcommand{\cB}{\mathcal{B}}
\newcommand{\cM}{\mathcal{M}}
\newcommand{\cP}{\mathcal{P}}
\newcommand{\cA}{\mathcal{A}}
\newcommand{\mR}{\mathbb{R}}
\newcommand{\cH}{\mathcal{H}}
\newcommand{\Z}{\mathbb{Z}}
\newcommand{\mZ}{\mathbb{Z}}
\newcommand{\mQ}{\mathbb{C}}

\newcommand{\mC}{\mathbb{C}}

\newcommand{\Q}{\mathbb{C}}
\newcommand{\Qq}{\mathbb{C}(q)}
\newcommand{\cO} {\mathcal{O}}

\newcommand{\cU}{\mathcal{U}}
\newcommand{\cZ}{\mathcal{Z}}
\newcommand{\mg}{\mathfrak{g}}
\newcommand{\mb}{\mathfrak{b}}
\newcommand{\cOn}{\cO(\mathfrak{gl}_n)}
\newcommand{\mh}{\mathfrak{h}}
\newcommand{\la}{\lambda}
\newcommand{\cC}{\mathcal{C}}
\newcommand{\cE}{\mathcal{E}}
\newcommand{\cF}{\mathcal{F}}

\newcommand{\mV}{\mathbb{V}}
\newcommand{\op}{\operatorname}
\newcommand{\KER}{\operatorname{ker}}
\newcommand{\END}{\operatorname{End}}
\newcommand{\HOM}{\operatorname{Hom}}
\newcommand{\RES}{\operatorname{Res}}
\newcommand{\MOD}{\operatorname{-mod}}

\newcommand{\gMOD}{\operatorname{-gmod}}

\newcommand{\DIM}{\operatorname{dim}}
\renewcommand{\H}{C}
\newcommand{\GL}{\mathfrak{gl}}
\newcommand{\SL}{\mathfrak{sl}}
\renewcommand{\blacktriangle}{\triangle}

\begin{document}
\ifx\href\undefined\else\hypersetup{linktocpage=true}\fi 
\maketitle
\begin{abstract}
  The purpose of this paper is to study categorifications of tensor products
  of finite dimensional modules for the quantum group $U_q(\mathfrak{sl}_2)$. The main
  categorification is obtained using certain Harish-Chandra bimodules for the
  complex Lie algebra $\mathfrak{gl}_n$. For the special case of simple
  modules we naturally deduce a categorification via modules over the
  cohomology ring of certain flag varieties. Further geometric
  categorifications and the relation to Steinberg varieties is discussed. We
  also give a categorical version of the quantised Schur-Weyl duality and an interpretation of the (dual) canonical bases and the (dual) standard
  bases in terms of projective, tilting, standard and simple Harish-Chandra bimodules.
\end{abstract}

\tableofcontents

\section*{Introduction} 
A categorification program of the simplest quantum group,
$U_q(\mathfrak{sl}_2)$, has been formulated in \cite{BFK}. Its ultimate goal is
to construct a certain tensor $2$-category with a Grothendieck ring equivalent
to the representation category of $U_q(\mathfrak{sl}_2)$. In \cite{BFK},
the categorification problem was studied for modules over the classical algebra
$\cU(\mathfrak{sl}_2)$. There, two
versions of categorification of the $\cU(\mathfrak{sl}_2)$-action on
$V_1^{\otimes n}$, the $n$-th power of the two-dimensional fundamental
representation for $\mathfrak{sl}_2$, were obtained: one using certain
singular blocks of the category $\cO(\mathfrak{gl}_n)$ of highest weight
$\mathfrak{gl}_n$-modules, another one via certain parabolic subcategories of the
regular block of $\cO(\mathfrak{gl}_n)$. The quantum versions of these two
categorifications were conjectured in \cite{BFK} and established for
the parabolic case in \cite{Sduke} using the graded version of the category
$\cO(\mathfrak{gl}_n)$ from \cite{Sperv} and \cite{BGS} and graded lifts of
translation functors. In the present paper, among other results, we develop a categorification of
the $U_q(\mathfrak{sl}_2)$-action on $V_1^{\otimes n}$ using (a graded version
of) certain singular blocks of $\cOn$ (see Section~\ref{Ograd}). It was also
conjectured in \cite{BFK} that the two categorifications, the one using parabolic
subcategories of the regular block of $\cOn$ and the one using singular blocks of $\cOn$, are related by the
Koszul duality functor of \cite{BGS}. In fact, it was shown in \cite{Steen}
that the key functors used in both categorifications, namely graded versions
of translation functors and Zuckerman functors, are Koszul dual to each
other. This completes the general picture of categorifications of the
$U_q(\mathfrak{sl}_2)$-action in $V_1^{\otimes n}$ initiated in \cite{BFK} and
opens a way for further steps in the general program. \\

The main purpose of the
present paper is to study the categorification of the tensor products of
arbitrary finite dimensional irreducible $U_q(\mathfrak{sl}_2)$-modules (as recalled in
Section~\ref{quantum-sl2}), \ie of modules of the form 
\begin{eqnarray}
  \label{tensor}
  V_{\bf d}=V_{d_1}\otimes V_{d_{2}}\otimes\cdots\otimes
V_{d_r}, 
\end{eqnarray}
where $V_{d_i}$ denotes an irreducible representation of dimension
$d_i+1$. Our basic observation is that $V_{\bf d}$ admits a categorification via
blocks of the category of Harish-Chandra bimodules $\cH$ for
$\mathfrak{gl}_n$, where $n=\sum_{i=1}^rd_i$. To make this more precise, we go
back and recall, from \cite{BFK}, that the $\cU(\mathfrak{sl}_2)$-module 
$V_{1}^{\otimes n}$ has a categorification via the category
\begin{eqnarray}
\label{O}
  \bigoplus_{i=0}^{n}{}_{\omega_i}\cO(\mathfrak{gl}_n),
\end{eqnarray}
where $_{\omega_i}\cO(\mathfrak{gl}_n)$ is the block corresponding to an
integral weight $\omega_i$ with stabiliser $S_i\times S_{n-i}$. We establish a
categorification of the corresponding $U_q(\mathfrak{sl}_2)$-module using a
graded version of the category~\eqref{O}. The categorification of an arbitrary tensor
product $V_{\bf d}$ is then given by a graded version of  
\begin{eqnarray*}
  \cH_\mu:=\bigoplus_{i=0}^n{}_{\omega_i}\cH_\mu, 
\end{eqnarray*}
where $\omega_i$, $i=0,\ldots, n$ are as above, and $\mu$ is a dominant integral weight with stabiliser isomorphic to
 \begin{eqnarray*}
   S_{\bf d}=S_{d_1}\times S_{d_{2}}\times\cdots\times S_{d_r}.
 \end{eqnarray*}
The connection of our categorification with the one in \cite{BFK} is given by
a functor, introduced in \cite{BG}, namely 
the functor of tensoring with some Verma module, $M(\mu)$:  
\begin{eqnarray}
\label{BGGfunc}
  {}_\la F_\mu:\quad {}_\la\cH_\mu^1&\rightarrow& {}_\la\cOn\\
 X&\mapsto& X\otimes_{\cU(\mathfrak{gl}_n)}M(\mu),\nonumber
\end{eqnarray}
for dominant integral $\la$ and $\mu$. This functor defines an equivalence of
categories for regular $\mu$. If $\mu$ is not regular, then ${}_\la F_\mu$
yields an equivalence with a certain full subcategory of ${}_\la\cO$. (Note that,
although 
${}_\la\cH_\mu^1$ is in fact a full subcategory of ${}_\la\cH_\mu$, they
have the same Grothendieck group.)

In the present paper we study a categorification based on singular blocks of
Harish-Chandra bimodules. By analogy with \cite{BFK} one expects the existence
of a second categorification that uses parabolic subcategories, related to the
first by (a certain generalisation of) the Koszul duality functor. (The first
steps in this direction will appear in \cite{MOS}.)\\
To achieve our main goal, we first construct a categorification of tensor
products of the form~\eqref{tensor} for the classical algebra $\cU(\mathfrak{sl}_2)$ in
Section~\ref{catO}. After recalling and further developing the graded version of the category $\cOn$
in Section~\ref{Ograd} we obtain the categorification of the tensor
products~\eqref{tensor} for the quantum algebra $U_q(\mathfrak{sl}_2)$ in
Section~\ref{mainresult}. 

Our approach to the grading of the category $\cO$ and the category of
Harish-Chandra modules is based on the Soergel functor (\cite{Sperv}) 
\begin{eqnarray*}
  _\la\mV:{}_\la\cO\rightarrow C^\la,
\end{eqnarray*}
where $C^\la$ is the algebra of endomorphisms of the unique indecomposable
projective tilting module in ${}_\la\cO$. One of the main results of
\cite{Sperv} is an explicit description of
$C^\la$ as a ring of invariants in coinvariants. The latter is known to be
isomorphic to the cohomology ring of the partial flag variety $\cF_\la$
corresponding to $\la$. In particular, there is a natural $\mZ$-grading on
$C^\la$. Moreover, it was shown in \cite{Sperv} that the functor
${}_\la\mV$ intertwines the translation
functors $\mathds{T}_\la^\mu:{}_\la\cO\rightarrow{}_\mu\cO$ in the category
$\cO$ with the induction $C^\mu\otimes_{C^\la}\bullet:C^\mu\rightarrow C^\la$
or with the restriction functor $\RES_\mu^\la:C^\la\MOD\rightarrow
C^\mu\MOD$, depending on whether the stabiliser of $\mu$ is contained
in the stabiliser of $\la$ or vice versa. This fact is employed for
constructing the functors that provide a categorification of the
$\cU(\mathfrak{sl}_2)$-action. The natural grading of the algebra $C^\la$ is
used to define graded lifts of the previous functors which yield the quantum counterpart. \\

The categorification of tensor products~\eqref{tensor} of
$U_q(\mathfrak{sl}_2)$-modules via Harish-Chandra bimodules gives rise to
special bases resulting from indecomposable tilting, simple, standard and dual
standard modules. In Section~\ref{SchurWeylduality} we prove in the special
case of $V_1^{\otimes n}$ (\ie $d_1=d_2=\cdots =d_n=1$) that, at the
Grothendieck level, these bases can be identified as canonical, dual
canonical, standard and dual standard bases in $V_1^{\otimes n}$. This is
established using the Kazhdan-Lusztig theory, Schur-Weyl duality and the
graphical calculus for tensor products. A generalisation of these results to
arbitrary tensor products requires a more extensive study of the category of
Harish-Chandra bimodules which will be postponed to a future paper.\\ 

Our main theorem about the categorification of the tensor products of the form~\eqref{tensor}
in the special case of a single factor (\ie $r=1$), combined with the
properties of Soergel's functor leads to a ``geometric'' categorification of
irreducible representations of $U_q(\mathfrak{sl}_2)$ which we explain in
Section~\ref{geometry}. We conclude this section with a discussion and a
conjecture about a possible geometric categorification of the general tensor
product, based on the Borel-Moore homology of the generalised Steinberg
varieties $X^{\la,\mu}$ of triples (see \cite{RD1}). The proposed geometric
approach is connected with our geometric categorification of simple
$U_q(\mathfrak{sl}_2)$-modules, since in the special case $r=1$,
the generalised Steinberg variety degenerates into the partial flag variety
$\cF_\la$. It is important to note that the generalised Steinberg varieties
$X^{\la,\mu}$ are precisely the tensor product varieties of Malkin (\cite{Malkin}) and
Nakajima (\cite{Nakajima}) in the special case of tensor products of finite
dimensional irreducible $U_q(\mathfrak{sl}_2)$-modules (see
\cite{Savage}). The geometric categorification that we propose is a natural
next step in the geometric description of tensor products of the form~\eqref{tensor} of
$U_q(\mathfrak{sl}_2)$-modules via Steinberg varieties, which was started in
\cite{Savage}. We also remark that the generalised Steinberg varieties
naturally appear as characteristic varieties of Harish-Chandra bimodules
(\cite{Borho}). These facts strongly indicate that the geometric
categorification (via Borel-Moore homology), its relation to the algebraic categorification (by means of
Harish-Chandra bimodules), and to the theory of characteristic varieties is a
very rich area for future research. \\

Finally, we would like to mention that the categorification of the
representation theory of $U_q(\mathfrak{sl}_2)$ has powerful applications to
different areas in mathematics and physics. In particular, the general notion
of an $\mathfrak{sl}_2$-categorification via abelian categories has been
introduced in \cite{CRsl2} and effectively used to solve various outstanding
problems in the representation theory of finite groups. Also, in the same way
as the representation theory of $U_q(\mathfrak{sl}_2)$, especially in the
framework of the tensor products~\eqref{tensor}, is applied to invariants of
knots in three dimensions, the categorification of all these structures is
believed to produce invariants of $2$-knots in four dimensions. So far,
invariants of link cobordisms were obtained via a categorification of the
Jones polynomial in \cite{Jacobsson} and \cite{Khovanovcob}. Their representation theoretic interpretation requires further
steps in the ``categorification program'' of the representation theory of the
quantum group $U_q(\mathfrak{sl}_2)$.    

\subsection*{Acknowledgements}
We would like to thank Henning Haahr Andersen, Daniel Krashen, Wolfgang
Soergel and  Joshua Sussan for interesting discussions. We are in
particular indebted to Volodymyr Mazorchuk for sharing his ideas, several
discussions and comments on a preliminary version of this paper. The first
author was supported by the NSF grants DMS-0070551 and DMS-0457444. The second
author acknowledges NSF support in form of grant DMS-0407784. The third author
was supported by CAALT and EPSRC. Part of this research was carried out when she
visited Yale university in 2004. She would like to thank all the members of
the department of Mathematics for their hospitality. 

\section{Quantum $\mf{sl}_2$ and its finite dimensional representations} 
\label{quantum-sl2} 
\subsection{Definitions and preliminaries}
  We start by recalling some basics about the quantised enveloping algebra of
  $\SL_2$. For details we refer for example to \cite{CP}, \cite{Jquant}. 
  Let $\Q(q)$ be the field of rational functions with complex coefficients in an indeterminate 
  $q$. 

  \begin{definition}
  \label{defquantumsl2}
  The quantum group $U_q(\mf{sl}_2)$ is an 
  associative algebra over $\Q(q)$ with generators $E,F,K,K^{-1}$ 
  and relations 
  \begin{align*} 
  KE&= q^2 EK,& KF&= q^{-2}FK,& KK^{-1}&=1=K^{-1}K, \\
  &&EF-FE&= \frac{K-K^{-1}}{q-q^{-1}}. 
  \end{align*} 
  \end{definition} 
We will also denote $U_q(\mf{sl}_2)$ by $U$. For $a\in\mZ$ denote $E^{(a)}= \frac{E^{a}}{[a]!}$ and $F^{(a)}= \frac{F^{(a)}}{[a]!}$ 
where $[a]!= [a][a-1]\dots [1]$ and $[a]=\frac{q^a -q^{-a}}{q-q^{-1}}$. 
Let $[a,b]=\frac{[a]!}{[b]![a-b]!}$ for $0\le b \le a$.
  Let $\overline{\phantom{a}}$ be the $\Q$-linear involution of $\Q(q)$ 
  which changes $q$ into $q^{-1}$. A $\Q$-linear (anti-)automorphism $\phi$ of
  $U$ is called {\it $\Q(q)$-antilinear} if $\phi(fx)= \overline{f}\phi(x)$,
  for $f\in \Q(q)$, $x\in U$.
 
The algebra $U$ has an antilinear anti-automorphism $\tau$ given by 
  \begin{align}
  \label{tau}  
  \tau(E)&=q  F K^{-1},&
  \tau(F)&= q E K,& 
  \tau(K)&=K^{-1}. 
  \end{align} 
   
Let $\sigma: U\to U$ be the 
$\Q(q)$-linear algebra involution defined by  
\begin{align}
\label{sigma} 
\sigma(E)&= F,& \sigma(F)&= E,& \sigma(K)&= K^{-1}. 
\end{align}

Let $\psi: U\to U$ be the $\Q(q)$-antilinear algebra involution
 \begin{align*} 
 \psi(E)&= E,&\psi(F)&= F,&\psi(K)&= K^{-1}. 
 \end{align*} 

The algebra $U$ is also a Hopf algebra with the comultiplication
$\blacktriangle$ given by 
\begin{equation}
\label{comult}
\begin{split}
\blacktriangle(E)=1\otimes E + E\otimes K^{-1},\quad&\quad\quad\blacktriangle(F)=K\otimes
F+F\otimes 1,\\
\blacktriangle(K^{\mp 1})&=K^{\mp1}\otimes K^{\mp1}.
\end{split}
\end{equation}
(The antipode $S$ is defined as $S(K)=K^{-1}$, $S(E)=-EK$ and $S(F)=-K^{-1}F$.) 

\subsection{The $n+1$-dimensional representation $V_n$}
\label{finitedim}
For any positive integer $n$, the algebra $U$ has a unique (up to isomorphism) irreducible 
$(n+1)$-dimensional representation $V_n$ 
on which $K$ acts semi-simply with powers of $q$ as eigenvalues. The $V_n$
constitute a complete set of representatives for the iso-classes of simple
$U$-modules of type I. We can choose in $V_n$ a basis of 
weight vectors $\{ v_0, v_1, \dots , v_n \}$ 
such that $U$ acts as follows 
 \begin{align} 
\label{canbasis}
 K^{\pm 1} v_k
&= q^{\pm (2k-n)} v_k,&
&Ev_k = [k+1] v_{k+1},&
&Fv_k = [n-k+1] v_{k-1}.  
  \end{align} 
We call this basis {\it the canonical basis of} $V_n,$ 
since it is a special, though rather trivial, 
case of the Lusztig-Kashiwara canonical bases 
in finite-dimensional irreducible modules over quantum groups (see
\cite{Lu1}, \cite{Kas}). Let $\cU(\mathfrak{sl}_2)$ denote the universal
enveloping algebra of the semisimple complex Lie algebra $\mathfrak{sl}_2$. We denote by $\ov{V}_n$ the $n+1$-dimensional irreducible representation of $\mathfrak{sl}_2$.  

Given a $\Qq$-vector space $V$, a $\Q$-bilinear form $V\times V\to \Q(q)$ 
is called {\it semi-linear} if it is 
$\Q(q)$-antilinear in the first variable and $\Q(q)$-linear in the 
second, \ie
\begin{align} 
<fx,y>&= \overline{f}<x,y>, &<x,fy>&= f<x,y>&f\in \Qq, x,y\in U. 
\end{align} 
On $V_n$, there is a (unique up to scaling) nondegenerate semi-linear form 
\begin{equation}
<,>: V_n \times V_n \to \Q(q) 
\end{equation}
which satisfies 
$<xu,v> = <u,\tau(x)v>$ for any $x\in U_q(\mf{sl}_2)$ and 
$u,v\in V_n$.  
In the basis $\{ v_k\}_{0\leq k\leq n}$ the form is given by 
\begin{equation} 
\label{semilin} 
<v_k, v_l> = \delta_{k,l} q^{k(n-k)}[n,k] 
\end{equation} 
Define the {\it dual canonical basis} $\{ v^k\}_{0\leq k\leq n}$ of $V_n$ by 
$<v_l,v^k> = \delta_{k,l}q^{k(n-k)}$. Then $v_k = [n,k]v^k$ and 
the action of $E$, $F$ and $K$ in the dual canonical basis is 
\begin{align}
\label{dualcan} 
K^{\pm 1} v^k&=q^{\pm (2k-n)}v^k,&  
E v^k&=[n-k] v^{k+1},& 
F v^k&=[k] v^{k-1}. 
\end{align} 
The involutions $\psi$ and $\sigma$ of $U$ give rise to endomorphisms of $V_n$
as follows: Let $\sigma_n: V_n \to V_n$ be the $\Q(q)$-linear map defined by 
$\sigma_n(v_n)= v_0$ and $\sigma_n(x a) = \sigma(x)\sigma_n(a)$ for 
$x\in U$ and $a\in V_n$. Then 
\begin{eqnarray}
\label{omegan}
  \sigma_n(v_k) = v_{n-k} \text{ for any }k. 
\end{eqnarray}
Let  $\psi_n: V_n \to V_n$ be the $\Q$-linear map defined by 
$\psi_n(v_n) = v_n$ and $\psi_n(xa) = \psi(x) \psi_n(a)$ for 
$x\in U$ and $a\in V_n$. Then 
$\psi_n$ is $\Q(q)$-antilinear 
and 
\begin{eqnarray}
\label{psin}
  \psi_n(v_k) = v_k \text{ for any } k.
\end{eqnarray}

\subsection{Tensor products of finite dimensional representations}
Given a positive integer $n$ and a composition ${\bf d}=(d_1,d_2,\ldots,
d_r)$ of $n$ we use the comultiplication \eqref{comult} to define the $U$-module 
\begin{equation*}
V_{\bf d}=V_{d_1}\otimes V_{d_{2}}\otimes\cdots\otimes
V_{d_r}.
\end{equation*}
Let $\ov{V}_{\bf d}=\ov{V}_{d_1}\otimes
\ov{V}_{d_{2}}\otimes\cdots\otimes \ov{V}_{d_r}$ be the corresponding $\mathfrak{sl}_2$-module. 
The {\it standard basis} (and {\it dual standard basis}) of $V_{\bf{d}}$ is
given by $\{v_{\bf a}=v_{a_1}\otimes\cdots\otimes v_{a_r}\}$, 
($\{v^{\bf a}=v^{a_1}\otimes\cdots\otimes v^{a_r}\}$ respectively), where
${\bf a}$ runs through all sequences ${\bf a}=(a_1,a_2,\ldots, a_r)\in\mZ^r$
such that $0\leq a_j\leq d_j$ for $1\leq j\leq r$. Both bases are
$\mC(q)$-bases for $V_{\bf d}$. 
Likewise, we get the standard and the dual standard $\mQ$-bases for $\ov
V_{\bf d}$. We denote by $S_{\bf d}=S_{d_1}\times S_{d_2}\times
\cdots\times S_{d_r}$ the Young subgroup, corresponding to ${\bf d}$, of the symmetric group $S_n$.

\section{The Grothendieck group of $\cO$ and of the category of Harish-Chandra bimodules} 
\label{catO}
\subsection{Preliminaries}
For an abelian category $\mc{B}$ let ${\bf
  G}(\mc{B})=\mQ\otimes_\mZ[\mc{B}]$, where $[\mc{B}]$ denotes the Grothendieck 
group of $\mc{B}$. The latter is the (free) abelian group generated by symbols $[M]$, as 
$M$ ranges over all objects of $\mc{B}$, subject to relations 
$[M_2]=[M_1]+ [M_3]$ for all short exact sequences 
\begin{equation} 
0 \to M_1\to M_2\to M_3\to 0
\end{equation} 
in $\mc{B}$.
An exact functor $F:\mc{B}\rightarrow\mc{B}'$ between abelian
categories $\mc{B}$ and $\mc{B'}$ induces a $\mC$-linear map $F^{\bf G}: {\bf  G}(\mc{B})\rightarrow {\bf  G}(\mc{B}')$.\\   

For a complex Lie algebra $\mg$ we denote by $\cU(\mg)$ its universal
enveloping algebra with centre $\cZ(\mg)$. We fix an integer $n\geq 2$ and set
$\mg=\mathfrak{gl}_n$. We fix a triangular decomposition
$\mg=\mathfrak{n}_-\oplus \mh\oplus\mathfrak{n}_+$ of $\mg$.  
Let $\cO=\cO(\GL_n)$ denote the corresponding Bernstein-Gelfand-Gelfand category
for $\mg=\mathfrak{gl}_{n}$ (see \cite{BGG}), \ie the category of finitely
generated $\mg$-modules which are $\mh$-diagonalisable and locally
$\cU(\mathfrak{n}_+)$-nilpotent. Let $W=S_{n}$ denote the Weyl group generated by the
simple reflections $s_i$, $1\leq i\leq n-1$, where $s_is_j=s_js_i$ if
$|i-j|>1$.    
Let $\rho$ be the half-sum of positive roots. For $w\in W$ and
$\la\in\mh^\ast$ let $w\cdot\la=w(\la+\rho)-\rho$. Let $W_\la=\{w\in W\mid
w\cdot\la=\la\}$ be the stabiliser of $\la$ with respect to this action. We
denote by $W^\la$ the set of (with respect to the length function) shortest
coset representatives in $W/W_\la$.  We denote by $w_0$ the longest element in
$W$ and by $w_0^\la$ the longest element in $W_\la$. A weight $\la\in\mh^\ast$
is {\it integral}, if $\langle \la,\check{\alpha}\rangle\in\mZ$ for any coroot
$\check{\alpha}$. We call a weight $\lambda\in\mh^\ast$ {\it (strictly)
  dominant} if $\langle\la +\rho,\check{\alpha}\rangle\geq 0$ (or $\langle
\la,\check{\alpha}\rangle\geq 0$ resp.) for any coroot $\check{\alpha}$
corresponding to a positive root $\alpha$ such that $\langle \la,\check{\alpha}\rangle\in\mZ$. Note that this terminology is not commonly used in the literature, usually a
 weight is called dominant and integral if it is strictly dominant and
 integral in our terminology.
 With this notion, the action of the centre of $\cU(\mathfrak{gl}_n)$ gives a
 block decomposition $\cO=\oplus\;{_\la}\cO$, where the sum runs through the set of
 dominant weights $\la$. The finite dimensional simple
 objects in $\cO$, however, are naturally indexed by strictly dominant and
 integral 
 weights. 
 For $\la\in\mh^\ast$ let $M(\la)=\cU(\mg)\otimes_{\cU(\mb)}\mC_{\la}$ denote the Verma module
with highest weight $\la$. Let $P(\la)$ be its projective cover
with simple head $L(\la)$. 
If $\la\in\mh^\ast$ is dominant and integral, then $_\la\cO$ denotes the block containing
all $M(x\cdot\la)$ with $x\in W$. (Based on the functor~\eqref{BGGfunc} we prefer the notation ${}_\la\cO$ to the
more common notation $\cO_\la$ in order to be consistent with the fact that the objects
of $\cO$ are {\it left} $\mg$-modules.) We denote by $\op{d}$ the usual contravariant
duality on $\cO$ preserving the simple objects. Recall that a module in
${}_\la\cO$ having a Verma flag and a dual Verma flag is called a {\it tilting
  module}. Let $T(x\cdot\la)\in{}_\la\cO$ be the indecomposable tilting module
with $M(x\cdot\la)$ occurring as a submodule in any Verma flag. (For the
classification we refer to \cite{CI}, for the general theory to \cite{DR1},
for example.) 

\subsection{Tensor products of finite dimensional $\cU(\SL_2)$-modules and
  Harish-Chandra bimodules}
The purpose of this section is to associate with any finite tensor product
$\ov{V}_{\bf d}$ of
finite dimensional  $\cU(\SL_2)$-modules an abelian category $\mc{B}$ of Harish-Chandra
bimodules together with an isomorphism $\ov\Phi:{\bf G}(\mc{B})\cong
\ov{V}_{\bf d}$ of $\mQ$-vector
spaces. In the following section we will define exact endofunctors on $\mc{B}$
which give rise to a $\cU(\SL_2)$-module
structure on ${\bf G}(\mc{B})$ and show that the morphisms $\ov\Phi$ become
isomorphisms of $\cU(\SL_2)$-modules.   

We choose an ONB $\{e_i\}_{1\leq i\leq n}$ of
  $\mR^{n}$ and identify $\mC\otimes_\mR{\mR^{n}}$ with $\mh^\ast$ such that
  $R_+=\{e_i-e_j\mid i<j\}$ is the
  set of positive roots. The simple reflection $s_i\in W$ acts by permuting
  $e_i$ and $e_{i+1}$. For ${\bf a}=(a_1,a_2,\ldots,a_n)\in\mathbb{Z}^n$ let
  $M({\bf a})$ also denote the Verma module with highest weight
  $\sum_{i=1}^{n}a_ie_i-\rho$. 
  If $I=\{i_1,i_2\cdots i_r\}\subseteq \{1,\ldots n\}$ then we denote $W_I=\langle
s_i\mid i\notin I\rangle$. Corresponding to $I$ we fix an integral block
$_{I,n}\cO={_{\la_I}}\cO$, such that $\{w\in W\mid w\cdot\la_I=\la_I\}=W_I$ and
$\la_I$ is minimal with this property. In the following we will also just write
$_{i_1,i_2,\cdots, i_r;n}\cO={}_{I;n}\cO$, $W_I=W_{i_1,i_2,\cdots, i_r}$, $W^I=W^{i_1,i_2,\cdots, i_r}$ etc. In particular, we have the ``maximal singular'' blocks ${_{\omega_i}}\cO={_{i;n}}\cO$ where
$0\leq i\leq n$, and $\omega_i$ is the $i$-th fundamental weight. (Note
that by definition $W_\emptyset=W_{n}=W$. We remark also that ${}_{1;2}\cO$ is in fact
regular.) Let ${_{i;n}}\cO$ be the category containing only the zero module for
$n<i$ or $i<0$. Note that $M({\bf a})\in{_{i;n}}\cO$ if and only if
  $a_j\in\{0,1\}$ for all $j$ and ${\bf a}$ contains exactly $i$ ones. For
  $x\in W/W_I$ we denote by ${\bf a}(x)={\bf a}_I(x)\in\mZ^n$ the sequence such that 
  \begin{eqnarray}
    \label{bij}
  M(x\cdot\la_I)=M({\bf a}(x)).
  \end{eqnarray}
By definition, the vector space $\ov V_1^{\otimes n}$ has a basis of the form
$\{v_{\bf a}\}$, where ${\bf a}$ runs through all $\{0,1\}$-sequences of
length $n$. On the other hand ${\bf G}({_{i;n}}\cO)$ has a basis of the form
$\{[M({\bf a})]\}$, where  ${\bf a}$ runs through all $\{0,1\}$-sequences of
length $n$, containing exactly $i$ ones. Therefore, (see also
\cite[(34)]{BFK}) there is an isomorphism of vector spaces
  \begin{eqnarray}
\label{phi}
  {\bf{G}}\big(\bigoplus_{i=0}^{n}{_{i;n}}\cO\big)&\cong&\overline{V}_1^{\otimes n}\\
 1\otimes\big[M({\bf a})\big]&\mapsto&v_{a_1}\otimes\cdots\otimes v_{a_n}.\nonumber
  \end{eqnarray}
Before generalising this to arbitrary tensor products we give an 
\begin{exs}
{\rm
  If $n=2$ then ${_{0;2}}\cO$ and ${_{2;2}}\cO$ are semisimple with one simple
  object $L((0,0))=M((0,0))=P((0,0))=T((0,0))$ (or $L((1,1))=M((1,1))=P((1,1))=T((1,1))$
  respectively), whereas ${_{1;2}}\cO$ has the two simple objects $L((1,0))$ and
  $L((0,1))=M((0,1))=T((0,1))$ and $[M((1,0))]=[P((1,0))]=[L((1,0))]+[L((0,1))]$ and
  $[P((0,1))]=[M((1,0))]+[M((0,1))]=[L((0,1))]+[L((1,0))]+[L((0,1))]$. Hence
  ${\bf{G}}\big(\bigoplus_{i=0}^{2}{_{i;2}}\cO\big)\cong\overline{V}_1^{\otimes
  2}$ as vector spaces.} 
\end{exs}

Let $\cH=\cH(\mg)$ denote the category of {\it Harish-Chandra bimodules} for
$\mg$. That is the full subcategory, inside the category of finitely generated
$\cU(\mg)$-bimodules of finite length, given by all objects which are locally finite with
respect to the adjoint action of $\mg$ (see \eg \cite{BG}, or for an overview 
\cite[Kapitel 6]{Ja2}). As for the category $\cO$, the action of $\cZ(\mg)$ gives a block decomposition
$\cH=\oplus_\la\cH_\mu$, where $\la$, $\mu\in\mh^\ast$ are dominant
weights. More precisely it is given as follows: Let $\la\in\mh^\ast$ be
dominant. We denote by $\KER\chi_\la$ the $\cZ(\mg)$-annihilator
of the Verma module $M(\la)$. Note that $\KER\chi_\la$ is a maximal ideal in
$\cZ(\mg)$. A Harish-Chandra bimodule $X$ is an object of $_\la\cH_\mu$ if and only if
$(\KER\chi_\la)^mX=0=X(\KER\chi_\mu)^m$ for large enough
$m\in\mZ_{>0}$.

For any two $\mg$-modules $M$ and $N$, the space $\HOM_\mC(M,N)$ is naturally
a $\cU(\mg)$-bimodule. Let $\cL(M,N)$ denote its maximal submodule which is
locally finite with respect to the adjoint action of $\mg$. Then the simple objects in
$_\la\cH_\mu$ are of the form $\cL(M(\mu), L(x\cdot\la))$ where $x$ is a
longest coset representative in $W_\mu\backslash W / W_\la$
(see \eg \cite[6.26]{Ja2}). However, $_\la\cH_\mu$ does not have enough
projective objects. Therefore, we consider the category $_\la\cH_\mu^1$ which
is by definition the full subcategory of $_\la\cH_\mu$ given by all
objects such that $X\KER\chi_\mu=0$. (Note that the simple objects stay the
same.) Let $\la\in\mh^\ast$ be integral and dominant. If $\mu$ is integral, regular and
dominant then the functor $F_\mu=\bullet\;\otimes_{\cU(\mg)}M(\mu)$ defines an
equivalence of categories 
\begin{eqnarray}
\label{BGequiv}
{}_\la F_\mu:\quad ^{}{}_\la\cH_\mu^1\cong {_\la}\cO.  
\end{eqnarray}
The inverse functor is given by $\cL(M(\mu),\bullet\;)$. If $\mu$ is singular,
then the functor ${}_\la F_\mu$ defines only an embedding. (All this is proved in
\cite{BG}, for an overview see also \cite{Ja2}). 

In this setup, one of the main ideas of this paper is that formula~\eqref{phi} is
in fact only a very special case of the following more general fact: 
  
\begin{prop}
\label{lemmatensor}
  Let ${\bf d}=(d_1,\ldots d_r)$ be a composition of $n$. Let $\mu\in\mh^\ast$
  be dominant and integral such that $W_\mu=S_{\bf d}$. There is an
  isomorphism of vector spaces
\begin{eqnarray*}
    \ov\Phi:\quad\quad{\bf
      G}\big(\bigoplus_{i=0}^n\;{_{\omega_i}}\cH_\mu^1(\mathfrak{gl}_n)\big)&\cong&\ov V_{\bf
      d},\\
1\otimes\big[\cL(M(\mu),M({\bf a}))\big]&\mapsto& v^{{\bf a}(\mu)}=v^{a(\mu)_1}\otimes v^{a(\mu)_2}\otimes\cdots\otimes v^{a(\mu)_r}, 
  \end{eqnarray*}
where $a(\mu)_j=|\{a_k=1\mid d_{j-1}< k\leq d_j\}|$ with $d_0=0$, $d_{r+1}=n$.
\end{prop}

Note that Proposition~\ref{lemmatensor} is in fact a generalisation of \eqref{phi}, because if ${\bf d}=(1,1,\cdots 1)$, then $_{\omega_i}\cH_{\bf
  d}^1(\GL_n)\cong{}_{i;n}\cO$ via the equivalence \eqref{BGequiv} and the
  isomorphism $\ov\Phi$ gives rise to the one from \eqref{phi}, because
  $v_k=v^k$ in $V_1$. The ``asymmetry'' with respect to the central characters associated to 
  $\bigoplus_{i=0}^n\;{_{\omega_i}}\cH_\mu^1(\GL_n)$ (\ie we consider bimodules with
  a fixed central character from the right hand side, but with a fixed {\it
  generalised} character from the left hand side, appears to be unnatural. In
  fact, we could also work with the category
  $\bigoplus_{i=0}^n\;{_{\omega_i}}\cH_\mu(\GL_n)$ instead, since its Grothendieck
  group coincides with the one of the previous category and the functors $\cE$
  and $\cF$ (which will be introduced later) can be extended naturally. However, proofs become much simpler if
  we make use of the functor $F_\mu$. Moreover, the interpretation of the
  canonical basis in terms of tilting objects (see Theorem~\ref{bases2} and Remark~\ref{finalremark}) also militate in favour of
  using the category $\bigoplus_{i=0}^n\;{_{\omega_i}}\cH_\mu^1(\GL_n)$.  
Before we prove the Proposition~\ref{lemmatensor} we give an
 
\begin{exs}
\label{Ex}
{\rm 
  We consider the Lie algebra $\mathfrak{gl}_3$. The block $_{0;3}\cO$ is
  semisimple with the only (simple and projective) Verma module
  $M((0,0,0))$. Likewise, the block $_{3;3}\cO$ is
  semisimple with the only (simple and projective) Verma module
  $M((1,1,1))$.  Each of the blocks $_{1;3}\cO$ and  $_{2;3}\cO$ contains
  exactly three Verma modules, namely $M((1,0,0))$, $M((0,1,0))$ and $M((0,0,1))$
  (or the Verma modules $M((1,1,0))$, $M((1,0,1))$ and $M((0,1,1))$ respectively). We have the
  following possibilities
  \begin{itemize}
  \item $W_\mu$ is trivial, \ie $\mu$ is regular: the simple objects in
  $\bigoplus_{i=0}^n\;{_{\omega_i}}\cH_\mu^1$ are exactly the 
  $\cL(M(\mu), L({\bf a}))$, where $L({\bf a})$ occurs is the head of one out
  of these eight Verma modules. This ``models'' the $8$-dimensional vector
  space $V_1^{\otimes 3}$.
  \item $W_\mu\cong S_1\times S_2$: the simple objects in
  $\bigoplus_{i=0}^n\;{_{\omega_i}}\cH_\mu^1$ are the  $\cL(M(\mu), L({\bf a}))$, where ${\bf
  a}\in\{(0,0,0),(1,0,0),(0,0,1),(1,0,1), (0,1,1),(1,1,1)\}$. This
  ``models'' the $6$-dimensional vector space $V_1\otimes V_2$. 
  \item $W_\mu\cong S_2\times S_1$: the simple objects in
  $\bigoplus_{i=0}^n\;{_{\omega_i}}\cH_\mu^1$ are the  
  $\cL(M(\mu), L({\bf a}))$, where ${\bf
  a}\in\{(0,0,0),(0,1,0),(0,0,1),(0,1,1),(1,1,0),(1,1,1)\}$. This ``models''
  the $6$-dimensional vector space $V_2\otimes V_1$.
   \item $W_\mu\cong S_3$: the simple objects in
  $\bigoplus_{i=0}^n\;{_{\omega_i}}\cH_\mu^1$ are the 
  $\cL(M(\mu), L({\bf a}))$, where ${\bf
  a}\in\{(0,0,0),(0,0,1),(0,1,1),(1,1,1)\}$. This ``models''
  the $4$-dimensional vector space $V_3$.     
  \end{itemize}
Note that  the simple objects in $\bigoplus_{i=0}^n\;{_{\omega_i}}\cH_\mu^1$ are of the form
  $\cL(M(\mu), L({\bf a}))$, where ${\bf
  a}=(a_1,a_2,a_3)$ such that the $a_i$'s are (weakly) increasing within the parts of
  the composition given by $\mu$. 
}
\end{exs}

\begin{proof}[Proof of Proposition~\ref{lemmatensor}]
  Recall that the simple objects in $_{\omega_i}\cH_{\bf d}^1(\GL_n)$
  are of the form
  $\cL(M(\mu), L(x\cdot\omega_i))$ where $x$ runs
  through the set of longest coset representatives $D$ for the double cosets
  ${W_\mu}\backslash W/W_i$.  On the other hand $x\in D$
  if and only if $x\cdot\omega_i=\sum_{j=1}^{n}a_je_j-\rho$, where $a_k\leq a_j$ if
  $d_{l-1}<k<j\leq d_l$ for some $l\in\{1,2,\ldots, r\}$ and the number of $j$'s
  such that $a_j=1$ is $i$. That is, the simple objects in
  $\bigoplus_{i=0}^n\;{_{\omega_i}}\cH_\mu^1(\GL_n)$ are exactly the bimodules
  of the form
  $\cL(M(\mu), L({\bf a}))$, where ${\bf
  a}=(a_1,a_2,\ldots, a_n)$ is a $\{0,1\}$-sequence with exactly $i$ ones such that the $a_i$'s are (weakly) increasing
  within the parts of the composition given by $\mu$.  Let now $y\in W_\mu$ and $x\in W_\mu\backslash W / W_i$
  be a longest coset representative. The following formula holds (see
  \cite[Lemma 2.5]{Jo1}):
  \begin{eqnarray}
\label{Joseph1}
    \cL(M(\mu),M(x\cdot\omega_i))&=&\cL(M(y\cdot\mu),M(x\cdot\omega_i))\\
    &\cong&\cL(M(\mu),M(y^{-1}x\cdot\omega_i)).\nonumber
  \end{eqnarray}

Therefore, the isomorphism classes
  $[\cL(M(\mu), M(x\cdot\omega_i))]$ $(x\in D)$ give rise to a basis
  of $\bf{G}\big({}_{\omega_i}\cH_{\bf d}^1(\GL_n)\big)$. The statement of the
  proposition follows.
\end{proof}
In the situation of Proposition~\ref{lemmatensor} let $x\in
  W_\mu\backslash W/W_i$ be a longest coset representative. Set
  $[\Delta({\bf a}(x))]=\sum [M]$, where $[M]\in\{[\cL(M(\mu),M(yx\cdot\la)]\mid y\in
  W_\mu\}$. (Note the difference between the symbols $\Delta$ and
  $\blacktriangle$, the latter denoting the comultiplication.) From the
  formula~\eqref{Joseph1} it follows that the $1\otimes[\Delta({\bf a}(x))]$, where
  $x$ runs through all longest coset representatives in
  $W_\mu\backslash W/W_i$ form a basis of ${\bf G}\big({}_{\omega_i}\cH_{\bf
  d}^1(\GL_n)\big)$ (although the $[\Delta({\bf a}(x))]$ do not form a basis
  of $[{}_{\omega_i}\cH_{\bf d}(\GL_n)\big]$ in general, see \eg Remark~\ref{properly}~\eqref{properly2}).\\

The following holds 
\begin{corollary}
\label{standards}
With the assumption of Proposition~\ref{lemmatensor} we have
\begin{eqnarray*}
    \ov\Phi\big(1\otimes[\Delta({\bf a})]\big)=v_{{\bf a}(\mu)}:=v_{a(\mu)_1}\otimes v_{a(\mu)_2}\otimes\cdots\otimes v_{a(\mu)_r}, 
  \end{eqnarray*}
where $a(\mu)_j=|\{a_k=1\mid d_{j-1}< k\leq d_j\}|$ with $d_0=0$, $d_{r+1}=n$.
\end{corollary}
\begin{proof}
The statement follows immediately from the definitions of $\ov\Phi$ and $[\Delta({\bf
  a})]$ together with the formula~\eqref{Joseph1} and the equality $v_k=[n,k]v^k$ in $V_n$. 
\end{proof}

\begin{remark}
\label{properly}
  {\rm  
    \begin{enumerate}
    \item 
\label{properly1}
In general, the category
$\bigoplus_{i=0}^n\;{_{\omega_i}}\cH_\mu^1$ is not necessarily a highest weight category,
in the sense of \cite{CPS}. However (see \cite{KoMa}), it is equivalent to a module category over a {\it properly
  stratified algebra} (as introduced in \cite{Dl}, generalising the notion of
quasi-hereditary algebras and highest weight categories). In particular, this algebra
might have infinite global dimension. If $A$ is a properly stratified algebra
 then, by definition, the projective $A$-modules have a filtration such
 that subquotients are so-called {\it standard modules} $\Delta({\bf a})$; any
 standard module $\Delta(\bf{a})$ has a filtration with subquotients, each
 isomorphic to the (same) {\it proper costandard module} $\ov\Delta{(\bf
   a)}$. If the standard modules coincide with proper standard modules, then
 $A$ is quasi-hereditary. For example, for any block of category $\cO$, the
 (proper) standard modules are given by the Verma modules. In $_\la\cH^1_\mu$, the
 proper standard modules are given by bimodules of the form
 $\cL(M(\mu), M(\bf{a}))$, whereas the standard modules are certain modules
 $\Delta({\bf a})$, where its image in the Grothendieck group $[\Delta({\bf
a})]$ is as above (for definition and general theory see \eg \cite{Dl}, for
the special situation see \eg \cite{MS}). 
\item \label{properly2}
As an example let us consider the case $\mg=\GL_2$. The category $_0\cH_{-\rho}^1$
  has one indecomposable projective object, $P=\cL(M(-\rho), P(-2\rho))$, and
  one simple object, 
  $S=\cL(M(-\rho),M(-2\rho))$. Then $P$ is the unique standard module, whereas
  $S$ is the unique proper standard module and $P$ is a self-extension of
  $S$. The category $_0\cH_{-\rho}^1$ is equivalent to $\mC[x]/(x^2)\MOD$ (by
  \cite[Endomorphismensatz]{Sperv}),
  in particular, it has infinite global dimension. Note that the functor
  ${}_0F_{-\rho}$ maps $S$ to the {\it dual} Verma module $\op{d}M(0)$ with highest
  weight zero.  
\end{enumerate}
}
\end{remark}

\subsection{The submodule $V_n$ inside $V_1^{\otimes n}$}
There is a unique direct summand isomorphic to $V_n$ inside
$V_1^{\otimes n}$. The inclusion is given by the map 
\begin{eqnarray}
\label{inclusion}
 i_n: V_n&\hookrightarrow& V_1^{\otimes n},\nonumber\\
v_k&\mapsto& \sum_{|{\bf a}|=k} q^{{\bf a}^-}v_{\bf a}, 
\end{eqnarray}
  where $|{\bf a}|=\sum_{i=n}^r a_i$ and $\bf{a}^-$ is the cardinality of
  the set $\{(i,j)\mid 1\leq i<j\leq r, a_i>a_j\}$. A split of this inclusion
  is given by the projection 
  \begin{eqnarray}
\label{pi}
  \pi_n: V_1^{\otimes n}&\rightarrow& V_n,\nonumber\\
v^{\bf a}&\mapsto& q^{-{\bf a}^+}v^{|\bf{a}|}    
  \end{eqnarray}
 where  $\bf{a}^+$ is the cardinality of
  the set $\{(i,j)\mid 1\leq i<j\leq n, a_i<a_j\}$ (see \eg \cite[Proposition 1.3]{FK}). The
  composition $i_n\circ\pi_n$ is the {\it Jones-Wenzl projector}. \\
We define a linear map $F: {\bf
 G}\big(\bigoplus_{i=0}^n\;{_{\omega_i}}\cH_\mu^1(\GL_n))\rightarrow {\bf
      G}\big(\bigoplus_{i=0}^n\;{_{\omega_i}\cO})$ as follows: For $x\in {W_\mu}
  \backslash W/W_i$ a longest coset representative let $F(1\otimes[\Delta({\bf
  a}(x))])=\sum1\otimes [M]$,
  where the sum runs over all $[M]\in\{[M(yx\cdot\omega_i)]\mid y\in W_\mu\}$. Then the following holds: 
\begin{prop}
\label{JonesWenzl}
 Let ${\bf d}=(d_1,\ldots d_r)$ be a composition of $n$. Let $\mu\in\mh^\ast$
  be dominant and integral such that $W_\mu=S_{\bf d}$. The following diagrams commute
  \begin{eqnarray*}
\xymatrix{{\bf
 G}\big(\bigoplus_{i=0}^n\;{_{\omega_i}}\cH_\mu^1(\GL_n))\ar@{->}[d]_{F}\ar@{->}[rr]^{\ov\Phi}\ar@{->}[d]&&\ov
 V_{\bf d},\ar@{->}[d]^{i_{d_1}\otimes\cdots\otimes i_{d_r}}\\
 {\bf
      G}\big(\bigoplus_{i=0}^n\;{_{\omega_i}\cO})\ar@{->}[rr]^{\ov\Phi}&&\ov
      V_1^{\otimes{d_1}}\otimes\cdots\otimes
      \ov V_1^{\otimes{d_r}}=\ov V_1^{\otimes n}
}
  \end{eqnarray*}
  \begin{eqnarray*}
\xymatrix{{\bf
 G}\big(\bigoplus_{i=0}^n\;{_{\omega_i}}\cH_\mu^1(\GL_n))\ar@{->}[d]_{F}\ar@{->}[rr]^{\ov\Phi}\ar@{->}[d]&&\ov
 V_{\bf d},\ar@{<-}[d]^{\pi_{d_1}\otimes\cdots\otimes \pi_{d_r}}\\
 {\bf G}\big(\bigoplus_{i=0}^n\;{_{\omega_i}\cO})\ar@{->}[rr]^{\ov\Phi}&&\ov
      V_1^{\otimes{d_1}}\otimes\cdots\otimes
      \ov V_{1}^{\otimes{d_r}}=\ov V_1^{\otimes n}
}
  \end{eqnarray*}
\end{prop}

\begin{proof}
Let $x\in {W_\mu}
  \backslash W/W_i$ be a longest coset representative. By Corollary~\ref{standards} we have 
  $i \circ \ov\Phi([\Delta({\bf a}(x))])=i (v_{\bf
  a(\mu)})$. Now $i_{d_j}(v_{\bf a}(\mu)_j)= \sum_{|{\bf b(j)}|={\bf
  a}(\mu)_j}v_{\bf b(j)}$. This means
\begin{eqnarray*}  
i\circ\ov\Phi([\Delta({\bf a}(x))])=
  \bigotimes_{j=1}^r\big(\sum_{|{\bf b(j)}|={\bf a}(\mu)_j} v_{\bf
  b(j)}\big).
\end{eqnarray*}
 On the other hand we get $\ov\Phi\circ F([\Delta({\bf a}(x))])=\sum\ov\Phi([M])$,
  where the sum runs over all $[M]\in\{[M(yx\cdot\omega_i)]\mid y\in
  W_\mu\}$. This means
  \begin{eqnarray*}
  \ov\Phi\circ F([\Delta({\bf a}(x))])=\bigotimes_{j=1}^r\big(\sum_{|{\bf b(j)}|={\bf
  a}(\mu)_j} v_{\bf {b(j)}}\big).
\end{eqnarray*}
 Hence, the first diagram commutes. 
  Let now $\pi={\pi_{d_1}\otimes\cdots\otimes \pi_{d_r}}$. The first diagram
  says that $i\circ\ov\Phi=\ov\Phi\circ F$. We  get $\ov\Phi=\pi\circ
  i\circ\ov\Phi=\pi\circ\ov\Phi\circ F$, hence the second diagram commutes as
  well.
\end{proof}

  We would like to remark that the map $F$ has in fact a categorical interpretation. To explain this we
  have to recall some facts on the properly stratified structure of ${_\la}\cH_\mu^1$ for integral dominant weights $\la$ and
  $\mu$. The indecomposable projective objects in
  ${_\la}\cH_\mu^1$ are exactly the bimodules $X$ such that ${}_\la F_\mu
  (X)\cong P(x\cdot\la)\in {_\la}\cO$ for some longest coset representative
  $x\in W_\mu\backslash W/W_\la$ (see \cite{BG} or \cite[6.17, 6.18]{Ja2}).
  For a longest coset representative $x\in W_\mu\backslash W/W_\la$ let $P(x)$
  be the indecomposable projective object ${_\la}\cH_\mu^1$, such
  that ${}_\la F_\mu(P(x))=P(x\cdot\la)$.   Let
  $P^{<x}=\oplus P(y)$, where the sum runs over all longest coset
  representatives $y\in W_\mu\backslash W/W_\la$, and $y<x$ in the Bruhat
  ordering. Then the standard module corresponding to $x$ is defined as 
  $\Delta(x)=P(x)/M$, where $M$ is the trace of  $P^{<x}$ in
  $P(x)$ (see \cite[Definition 3]{Dlab}). One can show that in fact
  $[\Delta(x)]=[\Delta({\bf a}(x))]$ as defined above, and, moreover, the images of the standard modules from ${}_\la\cH_\mu^1$
  under the functor ${}_\la F_\mu$ are objects in ${_\la}\cO$ having a Verma flag (\ie a
  filtration with subquotients isomorphic to Verma modules). A proof of this
  fact can be found for example in \cite{MS}. More precisely we have that
  $[{}_\la F_\mu\Delta(x)]=\sum [M]$, where the sum runs over all $[M]\in\{[M(yx\cdot\la)]\mid y\in
  W_\mu\}$ (see \cite[Proposition 2.18]{MS}, note that the proof there works
  also for singular $\la$). So, $F$ should be considered as a replacement
  for ${}_\la F_\mu^{\bf G}$ (which is not defined in general, since ${}_\la
  F_\mu$ is not necessarily exact). The properties of the derived functor of
  ${}_\la F_\mu$ were studied by J. Sussan, who constructed
  a categorification of the Jones-Wenzl projector in the setting of derived
  categories (\cite{Josh}).   

\subsection{Categorification of finite tensor products of finite dimensional $\cU(\SL_2)$-modules via representations of $\GL_n$.}
\label{section1}

The purpose of this section is to lift the isomorphisms $\ov\Phi$ from
Proposition~\ref{lemmatensor} to isomorphisms of $\cU(\SL_2)$-modules. In
other words, we {\it categorify} tensor
products of finite dimensional $\cU(\SL_2)$-modules.
Our main results (Theorem~\ref{HCnongrad} and Theorem~\ref{cat}) are {\it
  categorifications} of the modules $\ov V_{\bf d}$ and $V_{\bf d}$. Let us
  make this more precisely: By a {\it categorification of $\ov V_{\bf d}$} we mean an abelian category $\cC$ together with exact endofunctors
  $\cE$ and $\cF$ and an isomorphism (of vector spaces) $\Psi:$ ${\bf G(\cC)}\cong\ov V_{\bf
      d}$ such that  
  \begin{eqnarray}
\label{sl2}
\Psi([\cE M])=E\Psi([M]) &\mbox{and}&
\Psi([\cF M])=F\Psi([M])    
  \end{eqnarray}
 for any object $M\in\cC$. That is, if $\cC$ is a categorification of
  $\ov V_{\bf d}$ then the functors $\cE$ and $\cF$ define an
  $\SL_2$-action on ${\bf G(\cC)}$ such that $\Psi$ becomes an isomorphism of
  $\SL_2$-modules.\\
 
For $\la$, $\mu\in\mh^\ast$ dominant and integral let
$\mathds{T}_\la^\mu:{}_\la\cO\rightarrow{}_\mu\cO$ denote the translation functor (see
\cite{BG}, \cite{Ja2}). With the notations above we also denote $\mathds{T}^J_I=\mathds{T}_{i_1,i_2,\ldots,
  i_r}^{j_1,j_2,\ldots, j_s}:{}_I\cO\rightarrow {}_J\cO$ for $J=\{j_1,j_2\ldots
  j_s\}\subseteq\{1,2,\ldots n\}$. 
Let $k\in\Z_{>0}$. Recall from \cite{BFK} the projective functors
$\cE_i^{(k)}:{}_{i;n}\cO(\mathfrak{gl}_n)\rightarrow{}_{i+k;n}\cO(\mathfrak{gl}_n)$ given by
tensoring with the $k$-th exterior power of the natural representation for
$\mathfrak{gl}_n$ composed with the 
projection onto $\cO(\mathfrak{gl}_n)_{{i+k};n}$. In particular,
$\cE_i=\cE_i^{(1)}$ is given by tensoring with the natural representation
composed with the projection onto $\cO(\mathfrak{gl}_n)_{{i+1};n}$. Let
$\cF_{i}^{(k)}:{}_{i;n}\cO\rightarrow{}_{i-k;n}\cO$ be the adjoint of $\cE_{i-k}^{(k)}:{}_{i-k;n}\cO\rightarrow{}_{i;n}\cO$. Set $\cE=\oplus_{i=0}^n\cE_i$, 
$\cF=\oplus_{i=0}^n\cF_i$. Since the functor $F_\mu$ commutes with
  tensoring with finite dimensional (left) $\mg$-modules, the functors $\cE$
  and $\cF$ give rise to endofunctors of
  $\bigoplus_{i=0}^n\;{_{\omega_i}}\cH_\mu^1$ via restriction. We will denote
  these functors with the same letter $\cE$ and $\cF$ respectively. \\

We get the following result (which can also be viewed as a specialisation of Theorem~\ref{cat}): 

\begin{theorem}
\label{HCnongrad}
Let ${\bf d}=(d_1,\ldots, d_r)$ be a composition of $n$. Let $\mu\in\mh^\ast$
  be dominant and integral such that $W_\mu=S_{\bf d}$. Then the isomorphism 
  \begin{eqnarray*}
       \overline\Phi:\quad\quad{\bf
      G}\big(\bigoplus_{i=0}^n\;{_{\omega_i}}\cH_\mu^1(\GL_n))&\cong&\ov V_{\bf d},\\
\big[\cL(M(\mu),M({\bf a}))\big]&\mapsto&v^{\bf b}=v^{b(\mu)_1}\otimes
      v^{b(\mu)_2}\otimes\cdots\otimes  v^{b(\mu)_r}, 
  \end{eqnarray*}
where $b_j(\mu)=|\{a_k=1\mid d_{j-1}< k\leq d_j\}|$ with $d_0=0$, $d_{r+1}=n$, together with the endofunctors $\cE$, $\cF$ of
  $\bigoplus_{i=0}^n\;{_{\omega_i}}\cH_\mu^1(\GL_n)$ give rise to a
  categorification of $\ov V_{\bf d}$.
\end{theorem}

\begin{proof} 
By Proposition~\ref{lemmatensor} and \cite[Theorem 1]{BFK}, we only have to check the formulas \eqref{sl2} for
  $\Psi=\ov\Phi$, that is  $\ov\Phi([\cE
  \bullet])=E\ov\Phi([\bullet])$,  and $\ov\Phi([\cF
  \bullet])=F\ov\Phi([\bullet])$. \\
Let now $M=\cL(M(\mu),M({\bf a}))$ for some
  $\{0,1\}$-sequence ${\bf a}$ of length $n$. Then we get 
\begin{equation*}
\ov\Phi([\cE
  M])=\ov\Phi([\cL(M(\mu),\cE M({\bf a}))])=\ov\Phi(\sum_{\bf a'}[\cL(M(\mu),
  M({\bf a'}))]), 
\end{equation*}
where the sum runs over all
  sequences ${\bf a'}$ built up from ${\bf a}$ by replacing one zero by a
  one (see \cite[(38)]{BFK}). From the definition of $\Phi$ and formula~\eqref{Joseph1} we get $\ov\Phi([\cE
  M])=\sum_{\bf a'}\alpha_{{\bf a},{\bf a'}} v^{{\bf a'}(\mu)}$, where the sum runs over all
  ${\bf a'}$ where $a'_{j_o}=b(\mu)_{j_o}+1$ for some $j_0$ and
  $a'_{j}=b(\mu)_j$ for $j\not=j_0$ and $\alpha_{{\bf a}, {\bf a'}}$ is the number of zeros
  occurring in $\{a_k\mid d_{j_0}-1<k\leq d_{j_0}\}$. On the other hand $E
  v^{\bf b}=\sum_{\bf a'}\beta_{{\bf a},{\bf a'}} v^{\bf a'}$, where the sum runs over all
  ${\bf a'}$ where $a'_{j_o}=b(\mu)_{j_o}+1$ for some $j_0$ and
  $a'_{j}=b(\mu)_j$ for $j\not=j_0$ and $\beta_{{\bf a},{\bf a'}}$ is defined by
  the equation $Ev^{b(\mu)_{j_0}}=\beta_{{\bf a},{\bf a'}}v^{b(\mu)_{j_0}+1}$  in
  $V_{d_{j_0}-{d_{j_0-1}+1}}$. By formula \eqref{dualcan} we actually have 
  $\beta_{{\bf a},{\bf a'}}={d_{j_0}}-{d_{j_0-1}+1-{b(\mu)_{j_0}}}$. The latter
  is, by definition of ${b(\mu)_{j_0}}$, exactly the number of  zeros
  occurring in $\{a_k\mid d_{j_0}-1<k\leq d_{j_0}\}$. We therefore get that $\ov\Phi([\cE
  M])=E\ov\Phi([M])$, even for all objects $M$. The similar arguments to show $\ov\Phi([\cF M])=F\ov\Phi([M])$
are omitted. 
\end{proof}

\section{The graded version of $\cO$ and the action of $U_q(\mathfrak{sl}_2)$.}
\label{Ograd}
Our main goal is to ``categorify'' the $U$-modules $V_{\bf
  d}$. In Theorem~\ref{HCnongrad}, we already obtained a categorification of the $\cU(\mathfrak{sl}_2)$-modules $\ov{V_{\bf
  d}}$ via certain Harish-Chandra bimodules. The main idea is to introduce a
  graded version, say $\cH_\mu^{gr}$, of
  $\bigoplus_{i=0}^n\;{_{\omega_i}}\cH_\mu^1(\GL_n)$. Then the complexified
  Grothendieck group ${\bf
      G}(\cH_\mu^{gr})$ becomes a
  $\mC(q)$-module, where $q$ acts by shifting the grading. The second step
  will be to introduce exact functors ${\bf E}$, ${\bf F}$, ${\bf K}$, ${\bf K}^{-1}$ which induce an
  $U$-action on ${\bf
      G}(\cH_\mu^{gr})$ such that  ${\bf
      G}(\cH_\mu^{gr})\cong V_{\bf d}$ as $U$-modules. In the present section we consider first the special case $V_1^{\otimes n}$ and work
  with category $\cO$ instead of Harish-Chandra bimodules. The general case
  will easily be deduced afterwards.

\subsection{Graded algebras and modules}
For a ring or algebra
$A$ we denote by $A\MOD$ (resp. $\op{mod-}A$) the category of finitely generated left (resp. right)
$A$-modules. If $A$ is $\mZ$-graded then we denote the corresponding categories of
graded modules by $A\gMOD$ and $\op{gmod-}A$ respectively. For $k\in\mZ$ we denote
by $\langle k\rangle:\op{gmod-}A\rightarrow\op{gmod-}A$ the functor of shifting the
grading by $k$, \ie $(\langle k\rangle M)_i=M_{i-k}$ for $M=\oplus_{i\in\mZ}
M_i\in\op{gmod-}A$. We will normally write $M \langle k\rangle$ instead of
$\langle k\rangle M$. 
For a finite dimensional non-negatively graded algebra $A$ such that $A_0$ is
semisimple, the Grothendieck group $[\op{gmod-}A]$ is the free $\Z[q,q^{-1}]$-module, freely 
generated by a set of isomorphism classes of simple objects in $\op{gmod-}A$ which
give rise to a basis of
$[\op{mod-}A]$ after forgetting the grading. The $\Z[q,q^{-1}]$-module structure 
comes from the grading, namely, for an object $M$ of $\op{gmod-}A$ we define 
$q^i[M] = [M\langle i\rangle]$.

\subsection{The graded version of $\cO$}
\label{graded O}
We first recall (from \cite{BGS}) the graded version of each integral
block of category $\cO$ which is defined by introducing a (Koszul-)grading on the
endomorphism ring of a minimal projective generator in this block. Fixing $P_\la=\bigoplus_{x\in W^\la}P(x\cdot\la)\in{}_\la\cO$ a minimal projective generator defines an
equivalence of categories $\epsilon_\la:{}_\la\cO\cong\op{mod-}\END_\mg(P_\la)$ via
the functor $M\mapsto\HOM_\mg(P_\la,M)$ (see \eg \cite{Bass}). We fix such pairs $(P_\la,
\epsilon_\la)$ for any dominant integral $\la$. Now we have to introduce a
grading on $\END_\mg(P_\la)$. This is done in \cite{BGS} using the connection
between category $\cO$ and modules over the cohomology ring of certain partial
flag varieties (as described in \cite{Sperv}). Which partial flag we have to take depends on
the block of $\cO$ we consider. More precisely we have to do the following: In
each integral block ${}_\la\cO$, there is a unique indecomposable self-dual
projective-injective (tilting) 
module $T_\la$. We consider Soergel's functor
$\mV_\la:{}_\la\cO\rightarrow\op{mod-}\END_\mg(T_\la)$, $M\mapsto\HOM_\cO(T_\la,M)$
which has been introduced in \cite{Sperv}.  By
\cite[Endomorphismensatz]{Sperv}, the algebra  $\END_\cO(T_\la)$ is
canonically isomorphic to the invariants $\H^{W_\la}$
inside the coinvariants $\H=S(\mh)/ (S(\mh)_+^W)$, \ie it is isomorphic to
the cohomology ring of the partial flag variety
corresponding to $\la$ inside the cohomology ring of the full
flag variety (see \eg \cite{Hiller}). 

With the notation from above we also denote by
$T_{i_1,i_2,\cdots ,i_r}\in{}_I\cO$ the unique indecomposable projective-injective
(tilting) module with corresponding functor
$\mV_I:{}_I\cO\rightarrow \op{mod}-\END_\cO(T_I)$ and denote
$C^I=\END_\cO(T_I)$.  We will abbreviate $C^{W_\la}$ as
$C^\la$. 

If we fix a $\mZ$-grading on the algebra $S(\mh)$ by putting $\mh$ in degree
two, then $C$, and hence also any $C^\la$, inherit a grading (which is known
to coincide with the cohomology grading, see \eg \cite{Hiller} and references therein). For any $x\in W^\la$, the module $\mV P(x\cdot\la)$
has a graded lift  (see \cite{BGS}); \ie there exists a module $M\in C^\la\gMOD$ which is
isomorphic to $\mV P(x\cdot\la)$ after forgetting the grading. Since $\END_{C^\la}(\mV
P(x\cdot\la))\cong\END_\cO(P(x\cdot\la))$ (\cite[Struktursatz]{Sperv}), the module $\mV P(x\cdot\la)$ is
indecomposable. Hence a graded lift is unique up to isomorphism and grading
shift  (see \eg \cite[Lemma
2.5.3]{BGS}). We fix for any $x\in W^\la$ a graded lift of $\mV
  P(x\cdot\la)$ such that its lowest degree is $-l(x)$. By abuse of language
we denote this graded lift also  $\mV P(x\cdot\la)$. This defines a
grading on $\mV P_\la=\oplus_{x\in W^\la}P(x\cdot\la)$ and induces a grading
on $\END_\cO(P_\la)=\END_{C^\la}(\mV P_\la)$ by Soergel's structure theorem (\cite{Sperv}). In fact it turns
$\END_\cO(P_\la)$ into a (non-negatively graded) Koszul algebra (see
\cite{BGS}) which we denote by ${}_\la A$. As usual we also write ${}_iA$ instead of ${}_{\omega_i}A$.

\subsection{Graded lifts of modules and functors}
Now we have a graded version for any integral block of $\cO$. We need graded
lifts of modules in ${}_\la\cO$ which are defined as follows: 
\begin{definition}
{\rm
Let $\la$, $\mu$ be dominant and integral weights. Let
$f_\mu:\op{gmod-}{}_\mu A\rightarrow\op{mod-}{}_\mu A$ denote the functor which
forgets the grading.   

\begin{enumerate}
\item Let $\ov M\in{}_\mu\cO$. A {\it
    graded lift of $\ov{M}$} is a module $M\in\op{gmod-}{}_\mu A$ such that
    $f_\la(M)\cong\epsilon_\mu(\ov M)$.
 
\item Let $\overline{F}:{}_\la\cO\rightarrow{}_\mu\cO$ be a functor. A {\it graded lift of $\overline{F}$} is
a functor $F:\op{gmod-}{}_\la A\rightarrow\op{gmod-}{}_\mu A$ such that 
\begin{enumerate}[(i)]
\item $F\langle k\rangle\cong\langle k\rangle F$, and
\item $f_\mu F\cong\epsilon_\mu \overline{F}\epsilon_\la^{-1} f_\la$.
\end{enumerate}
\end{enumerate}
}
\end{definition}

If $\ov M\in{}_\la\cO$ is indecomposable, then a graded lift of $\ov M$, if it
exists, is unique up to isomorphism and grading shift (see \eg \cite[Lemma
2.5.3]{BGS}). In general, for an arbitrary module  $\ov M\in{}_\la\cO$, a graded lift does not have to exist (see \eg
\cite[Section 4]{Stgrad}). However, it
is known that indecomposable projective modules, (dual) Verma modules and simple
modules have graded lifts, see \cite[Section 3.11]{BGS}, \cite[Section
3]{Stgrad}. (For a more general setup for quasi-hereditary algebras we refer to
\cite{BinZhu}). We denote by $\tilde
P(x\cdot\la)$, $\tilde M(x\cdot\la)$, $\tilde
L(x\cdot\la)$ the (uniquely defined up to isomorphism) graded lifts of the modules  $P(x\cdot\la)$, $M(x\cdot\la)$,
$L(x\cdot\la)$ with the property that their heads are concentrated in degree
zero. Let $\tilde\nabla(x\cdot\la)$ be the graded lift of the dual
Verma module $\nabla(x\cdot\la)=\op{d}\Delta(x\cdot\la)$ such that the socle is concentrated in degree
zero. Let $\tilde I(x\cdot\la)$ be the injective hull of $\tilde
L(x\cdot\la)$. Let
$\op{d}=\HOM_\mC(\bullet, \mC)$ denote a graded lift of the duality functor such that $\op{d}(L)=L$ for any
simple module $L$ concentrated in degree zero (for properties see
\eg \cite{Stgrad}).
          
\subsection{Projective functors and the cohomology ring of the flag variety}
The purpose of this section is to give the tools for a construction of graded
lifts of the functors $\cE$ and $\cF$. The following result gives first of all a categorical interpretation of the divided
 powers of $E$ and $F$ and secondly provides an alternative description of the
 functors $\cE$ and $\cF$ which makes it possible to   
 relate them to the cohomology ring of the flag variety in Proposition~\ref{EFandV}
 afterwards. From these results, the desired graded lifts of functors can be
 constricted easily.  
\begin{prop}
\label{EFastransl} Let $0\leq i\leq n$ and $k\in\Z_{>0}$. 
\begin{enumerate}[(a)]
\item There are isomorphisms of projective functors
\begin{eqnarray}
\label{Es}
  \bigoplus_{j=1}^{k!}\cE_i^{(k)}&\cong&\cE_{i+k-1}\cdots\cE_{i+1}\cE_i,\\  
\label{Fs}  
 \bigoplus_{j=1}^{k!}\cF_{i}^{(k)}&\cong&\cF_{i+k+1}\cdots\cF_{i-1}\cF_i. 
\end{eqnarray}
\item Let $k\in\mZ_{>0}$. There are isomorphisms of indecomposable projective functors
  \begin{eqnarray}
\label{ET}
  \cE_i^{(k)}&\cong& \mathds{T}_{i,i+k}^{i+k}\,\mathds{T}_{i}^{i,i+k},\\  
\label{FT}
  \cF_{i}^{(k)}&\cong&\mathds{T}_{i,i-k}^{i-k}\,\mathds{T}_{i}^{i,i-k}. 
  \end{eqnarray}
\end{enumerate}

\end{prop}
\begin{proof}
  By adjointness properties it is enough to prove the formulas \eqref{Es} and
  \eqref{ET}.  Recall (see \eg \cite[4.6 (1)]{Ja2}) that for any $\la\in\mh^\ast$ and any finite
  dimensional $\mg$-module $E$ we have 
$[M(\la)\otimes E]=[\oplus_{\nu} M(\la+\nu)]$,
  where  $\nu$ runs through the multiset of weights of $E$. From the definition of $\cE^{(k)}_i$ we get in particular 
$[\cE^{(k)}_i(M({\bf a}))]=[\oplus_{\bf a'}M({\bf a'})]$, where ${\bf
  a'}=(a_1',a_2',\ldots, a_n')$ runs
  through the set of $\{0,1\}$-sequences containing exactly $i+k$ ones and
  where $a_j=1$ implies $a_j'=1$ for $1\leq j\leq n$. In particular  
  $[\oplus_{j=1}^{k!}\cE^{(k)}_i(M(\la_i))]=[\cE_{i+k-1}\cdots\cE_{i+1}\cE_{i}(M(\la_i))]$
  (see \cite[Proposition 6]{BFK}). Hence, \eqref{Es} follows from the
  classification of projective functors (\cite{BG}).
To prove the second part let first be $k=1$. The formula \cite[Proposition 6]{BFK} shows that
  $\cE_i(M(\omega_i))$ has a Verma flag with subquotients isomorphic to
  $M(x\cdot\omega_{i+1})$, where $x\in W_{i}/W_{i,i+1}$ is a shortest coset
  representative. The same is true for the module
  $\mathds{T}_{i,i+1}^{i+1}\,\mathds{T}_{i}^{i,i+1}\,M(\omega_i)$ (using \cite[4.13 (1), 4.12 (2)]{Ja2}). Since both modules
  are projectives, they are isomorphic. The classification theorem of
  projective functors (\cite[Section 3]{BG}) provides the required isomorphism of
  functors for $k=1$.  Using again the
  formulas \cite[4.13 (1), 4.12 (2)]{Ja2} we
  easily get
  $[\cE_i^{(k)}(M(\la_i))]=[\mathds{T}_{i,i+k}^{i+k}\,\mathds{T}_{i}^{i,i+k}\,M(\la_i)]$. The
  classification theorem of projective functors (\cite{BG}) implies
  $\cE_i^{(k)}\cong \mathds{T}_{i,i+k}^{i+k}\,\mathds{T}_{i}^{i,i+k}$. Note that $M(\la_{i+k})$ occurs with multiplicity
  one in $\cE_i^{(k)}\,M(\la_i)$, hence $\cE_i^{(k)}(M(\la_i))$ is
  indecomposable, and therefore so is $\cE_i^{(k)}$ (again by the
  classification of projective functors \cite{BG}).    
\end{proof}

We have restriction functors $\RES_J^I:\H^I-\op{mod}\rightarrow \H^J-\op{mod}$ if
$J\subseteq I$ and $\RES_\mu^\la:\H^\la-\op{mod}\rightarrow \H^\mu-\op{mod}$
if $W_\la\subseteq  W_\mu$.

\begin{prop}
\label{EFandV}
  Let $\la$, $\mu$ be dominant and integral weights such that $W_\la\subseteq
  W_\mu$. There are isomorphisms of functors
  \begin{eqnarray*}
    \mV_\mu\, \mathds{T}_\la^\mu&\cong&\RES^\la_\mu \mV_\mu,\\
    \mV_\la\,  \mathds{T}_\mu^\la&\cong&C^\la\otimes_{C^\mu}\mV_\mu(\bullet).
  \end{eqnarray*}
In particular
  \begin{eqnarray*}
    \mV_{i+k}\,\cE^{(k)}_i&\cong& \RES_{i+k}^{i,i+k}\, \H^{i,i+k}\otimes_{\H^i}\mV_{i}(\bullet),\\
    \mV_{i-k}\,\cF^{(k)}_i&\cong& \RES_{i-k}^{i,i-k}\, \H^{i,i-k}\otimes_{\H^i}\mV_{i}(\bullet),
  \end{eqnarray*}
for any $k\in\mZ>0$. 
\end{prop}

\begin{proof}
This is \cite[Theorem 12, Proposition 6]{S2} together with Proposition~\ref{EFastransl}.
\end{proof}

 To keep track of the grading, we first need the following
well-known, but crucial fact that $C^I$ is a free $C^J$-module of finite rank
whenever $J\subseteq I$. More precisely we need the following statement: 

\begin{lemma}
\label{free}
  Let $1\leq i\leq (n-1)$. There are isomorphisms of graded $C^i$-modules 
  \begin{eqnarray*}
    C^{i,i+1}&\cong&\bigoplus_{r=0}^{n-i-1} C^i\langle 2r\rangle,\\
    C^{i,i-1}&\cong&\bigoplus_{r=0}^{i-1} C^i\langle 2r\rangle.
  \end{eqnarray*}   
\end{lemma}
\begin{proof}
  By classical invariant theory (\cite[II.3]{Hiller}), $C^{i,i\pm 1}$ is a free $C^i$-module
  of rank $|W_i/W_{i,i\pm 1}|$, and a basis can be chosen homogeneous in the
  degrees length of $x$, where $x$ runs through the set of shortest coset
  representatives from $W_i/W_{i,i\pm 1}$. 
\end{proof}

The following adjoint pairs of functors will be used later

\begin{prop}
\label{adjunctions}
  Let $\la$, $\mu$ be dominant and integral. Assume $W_\la\subseteq W_\mu$. Then there are pairs of adjoint
  functors 
  \begin{eqnarray*}
    (C^{\la}\otimes_{C^{\mu}}\bullet\;,\RES^{\la}_{\mu}\;)&{\mbox and}&(\RES^{\la}_{\mu}\;,C^{\la}\otimes_{C^{\mu}}\bullet\;)
  \end{eqnarray*}
between $C^{\la}\MOD$ and $C^{\mu}\MOD$;
and 
\begin{eqnarray}
\label{adjg}
   (C^{\la}\otimes_{C^{\mu}}\bullet\;,\RES^{\la}_{\mu}\;)&{\mbox
   and}&(\RES^{\la}_{\mu}\;,C^{\la}\otimes_{C^{\mu}}\bullet\;\langle -{\op{max}}\rangle)
\end{eqnarray}
considered as functors between $C^{\la}\gMOD$ and $C^{\mu}\gMOD$. Here
$\op{max}\in I$ denotes the maximal element in $I$ where $C^{\la}\cong\oplus_{i\in I}
C^{\mu}\langle i\rangle$ as graded $C^{\mu}$-module.
\end{prop}

\begin{proof}
By Proposition~\ref{EFandV}, the pairs of adjoint functors $(\mathds{T}^{\la}_\mu,\mathds{T}^{\mu}_\la)$ and
$(\mathds{T}_\la^\mu,\mathds{T}_\mu^\la)$ give rise to the pairs
$(C^{\la}\otimes_{C^{\mu}}\bullet\;,\RES^{\la}_{\mu}\;)$ and
$(\RES^{\la}_{\mu}\;,C^{\la}\otimes_{C^{\mu}}\bullet\;)$, since $\mV$ is a quotient functor.
 However, for the
  graded version we have to be more explicit.
  For $M\in C^{\mu}\MOD$ we consider the inclusion $i_M:
  M\hookrightarrow C^{\la}\otimes_{C^\mu}M$, $m\mapsto 1\otimes m$. This is a
  $C^{\mu}$-morphism and defines a map 
  \begin{eqnarray*}
    \Phi_{M,N}:\;\HOM_{C^{\la}}(C^{\la}\otimes_{C^\mu}M,
    N)&\rightarrow&\HOM_{C^{\mu}}(M,\RES^{\la}_{\mu}N)\\
    f&\mapsto& f\circ i_M
  \end{eqnarray*}
for any $N\in C^{\la}\MOD$. The map is functorial in $M$ and $N$, and it is
injective, since $f\circ i_M=0$ implies $f(c\otimes m)=cf(1\otimes m)= c(f\circ
i_M(m))=0$ for any $m\in M$, $c\in C^{\la}$. To show the surjectivity it is
enough to compare the dimensions. Moreover, the functors are exact
and additive. Therefore, it is sufficient to consider the case $M=C^{\mu}$
(for general $M$ choose a free $2$-step resolution and use the Five-Lemma). In this case we have 
$\DIM\HOM_{C^{\la}}(C^{\la}\otimes_{C^\mu}C^{\mu},N)
=\DIM\HOM_{C^{\la}}(C^{\la},N)=\DIM N=
\DIM\HOM_{C^{\mu}}(C^{\mu},\RES^{\la}_{\mu}N)$. This proves the first
adjunction. Moreover, the adjunction is compatible with the grading, since $i_M$ is homogeneous of
degree zero for graded $C^\mu$-modules $M$. 

Let
  $p_N:C^{\la}\otimes_{C^{\mu}}N \cong \oplus_{i\in I} N \langle
  i\rangle\rightarrow N\langle \op{max}\rangle$ be the projection for $N\in C^\mu\gMOD$. It defines
  a natural morphism
  \begin{eqnarray*}
    \Phi'_{M,N}:\;
    \HOM_{C^{\la}}(M,C^{\la}\otimes_{C^{\mu}}N)&\rightarrow&\HOM_{C^{\mu}}(\RES^{\la}_{\mu}M,N)\\
f&\mapsto& p_N\circ f.
  \end{eqnarray*}
We show that $\Phi'$ is injective: Assume $p_N\circ f=0$. Then $f(M)\subseteq
  \bigoplus_{i\in I-{\op{max}}} N \langle i\rangle$. On the other hand, for
  any $n\in C^{\la}\otimes_{C^{\mu}}N$ there exists a $c\in C$ such that
  $cn\not\in\oplus_{i\in I-\{\op{max}\}} N \langle i\rangle$. Hence $f=0$ and
  the injectivity of $\Phi'_{M,N}$ follows. 
The map $\Phi'_{M,N}$ is surjective for $M=C^{\la}$, because 
  $\DIM\HOM_{C^\la}(C^{\la},C^{\la}\otimes_{C^{\mu}}N)=\DIM(C^{\la}\otimes_{C^{\mu}}
  N)=i\cdot\DIM N=\DIM\HOM_{C^{\mu}}(C^{\la},N)$, where $i=|W_\mu/W_\la|$ is
  the rank of $C^\la$ as $C^\mu$-module. The surjectivity in general follows then from
  the Five-Lemma. By construction, the isomorphism $\Phi'$ is homogeneous of
  degree $\op{max}$. The existence of the second adjunction from \eqref{adjg} follows.          
\end{proof}

\subsection{A functorial action of $U_q(\mathfrak{sl}_2)$}
We are prepared to introduce the graded lifts of our functor $\cE$ and $\cF$. We define 
\begin{eqnarray}
  \label{DefE}
  {\bf E}^{(k)}_i&=&\HOM_{C^{i+k}}\big(\mV P_{i+k},\RES_{i+k}^{i,i+k}
  C^{i,i+k}\otimes_{C^i}\mV P_i\langle -r_{i,k}\rangle\big)\\
 &{\mbox with}&r_{i,k}=\sum_{r=i}^{i+k-1}(n-r-1);\nonumber\\
 {\bf F}^{(k)}_i&=&\HOM_{C^{i-1}}\big(\mV P_{i-k},\RES_{i-k}^{i,i-k}
  C^{i,i-k}\otimes_{C^i}\mV P_i\langle -r'_{i,k}\rangle\big)\label{DefF}\\
&{\mbox with}& r'_{i,k}=\sum_{r=i}^{i+k-1}r-1.\nonumber
\end{eqnarray}
We have ${\bf E}^{(k)}_i\in\END_{C^i}(\mV P_i)\op{-gmod-}\END_{C^{i+k}}(\mV P_{i+k})$ via
$g.f.h=(\op{Id}\otimes g)\circ f\circ h$ for $g\in\END_{C^i}(\mV P_i)$,
$h\in\END_{C^{i+k}}$ and $f\in{\bf E}^{(k)}_i$. From the
definitions we may also consider ${\bf E}^{(k)}_i$ as an object in
${}_i A\op{-gmod-}{}_{i+k}A$. Tensoring with ${\bf E}^{(k)}_i$ defines a functor ${\bf
  E}^{(k)}_i: \op{gmod-}{}_i A\rightarrow\op{gmod-}{}_{i+k}A$ (which we denote, abusing
language, by the same
symbol). Analogously, ${\bf F}^{(k)}_i$
defines a functor ${\bf F}^{(k)}_i:\op{gmod-}{}_i A\rightarrow\op{gmod-}{}_{i-k}A$. Set ${\bf
  E}^{(k)}=\oplus_{i=0}^{n}{\bf E}^{(k)}_i$, considered as an endofunctor of
$\oplus_{i=0}^{n} \op{gmod-}{}_i A$. Similarly,  ${\bf
  F}^{(k)}=\oplus_{i=0}^{n}{\bf F}^{(k)}_i$. Let ${\bf K_i}=\langle 2i-n\rangle$
be the endofunctor of $\op{gmod-}{}_i A$ which shifts the degree by $(2i-n)$ and ${\bf
  K}=\oplus_{i=0}^{n}{\bf K}_i$. \\

The following result is the crucial step towards our main categorification
theorem (Theorem~\ref{cat}):
\begin{theorem}
\label{funcrel}
  \begin{enumerate}[(a)]
  \item \label{funcrel1} The functors ${\bf E}_i^{(k)}$ and ${\bf F}_i^{(k)}$ are graded lifts
  of $\cE_i^{(k)}$ and $\cF_i^{(k)}$ respectively. 
\item \label{funcrel2} The functors ${\bf E}$, ${\bf F}$, and ${\bf K}$ satisfy the relations
\begin{eqnarray*} 
  {\bf K}{\bf E}&\cong&{\bf E}{\bf K}\langle 2\rangle,\\ 
  {\bf K}{\bf F}&\cong&{\bf F}{\bf K}\langle -2\rangle,\\
  {\bf K} {\bf K}^{-1}\quad\cong&\op{Id}&\cong\quad{\bf K}^{-1}{\bf K}, \\
  {\bf E}_{i-1}{\bf F}_i\oplus\bigoplus_{r=0}^{n-i-1}\op{Id}\langle 
  n-1-2r-2i\rangle&\cong&{\bf F}_{i+1}{\bf
  E}_{i}\oplus\bigoplus_{r=0}^{i-1}\op{Id}\langle 2i-n-2r-1\rangle. 
  \end{eqnarray*} 
\item \label{funcrel3} In the Grothendieck group we have the equality 
\begin{eqnarray*}
(q-q^{-1})\big({\bf E}_{i-1}^{\bf G}{\bf
  F}_i^{\bf G}-{\bf F}_{i+1}^{\bf G}{\bf E}_i^{\bf G}\big)
&=&{\bf K}_i^{\bf G}-({\bf K}_i^{-1})^{\bf G}.
\end{eqnarray*}

Moreover ${\bf E}_{i-1}{\bf F}_i$ is a summand of ${\bf F}_{i+1}{\bf
  E}_i$ if $n-2i>0$. Likewise,  ${\bf F}_{i+1}{\bf
  E}_i$ is a summand of  ${\bf E}_{i-1}{\bf F}_i$ if $n-2i<0$.
 \end{enumerate}
\end{theorem}

\begin{proof}
  The first statement follows from the definitions, as well as ${\bf
  K}_{i+1}{\bf E}_i\cong {\bf E}_{i}{\bf K}_{i}\langle 2\rangle$, and ${\bf
  K}_{i-1}{\bf F}_i\cong {\bf F}_{i}{\bf K}_{i}\langle -2\rangle$. Obviously ${\bf K}$ is an
  auto-equivalence with inverse ${\bf K}^{-1}$. To prove the remaining isomorphisms
  we first claim that 
  \begin{eqnarray}
  \label{decomp}
    \cE_{i-1}\,\cF_i\cong\;G\oplus\; \bigoplus_{r=1}^{i-1}\op{Id}\quad
    &{\mbox{and}}&\quad\cF_{i+1}\,\cE_i\cong\; G\oplus\;\bigoplus_{r=1}^{n-i-1}\op{Id}
  \end{eqnarray}
for some indecomposable endofunctor $G$ of ${}_{i;n}\cO$. 
Let $M({\bf a})$ be the projective Verma module in ${}_{i;n}\cO$, \ie ${\bf
  a}=(1,\ldots,1,0,\ldots,0)$ with $i$ ones. Then 
\begin{eqnarray*}
  [\cE_{i-1}\cF_iM({\bf a})]&=&\big[\cE_{i-1}\oplus_{k=1}^i M({\bf a(k)})\big]\\
  &=&\big[\oplus_{k=1}^{i}\oplus_{l=i+1}^n M({\bf
  a(k,l)})]+[\oplus_{k=1}^i M({\bf a})],
\end{eqnarray*}
where ${\bf a(k)}, {\bf
  a(k,l)}\in\mZ^n$ such that  
\begin{eqnarray*}
  {\bf a(k)}_j&=&
  \begin{cases}
    0&\text{ if $i<j\leq n$ or $j=k$},\\
    1&\text{ otherwise},
  \end{cases}\\
  {\bf a(k,l)}_j&=&
  \begin{cases}
    0&\text{ if $i<j\not=l\leq n$ or $j=k$},\\
    1&\text{ otherwise}.
  \end{cases}
\end{eqnarray*}
On the other hand 
\begin{eqnarray*}
  [\cF_{i+1}\cE_iM({\bf a})]&=&\big[\cF_{i+1}\oplus_{l=i+1}^n M({\bf b(l)})\big]\\
  &=&\big[\oplus_{k=1}^{i}\oplus_{l=i+1}^n M({\bf
  b(l,k)})]+[\oplus_{k=i+1}^n M({\bf a})],
\end{eqnarray*}
where ${\bf b(l)}, {\bf
  b(l,k)}\in\mZ^n$ such that 
\begin{eqnarray*}
  {\bf b(l)}_j&=&
  \begin{cases}
    1&\text{ if $1\leq j\leq i$ or $j=l$},\\
    0&\text{ otherwise},
  \end{cases}\\
  {\bf b(l,k)}_j&=&
  \begin{cases}
    1&\text{ if $1\leq j\not=k\leq i$ or $j=l$},\\
    0&\text{ otherwise},
  \end{cases}\\
&=&
 \begin{cases}
    0&\text{ if $i<j\not=l\leq n$ or $j=k$},\\
    1&\text{ otherwise}.
  \end{cases}
\end{eqnarray*}
In particular, $[\cF_{i+1}\cE_iM({\bf
  a})]=\big[\oplus_{k=1}^{i}\oplus_{l=i+1}^n M({\bf
  a(k,l)})]+[\oplus_{k=i+1}^n M({\bf a})]$. 

Let $P$ be the projective cover of
$M((0,1,\ldots,1,0,\ldots, 0,1))\in{}_{i;n}\cO$. An easy calculation (in the Hecke algebra) shows that 
$[P]=\big[\oplus_{k=1}^{i}\oplus_{l=i+1}^n M({\bf
  a(k,l)})]+[M({\bf a})]$. The decompositions \eqref{decomp}
follow then from the classification theorem (\cite{BG}) of projective
functors. The functor $G$ is the indecomposable projective functor which sends
$M({\bf a})$ to $P$. 

Now we have to consider the graded picture. By Lemma~\ref{free} we get
isomorphisms of graded $C^i$-modules
\begin{eqnarray}
\label{Cdecomp}
  C^{i,i+1}\otimes_{C^{i+1}}C^{i,i+1}\otimes_{C^i}\mQ&\cong&
  \oplus_{l=0}^{i}\oplus_{k=0}^{n-i-1}\mQ\langle2k+2l\rangle,\nonumber\\
 C^{i,i-1}\otimes_{C^{i-1}}C^{i,i-1}\otimes_{C^i}\mQ&\cong&
  \oplus_{l=0}^{n-i}\oplus_{k=0}^{i-1}\mQ\langle2l+2k\rangle.
\end{eqnarray}
Since their lowest degrees coincide and $G$ is indecomposable, the
decompositions~\eqref{decomp} give rise to isomorphisms of endofunctors of
$\op{gmod-}{}_i A$ 
  \begin{eqnarray}
  \label{decompgrad}
    {\bf E_{i-1}}\,{\bf F_i}\cong\;{\bf G}\oplus\;
    \bigoplus_{r=0}^{i-1}\op{Id}\langle m_r\rangle\quad
    &{\mbox{and}}&\quad{\bf F}_{i+1}\,{\bf E}_i\cong\;{\bf
    G}\oplus\;\bigoplus_{r=0}^{n-i-1}\op{Id}\langle n_r\rangle,
  \end{eqnarray}
where ${\bf G}$ is a certain graded lift of $G$ and $m_r$, $n_r\in\mZ$. 
The formulas \eqref{Cdecomp}, together with the definition of the ${\bf E}_i$,
${\bf F}_i$ and the formula $\HOM(M,N\langle i\rangle)=\HOM(M,N)\langle
-i\rangle$ for graded morphisms between graded modules $M$, $N$, imply that, if $n-2i\geq 0$, then  
\begin{eqnarray*}
{\bf F}_{i+1}^{\bf G}{\bf E}_i^{\bf G}-{\bf E}_{i-1}^{\bf G}{\bf
  F}_i^{\bf G}
&=&
\big(\oplus_{k=0}^{n-i-1}\langle-(2k+2i)\rangle\langle-(-n+1)\rangle\big)^{\bf G}\\
&&-\;\big(\oplus_{k=0}^{i-1}\langle-(2(n-i)+2k)\rangle\langle n-1\rangle\big)^{\bf G}\\
&=&
\begin{cases}
0&\text{ if $n-2i=0$,}\\
\big(\oplus_{k=0}^{n-i-1}\langle -(2i+2k)\rangle\langle n-1\rangle\big)^{\bf G}&\text{ if $n-2i>0$,}
\end{cases}\\
&=& [n-2i]\op{id}.
\end{eqnarray*}
If $n-2i<0$, then 
\begin{eqnarray*}
{\bf E}_{i-1}^{\bf G}{\bf
  F}_i^{\bf G}-{\bf F}_{i+1}^{\bf G}{\bf E}_i^{\bf G}
&\cong&
\big(\oplus_{k=0}^{i-1}\langle -(2(n-i)+2k)\rangle\langle-(-n+1)\rangle)^{\bf G}\\
&-&\big(\oplus_{k=0}^{n-i-1}\langle-(2k+2i)\rangle\langle n-1\rangle\big)^{\bf G}\\
&=&\big(\oplus_{k=0}^{n-i-1}\langle -(n-2i+1+2k)\rangle\big)^{\bf G}\\
&=& [2i-n]\op{id}.
\end{eqnarray*}
In particular, 
\begin{eqnarray*}
(q-q^{-1})\big({\bf E}_{i-1}^{\bf G}{\bf
  F}_i^{\bf G}-{\bf F}_{i+1}^{\bf G}{\bf E}_i^{\bf G}\big)
&=&{\bf K}_i^{\bf G}-({\bf K}_i^{-1})^{\bf G},
\end{eqnarray*}
and the formula of part~\eqref{funcrel3} hold.

Since for $n-2i\geq 0$ the element ${\bf F}_{i+1}^{\bf G}{\bf E}_i^{\bf G}-{\bf E}_{i-1}^{\bf G}{\bf
  F}_i^{\bf G}$ is a (positive) sum of certain $(\langle k\rangle)^{\bf G}$ with $k\in\mZ$, the
  decomposition~\eqref{decompgrad} implies that 
  ${\bf E}_{i-1}{\bf F}_i$ is a summand of ${\bf F}_{i+1}{\bf
  E}_i$. Likewise,  ${\bf F}_{i+1}{\bf
  E}_i$ is a summand of  ${\bf E}_{i-1}{\bf F}_i$ if $n-2i<0$. If now $n-2i>0$ then the formula from~\eqref{funcrel3} implies  
\begin{eqnarray*}
    {\bf E}_{i-1}{\bf F}_i\oplus\bigoplus_{k=0}^{n-2i-1}\op{Id}\langle
    n-2i-1-2k\rangle\cong {\bf F}_{i+1}{\bf
  E}_i.
\end{eqnarray*}
If $n-2i<i$ then we get 
\begin{eqnarray*}
   {\bf F}_{i+1}{\bf E}_i\oplus\bigoplus_{k=0}^{2i-n-1}\langle
   2i-n-1-2k\rangle\cong{\bf E}_{i-1}{\bf F}_i. 
\end{eqnarray*}
From this we finally deduce the last formula in part~\eqref{funcrel2}.
\end{proof}

\section{Harish-Chandra bimodules and graded category $\cO$}
\label{mainresult}
\subsection{The categorification theorem}
In this section we deduce the main result which provides a categorification of
arbitrary finite tensor products $V_{\bf d}$ of finite dimensional
$U_q(\mathfrak{sl}_2)$-modules. \\
This will be in fact a rather easy consequence from the previous
Theorem~\ref{funcrel} using the embedding of categories
$_\la\cH_\mu^1\hookrightarrow {}_{\la}\cO$ for any reductive Lie algebra and dominant integral weights $\la$, $\mu$ from \cite{BG}. The image of this
functor is a full
subcategory of ${}_\la\cO$ given by  $\cP_\mu$-presentable objects; an object $M$
is {\it $\cP_\mu$-presentable} if there is an exact sequence of the
form $P_2\rightarrow P_1\rightarrow M\rightarrow 0$ such that $P_2$, $P_1$ are
direct sums of projective modules, each indecomposable summand isomorphic to
some $P(x\cdot\la)$ such that $x$ is a longest double coset representative in
$W_\mu \backslash W/W_\la$. Let $\cA_i^\mu=\cA_i^\mu(\mathfrak{gl}_n)$ be the
full subcategory of $\op{gmod-}{}_i A$ given by $\cP_\mu$-presentable objects; \ie by objects $M$ such that there is an exact
sequence (in $\op{gmod-}{}_i A$) of the
form $P_2\rightarrow P_1\rightarrow M\rightarrow 0$, where $P_2$, $P_1$ are
direct sums of projective modules, each indecomposable summand after
forgetting the grading isomorphic to
some $\epsilon_\la(P(x\cdot\la))$ where $x$ is a longest double coset
representative in $W_\mu\backslash W/W_i$. The category $\cA_i^{\mu}$ is
abelian, since the categories $_\la\cH^1_\mu$ are abelian. (It is not
completely trivial to describe this abelian structure of $\cA_i^{\mu}$ when
realized as a subcategory of $\op{gmod}-{}_i A$, see \eg \cite{MS}).        

The previous Theorem~\ref{funcrel} implies the following main result
\begin{theorem}[Categorification Theorem]
\label{cat}
Let ${\bf d}$ be a composition of $n$ and let $\mu\in\mh^\ast$ be dominant and
  integral such that $W_\mu\cong S_{\bf d}$. 
  There is an isomorphism of $U_q(\SL_2)$-modules
  \begin{eqnarray*}
{\bf G}(\bigoplus_{i=0}^n\cA_i^{\mu}(\mathfrak{gl}_n))\cong V_{\bf d}.    
  \end{eqnarray*}
where the left hand side becomes a $U_q(\SL_2)$-module structure via the
induced action of the exact functors ${\bf E}$, $\bf{F}$, and $\bf{K}$.  
\end{theorem}
\begin{proof}
  We only have to show that the functors in question preserve the category
  $\bigoplus_{i=0}^n\cA_i^{\mu}$, then the statement follows from the
  previous Theorem~\ref{funcrel}, Proposition~\ref{lemmatensor} and Theorem~\ref{HCnongrad}. Recall that the
  embedding $_\la\cH_\mu^1\hookrightarrow {}_{\la}\cO$ is given by $X\mapsto
  X\otimes_{\cU(\mg)} M(\mu)$. In particular, it commutes with
  translation functors. Hence $\bigoplus_{i=0}^n\cA_i^{\mu}$ is stable under ${\bf E}$ and ${\bf
  F}$. That it is also stable under grading shifts, in particular under ${\bf
  K}$, follows directly from the definitions. 
\end{proof}

Additionally, we are able to give a categorical interpretation of the
 involutions introduced in Section~\ref{quantum-sl2}. This will be the topic
 of the following subsections. 
\subsection{The anti-automorphism $\tau$ as taking left adjoints} 

The anti-automorphism $\tau$ (and its inverse) can be considered as the
operation of taking left (respectively right) adjoints: 

\begin{prop} \label{all-adjoint} 
There are pairs of adjoint functors 
\begin{eqnarray*}
({\bf E}, {\bf F K}^{-1}\langle 1\rangle),\quad ({\bf F}, {\bf EK}\langle 1\rangle), \quad
({\bf K},{\bf K}^{-1}),\quad ({\bf K}^{-1}
,{\bf K})
\end{eqnarray*}
and $(\langle -k\rangle, \langle k\rangle)$ for $k\in\mZ$.
\end{prop} 

\begin{proof}
  This follows directly from the definitions \eqref{DefE} and \eqref{DefF}
  using the Proposition~\ref{adjunctions} and Lemma~\ref{free}.
\end{proof}

\subsection{The Cartan involution $\sigma$ as an equivalence of categories}

The involution $\sigma$ (see \eqref{sigma}) has the following functorial interpretation:

\begin{prop}
\label{Cartan}
  There is an equivalence of categories 
  \begin{align*}
    {\hat\sigma}:\quad\bigoplus_{i=0}^n\op{gmod-}{}_i
    A\rightarrow\bigoplus_{i=0}^n\op{gmod-}{}_i A
  \end{align*}
such that 
\begin{align*}
  {\Sigma}({\bf E})&\cong{\bf F},&{\Sigma}{\bf F}&\cong{\bf E},&{\Sigma}({\bf
  K})&\cong{\bf K}^{-1},&{\Sigma}(\langle k\rangle)&\cong\langle k\rangle
\end{align*}
where $k\in\mZ$ and ${\Sigma}(F)={\hat{\sigma}} F{\hat{\sigma}}^{-1}$ for any endofunctor $F$ of
$\bigoplus_{i=0}^n\op{gmod-}{}_i A$. Moreover $\Sigma(G_1G_2)\cong\Sigma(G_1)\Sigma(G_2)$ and
${\Sigma}^2(G_1)\cong G_1$ for any  endofunctor $G_1$, $G_2$ of
$\bigoplus_{i=0}^n\op{gmod-}{}_i A$. 
\end{prop}

\begin{proof}
  Let $t:W\rightarrow W$ be the isomorphism given by $s_i\mapsto
  s_{n-i}$. By \cite[Theorem 11]{Sperv}, this induces an equivalence of
  categories 
  $\hat\sigma_i:\op{mod-}{}_i A\rightarrow\op{mod-}{}_{n-i} A$ such that $\hat{\sigma}_i
  M(x\cdot\omega_i)\cong M(t(x)\cdot\omega_{n-1})$ which lifts even to an
  equivalence of categories $\hat{\sigma}_i:\op{gmod-}{}_i A\rightarrow\op{gmod-}{}_{n-i}A$. 
  In particular $\hat{\sigma}_{i}{\bf K}_i\cong{\bf K}_{n-i}^{-1}\hat{\sigma}_i$. Set
  $\omega=\oplus_{i=0}^n\hat{\sigma}_i$. From the definitions we get ${\Sigma}(\langle
  k\rangle)\cong\langle k\rangle$. We also have $\hat{\sigma}_{i+1}\,{\bf
  E}_i\cong{\bf F}_{n-i}\,\hat{\sigma}_i$ after forgetting the grading, hence $\hat{\sigma}_{i+1}\,{\bf
  E}_i\cong{\bf F}_{n-i}\,\hat{\sigma}_i\langle j\rangle$ for some $j$, because the
  involved functors are indecomposable (see Proposition~\ref{EFastransl}). A
  direct calculation in the Grothendieck group shows that in fact $\hat{\sigma}_{i+1}\,{\bf
  E}_i\cong{\bf F}_{n-i}\,\hat{\sigma}_i$ and $\hat\sigma_{i-1}\,{\bf
  F}_i\cong{\bf E}_{n-i}\,\hat\sigma_i$. Since $t$ is an involution, we get that $\hat{\sigma}$ is an involution as well and so is $\Sigma$ by
  definition. The formula $\Sigma(G_1G_2)\cong\Sigma(G_1)\Sigma(G_2)$ is then also clear.  
\end{proof}

\subsection{The involution $\psi$ as a duality functor.} 
\label{psi}

For $\la\in\mh^\ast$ dominant and integral let $\op{d}=\HOM_\mQ(\cdot,\mQ):C^\la\gMOD\rightarrow
C^\la\gMOD$ be the graded duality, ie. $\op{d}(M)_i=\HOM_{\mQ}(M_{-i},\mQ)$. With
the conventions on the graded lifts of $\mV_\la P(x\cdot\la)$ we get in
particular that $\mV P(x\cdot\la)\cong\op{d}\mV P(x\cdot\la)$ is self-dual
(\cite[Lemma 9]{Sperv}). Hence, $\op{d}$ defines an isomorphism of graded algebras
${}_\la A\cong{}_\la A^{\op{opp}}$. We get a contravariant duality
$\op{d}_\la:\op{gmod-}{}_\la A\rightarrow{}_\la A\gMOD\cong\op{gmod-}{}_\la A$,
$M\mapsto\HOM_{\mQ}(M,\mC)$, where $\HOM_{\mQ}(M,\mC)$ has the dual
grading, that is $\HOM_{\mC}(M,\mC)_i=\HOM_{{}_\la A}(M_{-i}
,\mC)$. Let $\op{d}=\oplus_{i=0}^n \op{d}_{\omega_i}$ be the duality on
$\oplus_{i=0}^n\op{gmod-}{{}_i A}$. Put $\op{d}_i'=\langle
2i(n-i)\rangle\op{d}_{\omega_i}$. Of course, $\op{d}$ and $\op{d}'$ are
involutions. We first mention an important fact:
\begin{lemma}
\label{duality}
There are isomorphisms of functors 
\begin{eqnarray*}
  {\bf E}_i\op{d}_i'&\cong&\op{d}_{i+1}'{\bf E}_i,\\
  {\bf F}_i\op{d}_i'&\cong&\op{d}_{i-1}'{\bf F}_i,\\ 
  {\bf K}_i\op{d}_i'&\cong&\op{d}_i'{\bf K}_i^{-1},
\end{eqnarray*}
for any $0\leq i\leq n$. 
\end{lemma}

\begin{proof}
 Set $\op{d}_i=\op{d}_{\omega_i}$. After forgetting the grading we have
  $\op{d}_{i+1} {\bf E}_i\cong{\bf E}_{i} \op{d}_i$ (\cite[4.12
  (9)]{Ja2}). Since, considered as a functor from ${}_{i;n}\cO$ to
  ${}_{i+1;n}\cO$, the functor $\op{d}_{i+1}\cE_i\op{d}_i\cong\cE_i$ is indecomposable
  (Proposition~\ref{EFastransl}), we get $\op{d}_{i+1}{\bf E}_i\cong{\bf E}_{i} \op{d}_i\langle k\rangle$
  for some $k\in\mZ$ (\cite[Lemma 2.5.3]{BGS} applied to ${}_{i+1}A\otimes
  {}_i A^{\op{opp}}$). Hence, it is enough to prove that there exists some
  $M\in\op{gmod-}{}_i A$
  which is not annihilated by $\op{d}_{i+1}{\bf E}_i$ and satisfies $\op{d}_{i+1}'{\bf
  E}_i(M)\cong {\bf E}_{i} \op{d}_i'(M)$. Let $M\in\op{gmod-}{}_i A$ be the
  graded lift of the projective module $P(w_0\cdot\omega_i)$ with head
  concentrated in degree zero. From the definition of $\op{d}_i'$ it follows
  that $\op{d}_i'M\cong M$, since $P(w_0\cdot\omega_i)\in\cO$ is self-dual and
  the grading filtration of $M$ is of length $2i(n-i)$ (\cite[Theorem 3.11 (ii)]{BGS}).  By definition we have 
  \begin{eqnarray*}
     {\bf E}_i M&=&M\otimes_{{}_i A}\HOM_{C^{i+1}}(\mV_{i+1}
     P_{i+1},C^{i,i+1}\otimes_{C^i}\mV_i P_i\langle -(n-i-1)\rangle)\\
&=& M\otimes_{{}_i A}\HOM_{C^{i+1}}(\mV_{i+1} P_{i+1},C^{i,i+1}\otimes_{C^i}\mV_i T_i\langle -(n-i-1)\rangle).
  \end{eqnarray*}
By definition of the ${}_i A$-action we only have to consider 
$$M\otimes_{{}_i
     A}\HOM_{C^{i+1}}(\mV_{i+1} P_{i+1},C^{i,i+1}\otimes_{C^i}\mV_i T_i\langle
     -(n-i-1)\rangle).$$ 
From Lemma~\ref{free} we have isomorphisms 
     \begin{eqnarray*}
&&     \HOM_{C^{i+1}}(\mV_{i+1}P_{i+1},C^{i,i+1}\otimes_{C^i}\mV_i T_i\langle
     -(n-i-1)\rangle)\\
&\cong&\HOM_{C^{i+1}}(\mV_{i+1}P_{i+1},\oplus_{k=0}^i C^{i+1}\langle -(n-i-1)\rangle\langle-i(n-i)\rangle\langle 2k\rangle)
     \end{eqnarray*}
Since we are in fact only interested in knowing the head and the socle of $
{\bf E}_i M$ it is enough to consider 
\begin{eqnarray*}
&\HOM_{C^{i+1}}(\mV_{i+1}
     T_{i+1},\oplus_{k=0}^i C^{i+1}\langle -(n-i-1)\rangle\langle
     -i(n-i)\rangle\langle 2k\rangle)&\\
&=\HOM_{C^{i+1}}(C^{i+1},\oplus_{k=0}^i C^{i+1}\langle -(n-i-1)\rangle\langle
     -i(n-i)\rangle\langle 2k\rangle\langle (i+1)(n-i-1)\rangle).&
     \end{eqnarray*}
There are (up to scalars) unique morphisms of minimal (resp. maximal) degree namely of degree $s=-(n-i-1)-i(n-i)+(i+1)(n-i-1)=-i$ (and
$t=-i+2i+2(i+1)(n-i-1)=i+2(i+1)(n-i-1)$ respectively). Hence, the module  $
{\bf E}_i M={\bf E}_i\op{d}_i' M$ has minimal degree $-i$ and maximal degree
$t$. Therefore,  $\op{d_{i+1}}{\bf E}_i M$ has maximal degree $i$ and minimal
degree $-t$. By definition of $\op{d}_{i+1}'$, we get that $\op{d}_{i+1}'{\bf E}_i
M$ has maximal degree $i+2(i+1)(n-i-1)=t$ and minimal degree
$-t+2(i+1)(n-i-1)=-i$. This proves the first formula. The second follows then
by the adjointness properties from Proposition~\ref{all-adjoint} as follows:
The functor ${\bf E}_i$ has right adjoint ${\bf F}_{i+1}({\bf
  K}_{i+1})^{-1}\langle 1\rangle$ and the functor $\op{d}_{i+1}'{\bf E}_i\op{d}_i'$ has right adjoint 
\begin{eqnarray*}
\op{d}_{i}'{\bf K}_i\langle1\rangle{\bf F}_{i+1}\op{d}_{i+1}'
&\cong&\op{d}_{i}'{\bf F}_{i+1}\langle 2i-n+1\rangle\op{d}_{i+1}'\\
&\cong&\op{d}_{i}'{\bf F}_{i+1}\op{d}_{i+1}'\langle n-2i-1\rangle\\
&\cong&\op{d}_{i}'{\bf F}_{i+1}\op{d}_{i+1}'({\bf K}_{i+1})^{-1}
\langle1\rangle.
\end{eqnarray*}
Hence ${\bf F}_{i+1}\cong\op{d}_i'{\bf F}_{i+1}\op{d}_{i+1}'$, or
$\op{d}_i'{\bf F}_{i+1}\cong{\bf F}_{i+1}\op{d}_{i+1}'$. The last isomorphism
of the lemma follows from the isomorphisms 
${\bf K}_{i}\op{d}_{i}'={\bf K}_{i}\langle 2i(n-i)\rangle\op{d}_{i}
\cong\langle 2i(n-i)\rangle\op{d}_{i}{\bf K}_{i}^{-1}\cong\op{d}_{i}'{\bf K}_{i}^{-1}$.
\end{proof}

Let $\op{d}'=\oplus_{i=0}^n\op{d}_i':\oplus_{i=0}^n\op{gmod-}{{}_i A}$ be the
duality from above. For an endofunctor $F$ of $\oplus_{i=0}^n\op{gmod-}{{}_i A}$
let $\Psi (F)$ denote the functor $\op{d}'\,F\,\op{d}'$. The involution $\psi$ has the following functorial interpretation:
\begin{prop}
\label{invPsi}
  The functor $\Psi$ is an involution satisfying $\Psi({\bf E})\cong{\bf E}$,
  $\Psi({\bf F})\cong{\bf F}$, $\Psi({\bf K})\cong{\bf K}^{-1}$, and $\Psi(\langle
  k\rangle)\cong\langle -k\rangle$ for any $k\in \mZ$.
\end{prop}

\begin{proof}
  By definition, $\Psi$ is an involution satisfying $\Psi(\langle
  k\rangle)\cong\langle -k\rangle$. The rest follows from Lemma~\ref{duality}.
\end{proof}

\section{Schur-Weyl duality and special bases}
\label{SchurWeylduality}

Permuting the factors of the $\SL_2$-module $\overline{V}_1^{\otimes n}$ gives
rise to an additional $S_n$-module structure which commutes with the action of
the Lie algebra. In the quantised version we get an action of the Hecke
algebra corresponding to the symmetric group $S_n$. We would like to give a
categorical version of this bimodule $V_1^{\otimes n}$. 

\subsection{A categorical version of the Schur-Weyl duality}
Let $(W,S)$ be a Coxeter system. The corresponding Hecke algebra $\mathscr{H}(W,S)$
is the associative algebra (with $1$) over $\mZ[q,q^{-1}]$, the ring of Laurent polynomials
in one variable, with generators $H_s$, for $s\in S$ and relations
\begin{eqnarray*}
(H_s+q)(H_s-q^{-1})&=&0,\\
H_sH_tH_s\cdots H_t&=&H_tH_sH_t\cdots H_s,\quad\text{ if
$sts\cdots t=tst\cdots s$}\\
H_sH_tH_s\cdots H_s&=&H_tH_sH_t\cdots H_t,\quad\text{ if $sts\cdots s=tst\cdots t$}.
\end{eqnarray*}
In particular, if $x\in W$
with reduced expression $x=s_{i_1}\cdots s_{i_r}$ then
$H_x=H_{s_{i_1}}H_{s_{i_2}}\cdots H_{s_{i_r}}$ does not depend on the reduced
expression and $\{H_x\mid x\in W\}$ is a $\mZ[q,q^{-1}]$-basis of $\mathscr{H}(W,S)$.   
For any subset $S'$ of $S$ we
 get $W_{S'}\subset W$ and define
 $M^{S'}=\mathscr{H}(W,S)\otimes_{\mathscr{H}(W_{S'},S')}\mZ$, the corresponding
   permutation module. (Here, $H_s\in\mathscr(W_{S'}, S')$ acts on $\mZ$
   by multiplication with $q^{-1}$.) We denote by $\mathscr{H}$ the
   complexified Hecke algebra
   corresponding to the symmetric group $S_n$ and by
   $\cM^i=\mathbb{C}\otimes_\mZ M^i$ the complexified permutation
   module corresponding to $W=S_n$ and $W_{S'}=W_i$. Then the $\cM^i$ has a basis
   $\{M_x^i=1\otimes (H_x\otimes 1)\mid x\in
   W^i\}$. For more details we refer for example to \cite{SoKipp}. (Note that,
   beside the complexification, we work with left $\mathscr{H}$-modules, whereas the modules in
   \cite{SoKipp} are right $\mathscr{H}$-modules).  

In \cite{FKK}, the authors describe explicitly a well-known isomorphism of
$\mathscr{H}$-modules
\begin{eqnarray}
\label{KhKi}
  \alpha:\quad\bigoplus_{i=0}^n\cM^i&\cong& V_1^{\otimes n},\\
  M_x^i&\mapsto& v_{{\bf a}(x)}\nonumber
\end{eqnarray}
where ${\bf a}(x)$ is the $\{0,1\}$-sequence such that
$M(x\cdot\omega_i)=M({\bf a}(x))$. The algebra $\mathscr{H}$ acts on the right hand side via the so-called R-matrix. We refer to \cite[Proposition 2.1']{FKK} for details. (Note that
our $q$ is the $v$ there and our $H_i$ is $vT_i$ there.) 

On the other hand, there is an isomorphism of $\mQ(q)$-modules 
\begin{eqnarray}
\label{combinatorics}
  \beta\quad:\bigoplus_{i=0}^n\cM^i&\cong&\op{G}(\bigoplus_{i=0}^n\op{gmod-}{}_i A\big),\\
 M_{x}^i&\mapsto&1\otimes [\tilde M(x\cdot\omega_i)]\nonumber.
\end{eqnarray}
It induces an $\mathscr{H}$-action on the space on the right hand side. 

Our next task will be to ``categorify'' this action. For any simple reflection
$s\in W$, there is a {\it twisting functor} $T_s:\cO\rightarrow \cO$ which
preserves blocks, in particular induces $T_s:{}_\la\cO\rightarrow{}_\la\cO$ for any
integral dominant weight $\la$. These functors were studied for example in
\cite{AL}, \cite{AS}, \cite{KhM}. The most convenient description (for our
purposes) of these functors is given in \cite{KhM} in terms of partial
coapproximation: Let $M\in{}_\la\cO$ be projective. Let $M'\subset M$ be the
smallest submodule such that $M/M'$ has only composition factors of the form
$L(x\cdot\la)$, where $sx>x$. Then $M\mapsto M'$ defines a functor $T_s$ from the additive category of projective modules in
${}_\la\cO$ to ${}_\la\cO$. This functor extends in a unique way to a right exact
functor $T_s:{}_\la\cO\rightarrow{}_\la\cO$ (for details see \cite{KhM}). This
definition of $T_s$ has the advantage that it is immediately clear that this functor is
gradable. Moreover, we have the following:

\begin{prop}
\label{twist}
  For any simple reflection $s\in W$ and integral dominant weight $\la$, the twisting functor
  $T_s:{}_\la\cO\rightarrow{}_\la\cO$ is gradable. A graded lift is unique up to
  isomorphism and shift in the grading.  
\end{prop}

\begin{proof}
  We only have to show the uniqueness of a graded lift. Since $T_s$ is right
  exact, it is given by tensoring with some bimodule (see \eg \cite[2.2]{Bass}). By \cite[Lemma 2.5.3]{BGS} it is enough to show that $T_s$ is indecomposable. Let $G_s$ be
  the right adjoint functor of $T_s$. If $\la$ is
  regular then $G_sT_s\cong\ID$ on the additive category given by all
  projective modules in ${}_\la\cO$ (\cite[Corollary 4.2]{AS}). Since $T_s$ and $G_s$ commute
  with translation functors (\cite[Section 3]{AS}), the adjunction morphism
  defines an isomorphism $G_sT_s\cong\ID$ on the additive category given by all
  projective modules in $_\la\cO$ even for singular integral $\la$. 
  Hence $\END_\cO(T_sP)\cong\END_\cO(P)$ for any projective module $P\in{}_\la\cO$. In
  particular, $T_sP$ is indecomposable, if so is $P$. Assume now $T_s\cong
  F_1\oplus F_2$. For any indecomposable projective $P\in{}_\la\cO$ there exists
  $i(P)\in\{1,2\}$ such that $F_{i(P)}(P)=0$. Since $M(\la)$ is a submodule
  of any projective module $P$, we have $i(P)=i(M(\la))$ for any $P$. Hence
  $T_s$ is indecomposable when restricted to projective modules, hence also
  when considered as a functor on ${}_\la\cO$.      
\end{proof}

The next result lifts the action of the Hecke algebra to a
functorial action (this should be compared
with \cite[Proposition 1.1]{FKK} or \cite[page 86]{SoKipp}):

\begin{prop}
\label{SchurWeyl}
  Fix $i\in\{0,1,\ldots, n\}$. There are right exact functors ${\overline
  H}_j:{}_{i;n}\cO\rightarrow{}_{i;n}\cO$, $1\leq j<n$ satisfying the following
  properties 
  \begin{enumerate}[(a)]
  \item they are exact on the subcategory of modules having a filtration with
  subquotients isomorphic to Verma modules. 
 \item 
\label{twistb}
they have graded lifts ${\bf H}_j$ satisfying 
\pagebreak[0]
   \begin{align}
   \label{Rrel}
   \begin{split}
   &\big[{\bf H}_j(\tilde M(x\cdot\omega_i)\big]\\
   &=\begin{cases}
     [(\tilde M(s_jx\cdot\omega_i)]+(q^{-1}-q)[(\tilde M(x\cdot\omega_i)]&
\text{ if
     $s_jx<x$, $s_jx\in W^i$},\\
[(\tilde M(s_jx\cdot\omega_i))]&\text{ if $s_jx>x$, $s_jx\in W^i$},\\ 
q^{-1}[(\tilde M(x\cdot\omega_i))]&\text{ if $s_jx\notin W^i$.}
   \end{cases}
     \end{split}
   \end{align} 
\item  Let $w_0=s_{i_1}s_{i_2}\cdots s_{i_r}$ be a reduced expression for the
  longest element in $S_n$ and ${\bf H}_{w_o}={\bf H}_{i_1}{\bf H}_{i_2}\cdots
  {\bf H}_{i_r}$ the corresponding composition of functors. Then ${\bf
  H}_{w_o}$ is exact on modules with Verma flag and
\begin{alignat}{6}
\label{twist1}
 &{\bf H}_{w_o}\tilde M(x\cdot\la)&\quad&\cong&\quad&(\op{d} \tilde M(w_0x\cdot\la))\,\langle
  -l(w_0^i)\rangle,\\
\label{twist2}
&{\bf H}_{w_o}\tilde P(x\cdot\la)&&\cong&&\tilde T(w_0x\cdot\la)\,\langle -l(w_0^i)\rangle,\\
\label{twist3}
&{\bf H}_{w_o}\tilde T(x\cdot\la)&&\cong&&\tilde I(w_0x\cdot\la)\,\langle
  -l(w_0^i)\rangle.
\end{alignat}
where $\tilde T(x\cdot\la)$ denotes the graded lift of the tilting module
$T(x\cdot\la)$ such that $\tilde M(x\cdot\la)$ occurs as a submodule in a
Verma flag and $\tilde I(x\cdot\la)=\op{d}\tilde P(x\cdot\la)$ is the injective hull 
of the simple module $\tilde L(x\cdot\la)$.    
 \end{enumerate}
\end{prop}

\begin{proof}
 We claim that if we forget the grading, these functors are the twisting
 functors $T_{s_i}$. They are right exact by definition, are exact when restricted to the subcategory of modules having a
 Verma flag (\cite[Theorem 2.2]{AS}) and satisfy the relations \eqref{Rrel} if
 we forget the grading (\cite[6.5 and 6.6 or Lemma 2.1]{AL}). 

We know that these functors are
 gradable and indecomposable when restricted to an integral block (Proposition~\ref{twist}). Therefore, we just need a ``correct'' lift of these
 functors. We choose a graded lift $\tilde T_s$ of
 $T_s:{}_\la\cO\rightarrow{}_\la\cO$ such that 
 \begin{equation}
 \label{eq:lift}
   \begin{split}  
\tilde T_s \tilde
 M(\la)&\cong\tilde M(s\cdot\la)\quad\quad\quad\text{ if } s\in W^\la,\\
\tilde T_s \tilde
 M(\la)&\cong \tilde M(s\cdot\la)\langle -1\rangle\quad\text{ otherwise}.    
 \end{split}
\end{equation}

Let $sx>x$ then $T_sM(x\cdot\la)\cong M(sx\cdot\la)$ (\cite[Lemma 2.1]{AL}). Hence $\tilde
T_s\tilde M(x\cdot\la)\cong \tilde M(sx\cdot\la)\langle k\rangle$ for some
$k\in\mZ$. On the other hand we have an inclusion 
\begin{eqnarray}
  \label{eq:incl}
\tilde M(x\cdot\la)\langle
l(x)\rangle \hookrightarrow\tilde M(\la)
\end{eqnarray}
 for any $x\in W^\la$ (for example by \cite[Proposition 3.11.6]{BGS}).  
 \begin{itemize}
 \item Assume $sx>x$, $sx\in W^\la$. Then $\tilde M(sx\cdot\la)\langle
   l(x)+k\rangle\cong\tilde T_s\tilde M(x\cdot\la)\langle
l(x)\rangle\hookrightarrow\tilde T_s\tilde M(\la)$. From \eqref{eq:lift} it
   follows that  $\tilde M(sx\cdot\la)\langle
   l(x)+k\rangle\hookrightarrow\tilde M(\la)\langle -1\rangle$, hence
   $l(x)+k+1=l(sx)=l(x)+1$. That means $k=0$. We get $\tilde
T_s\tilde M(x\cdot\la)\cong \tilde M(sx\cdot\la)$.
\item Assume $sx\not\in W^\la$ (in particular $sx>x$), hence $\tilde
T_s\tilde M(x\cdot\la)\cong \tilde M(sx\cdot\la)\langle k\rangle$ for some
$k\in\mZ$ (\cite[6.5 and 6.6 or Lemma 2.1]{AL}).  Then 
\begin{equation*}
\tilde M(x\cdot\la)\langle
   l(x)+k\rangle\cong\tilde T_s\tilde M(x\cdot\la)\langle
l(x)\rangle\hookrightarrow\tilde T_s\tilde M(\la). 
\end{equation*}
From \eqref{eq:lift} it
   follows  $\tilde M(x\cdot\la)\langle
   l(x)+k\rangle\hookrightarrow\tilde M(\la)\langle -1\rangle$, hence $l(x)+k=l(x)-1$. That means $k=-1$. We get $\tilde
T_s\tilde M(x\cdot\la)\cong \tilde M(x\cdot\la)\langle -1\rangle$.
\item Assume $sx<x$, $sx\in W^\la$. From the
 translation principle it follows that there is a
 unique non-split extension 
 \begin{eqnarray*}
   0\rightarrow M(y\cdot\la)\rightarrow M\rightarrow M(sy\cdot\la)\rightarrow 0
 \end{eqnarray*}
whenever $y$, $sy\in W^\la$ and $sy>y$. (To see this one could first consider
the case where $\la$ is regular. If we write $sy=yt$ for some simple
reflection $t$ then the statement becomes familiar. The general statement
follows then by translation). 
If $\la$ is regular, the main result
of \cite{KhM} says that $T_s$ is adjoint to Joseph's completion functor (see
\cite{Jo1}), in particular $T_sM(z\cdot\la)$ is the cokernel of the inclusion
$M(z\cdot\la)\rightarrow M$ for $z\in\{y, sy\}$. In the graded picture we
have a unique non-split extension 
 \begin{eqnarray*}
   0\rightarrow M(y\cdot\la)\langle 1\rangle\rightarrow M\rightarrow M(sy\cdot\la)\rightarrow 0
 \end{eqnarray*}
whenever $y\in W$, $sy>x$ and then $\tilde T_s\tilde M(y\cdot\la)$ is
the cokernel of the inclusion $\tilde M(y\cdot\la)\langle 1\rangle\rightarrow
M$ if $y<sy$, and $\tilde T_s\tilde M(sy\cdot\la)$ is
the cokernel of the inclusion $\tilde M(sy\cdot\la)\langle 1\rangle\rightarrow
M\langle -1\rangle$ if $y>sy$ (compare \cite[Theorem 5.3, Theorem
3.6]{Stgrad}). We get the following formula for regular integral $\lambda$:
\begin{eqnarray}
\label{eqreg}
\big[\tilde T_s(\tilde M(x\cdot\la))\big]
   &=&[(\tilde M(sx\cdot\la))]+(q^{-1}-q)[(\tilde M(x\cdot\la))].
\end{eqnarray} 
To see this we just calculate 
$$[M\langle-1\rangle]-[\tilde M(x\cdot\la)\langle 1\rangle]=[\tilde
M(sx\cdot\la)]+[\tilde M(x\cdot\la)\langle -1\rangle]-[\tilde
M(x\cdot\la)\langle 1\rangle],$$ and the formula follows. 

To get the result for singular blocks we use translation functors.  Let $\la$,
 $\mu$ be dominant
 integral weights with $\mu$ dominant. Then the translation functors $\mathds{T}_\la^\mu$ and $\mathds{T}_\mu^\la$ are gradable. This follows from
 Proposition~\ref{EFandV} as follows. Since the functor $\mV_\mu$ induces an
 isomorphism $\HOM_\mg(P_\mu, \mathds{T}_\la^\mu P_\la)\cong\HOM_{C^\mu}(\mV_\mu P_\mu,
 \mV_\mu \mathds{T}_\la^\mu P_\la)$ (by the Struktursatz of \cite{Sperv}), a
 graded lift $\tilde{\mathds{T}}^\la_\mu$ of $\mathds{T}_\mu^\la$ is given by tensoring with the $(\END_{C^\mu}(\mV_\mu P_\mu),\END_{C^\la}(\mV_\la P_\la))$-bimodule $\HOM_{C^\la\gMOD}(\mV_\la P_\la,
 \RES^\mu_\la\mV P_\mu)$. We have $\tilde{\mathds{T}}^\la_\mu M(\mu)\cong M(\la)\langle
 k\rangle$  for some $k\in\mZ$. The inclusion~\eqref{eq:incl} implies $\tilde{\mathds{T}}^\la_\mu M(x\cdot\mu)\cong M(x\cdot\la)\langle
 k\rangle$ for any $x\in W^\la$. Without loss of generalities we may assume
 $k=0$. We have the following equalities
 \begin{eqnarray}
   [\tilde T_s\tilde M(x\cdot\la)]&=&[\tilde T_s\tilde{\mathds{T}}_\mu^\la\tilde
   M(x\cdot\mu)]\nonumber\\
&=&[\tilde{\mathds{T}}_\mu^\la\tilde T_s\tilde M(x\cdot\mu)]\label{eq1}\\
&=&[\tilde{\mathds{T}}_\mu^\la\tilde M(sx\cdot\mu)]+(q+q^{-1})[\tilde{\mathds{T}}_\mu^\la\tilde
   M(x\cdot\mu)]\label{eq2}\\
&=&[\tilde M(sx\cdot\la)]+(q+q^{-1})[\tilde
   M(x\cdot\la)].\nonumber
 \end{eqnarray}
The first and the last equality follow from the definitions. To see the equality
\eqref{eq1} observe that twisting functors and translation functors commute
(see \cite[Section 3]{AS}). With standard arguments one can check that $T_s \mathds{T}_\mu^\la$
is indecomposable, hence $\tilde T_s\tilde{\mathds{T}}_\mu^\la$ is isomorphic to $\tilde{\mathds{T}}_\mu^\la\tilde T_s$ up to a
shift in the grading. However, they agree on $\tilde M(\mu)$, hence they are
isomorphic. Finally \eqref{eq2} follows from the equation~\eqref{eqreg}. 
\end{itemize}   
This finishes the proof of part~\eqref{twistb}  of the proposition. After forgetting the grading we have ${\bf H}_{w_o}M(x\cdot\la)\cong\op{d}
M(w_0x\cdot\la)$ (this is \cite[(2.3) and Theorem 2.3]{AS} for the regular case, the general case
follows easily by translation). Since the left derived functor of ${\bf
  H}_{w_o}$ defines an equivalence of derived categories (\cite[Corollary
4.2]{AS}), by general arguments, we get isomorphisms
of modules  
\begin{alignat*}{6}
  &{\bf H}_{w_o}M(x\cdot\la)&\quad&\cong&\quad&\op{d}M(w_0x\cdot\la),\\
  &{\bf H}_{w_o}P(x\cdot\la)&&\cong&&T(w_0\cdot\la),\\
  &{\bf H}_{w_o}T(x\cdot\la)&&\cong&&I(w_0\cdot\la).
\end{alignat*}
For details we refer for example to \cite[Proposition~4.2]{GGOR}. In
particular, tilting modules and injective modules are gradable. We are left with checking the
graded version. Since all modules involved are indecomposable, we have 
\begin{alignat*}{6}
  &{\bf H}_{w_o}\tilde M(x\cdot\la)&\quad&\cong&\quad&(\op{d} \tilde M(w_0x\cdot\la))\langle k_x^1\rangle\\
  &{\bf H}_{w_o}\tilde P(x\cdot\la)&\quad&\cong&\quad&\tilde T(w_0\cdot\la)\langle k_x^2\rangle,\\
  &{\bf H}_{w_o}\tilde T(x\cdot\la)&\quad&\cong&\quad&\tilde I(w_0\cdot\la)\langle k_x^3\rangle.
\end{alignat*}
for some $k_x^i\in\mZ$. Set $s=-l(w_o^i)$. We claim that $\tilde
L(w_ox\cdot\la)\langle s\rangle$ occurs as a
composition factor in ${\bf H}_{w_o}\tilde M(x\cdot\la)$. This is clear from
the formulas~\eqref{Rrel}, hence $k_x^1=s$. The inclusion $\tilde
M(x\cdot\la)\hookrightarrow \tilde T(x\cdot\la)$ gives rise to an inclusion
$(\op{d}\tilde M(w_0x\cdot\la))\langle s\rangle\hookrightarrow {\bf H}_{w_o} T(x\cdot\la)$. This
implies $k_x^3=s$. The surjection $\tilde P(x\cdot\la)\rightarrow\tilde
M(x\cdot\la)$ gives rise to a surjection ${\bf H}_{w_o}\tilde P(x\cdot\la)\rightarrow(\op{d}\tilde
M(w_0x\cdot\la))\langle s\rangle$. Hence $\tilde L(w_0x\cdot\la)\langle
s\rangle$ occurs as a composition factor in ${\bf H}_{w_o}\tilde
P(x\cdot\la)$. This implies $k_x^2=s$. The proposition follows. 
\end{proof}

In the previous Proposition~\ref{SchurWeyl}, to categorify the action of the
Hecke algebra we restricted the (graded lifts) of the twisting functors
$\ov{H}_j$ to the category of modules with Verma flag, to force them to be
exact. J. Sussan studied a categorification of the Hecke algebra action
by considering the derived functors associated to the twisting functors
(\cite{Josh}). 

\subsection{The canonical, standard and dual canonical bases}

We will now combine the three pictures: the Hecke module, the Grothendieck
group of the graded version of certain blocks of $\cO$ and the $U_q(\SL_2)$
module $V_1^{\otimes n}$. We consider $V_1^{\otimes n}$ as a $U$-module via
the comultiplication ${\blacktriangle}$. 

We first have to introduce a bilinear form on $V_1^{\otimes n}$ for $n\geq 2$. There is a nondegenerate bilinear form 
$<,>$ on $V_1^{\otimes 2}$ defined by $<v_i\otimes v_j, v^k\otimes
v^l>=\delta_i^l\delta_j^k$. It satisfies 
\begin{equation*}
<\blacktriangle(x)(v_i\otimes
v_j), v^k\otimes v^l>=<v_i\otimes
v_j, \blacktriangle'(\sigma(x)) (v^k\otimes v^l)>,  
\end{equation*}
where $\blacktriangle'=\psi\otimes\psi\circ\blacktriangle\circ\psi$, explicitly \begin{equation*}
\blacktriangle'(E)=1\otimes E + E\otimes K, \quad\blacktriangle'(F)=K^{-1}\otimes F+F\otimes 1,
\quad \blacktriangle'(K^{\mp 1})=K^{\mp1}\otimes K^{\mp1}.  
\end{equation*}
Let  $(V_1^{\otimes n})'$ denote the $U_q(\mathfrak{sl}_2)$-module $V_1^{\otimes n}$
but with comultiplication $\blacktriangle'$. 
The form above can be extended to a bilinear form $< ,>:V^{\otimes n}\times (V^{\otimes n})'\rightarrow \mC$ by
putting   
\begin{eqnarray}
\label{form}
  \langle v_{\bf a}, v^{\bf b}\rangle= \prod_{i=1}^n\delta_{a_i,b_{n-i+1}}
\end{eqnarray}
Then it satisfies 
$<xu,v> = <u,\sigma(x)v>$ for any $x\in U_q(\mf{sl}_2)$ and 
$u,v\in V_1^{\otimes n}$. 

The module  $V_1^{\otimes n}$ has two distinguished $\Q(q)$-bases, namely
\begin{itemize}
\item the {\it standard basis}
  $\{v_{\bf a}=v_{a_1}\otimes v_{a_{2}}\otimes\cdots\otimes v_{a_n}\mid
  a_j\in\{0,1\}\}$,
\item the {\it canonical basis}  $\{v_{\bf a}^\Diamond=v_{a_1}\Diamond v_{a_2}\Diamond\cdots\Diamond
  v_{a_n}\mid a_j\in\{0,1\}\}$.
\end{itemize}

There are also two distinguished basis in the space $(V_1^{\otimes n})'$
namely 
\begin{itemize}
\item  the {\it dual standard basis}
  $\{v^{\bf a}=v^{a_1}\otimes v^{a_{2}}\otimes\cdots\otimes v^{a_n}\mid
  a_j\in\{0,1\}\}$,
\item the {\it dual canonical basis} $\{v^{\bf a}_\heartsuit=v^{a_1}\heartsuit v^{a_{2}}\heartsuit\cdots\heartsuit
  v^{a_n}\mid a_j\in\{0,1\}\}$.
\end{itemize}
The canonical and dual canonical bases were defined by Lusztig and Kashiwara (\cite{Lu1}, \cite{Lucantensor}, \cite{Kas}). Lusztig (\cite[Chapter 27]{Lubook}) defined a certain semilinear involution $\Psi$ on $V_1^{\otimes
n}$ which determines the canonical basis uniquely by the following two properties
\begin{enumerate}[(i)]
\item $\Psi(v_{\bf a}^\Diamond)=v_{\bf a}^\Diamond$.
\item $v_{\bf a}^\Diamond-v_{\bf a}\in\sum_{{\bf b}\not={\bf a}}q^{-1}\mZ[q^{-1}]v_{\bf b}$.
\end{enumerate}

Given the canonical basis, the dual canonical basis is defined by 
\begin{eqnarray}
\label{form2}
  \langle v_{\bf a}^\Diamond, v^{\bf b}_\heartsuit\rangle= \prod_{i=1}^n\delta_{a_i,b_{n-i+1}}.
\end{eqnarray}
(For an explicit graphical description of these bases we refer to
\cite{FK}.)\\

On the other hand, the permutation module $\cM^i$ has also several
distinguished $\mC[q,q^{-1}]$-bases, namely 

\begin{itemize}
\item the {\it standard basis}
  $\{M_x^i=1\otimes H_x\otimes 1\mid x\in W^i\}$,
\item the (positive) {\it self-dual basis}  $\{\underline {M}_x^i\mid x\in
  W^i\}$, 
\item the (negative) {\it self-dual basis}  $\{\underline{\tilde{M}}_x^{i}\mid
  x\in W^i\}$, 
\item the {\it ``twisted'' standard basis} $\{(M_x^i)^{Twist}:=q^{l(w_0^i)}H_{w_0}M_x^i\mid x\in W^i\}$,
\item the {\it ``twisted'' positive self-dual basis} $$\{(\underline M_x^i)^{Twist}:=q^{l(w_0^i)}H_{w_0}\underline
  {M}_x^i\mid x\in W^i\},$$ 
\item the {\it ``twisted'' negative self-dual basis}
  $$\{(\underline{\tilde{M}}_x^{i})^{Twist}:=q^{l(w_0^i)}H_{w_0}\underline{\tilde{M}}_x^{i}\mid
  x\in W^i\}.$$ 
\end{itemize}
These bases were defined by Kazhdan, Lusztig and Deodhar (see
  \cite{KLHecke}, \cite{Deodhar}). We use here the notation from \cite{SoKipp},
  except that we  have the upper
  index $i$ to indicate that $M_x^i\in \cM_x^i$ and we also use $q$ instead of
  $v$. The bases can be characterised as follows (\cite{KLHecke} in the notation of
\cite{SoKipp}): Let
$\Psi:\mathscr{H}\rightarrow\mathscr{H}$ be the $\mZ$-linear involution given
  by $H_x\mapsto (H_{x^{-1}})^{-1}$,
$q\mapsto q^{-1}$. It induces an involution on any $\cM^i$. Then the 
$\underline {M}_x^i$ are uniquely defined by

\begin{enumerate} [(i)]
\item $\Psi(\underline {M}_x^i)= \underline {M}_x^i$,
\item $\underline {M}_x^i- {M}_x^i\in \sum_{y\not=x} v\mZ[v]M_y^i$. 
\end{enumerate}

The basis elements $\underline{\tilde{M}}_x^i$ are characterised by

\begin{enumerate} [(i)]
\item $\Psi(\underline{\tilde{M}}_x^i)=\underline{\tilde{M}}_x^i$,
\item $\underline{\tilde{M}}_x^i- {M}_x^i\in \sum_{y\not=x} v^{-1}\mZ[v^{-1}]M_y^i$. 
\end{enumerate}
Note that what we call the ``twisted'' bases are in
fact bases, since $H_{w_0}$ is invertible in $\mathscr{H}$. 

The following theorem gives a categorical interpretation of all these bases:

\begin{theorem}
\label{bases2}
  \begin{enumerate}[(a)]
  \item \label{bases2.1}
  There is an isomorphism of $U_q(\SL_2)$-modules 
  \begin{eqnarray*}
    \Phi:\quad\op{G}(\bigoplus_{i=0}^{n} \op{gmod-}{}_i A)&\cong& V_1^{\otimes n}\\
    1\otimes\big[\tilde M({\bf a})\langle i\rangle\big]&\mapsto&q^iv_{\bf a}=q^i(v_{a_1}\otimes v_{a_2}\otimes\cdots\otimes
    v_{a_n}),
\end{eqnarray*}
where the $U_q(\mathfrak{sl}_2)$-structure on the left hand side is
induced by the functors ${\bf E}$, ${\bf F}$ and ${\bf K}$.
\item   \label{bases2.1b}
  There is an isomorphism of $U_q(\SL_2)$-modules 
  \begin{eqnarray*}
\Phi':{\bf{G}}\big(\bigoplus_{i=0}^{n}{\op{gmod}-{}_i A}\big)&\cong&({V}_1^{\otimes n})'\\
 1\otimes\big[\tilde\nabla({\bf a})\langle i\rangle\big]&\mapsto&q^{i}v^{\bf a}=q^iv^{a_1}\otimes\cdots\otimes v^{a_n}. 
\end{eqnarray*}
where the $U_q(\mathfrak{sl}_2)$-structure on the left hand side is
induced by the functors ${\bf E}'=\oplus_{i=0}^n{\bf E}_i'$, ${\bf
  F}'=\oplus_{i=0}^n{\bf F}_i'$ and ${\bf K}$, where 
\begin{eqnarray*}
{\bf E}_i'&=&\langle-2(i+1)(n-i-1)\rangle{\bf E}_i\langle 2i(n-1)\rangle\\
{\bf F}_i'&=&\langle-2(i-1)(n-i+1)\rangle{\bf F}_i\langle 2i(n-1)\rangle. 
\end{eqnarray*}
\item \label{bases2.2} 
The isomorphism $\beta$ defines bijections:
\begin{equation*}
 \begin{array}[htu]{cccc}
\mbox{positive self-dual
  basis}&\leftrightarrow&\mbox{ standard lifts of indec. projectives}\\
\underline M_{x}^i&\mapsto&[\tilde P(x\cdot\omega_i)],\\
\\
\mbox{standard basis}&\leftrightarrow&\mbox{standard lifts
    of Verma modules}\\
\underline M_{x}^i&\mapsto&[\tilde M(x\cdot\omega_i)].\\
\end{array}
\end{equation*}
\item \label{bases2.3} The bilinear form $\langle\quad,\quad\rangle$ can be realized as follows
  \begin{eqnarray*}
  <1\otimes [M],1\otimes [N]>=\sum_{i}\sum_{j}(-1)^j\DIM\HOM^j_\mathcal{D}(\op{d}N, {\bf H}M\langle -i\rangle)[j])\; q^i,  
  \end{eqnarray*}
where $\mathcal{D}$ denotes the bounded derived category of $\bigoplus
\op{gmod-}{}_i A$ with shift functors $[j]$, and
${\bf H}$ denotes the derived functor of the twisting functor ${\bf H}_{w_0}$
from Proposition~\ref{SchurWeyl}, but shifted in the grading such that the standard lifts
of Verma modules are sent to standard lifts of dual Verma modules. 

 \item \label{bases2.4} We have $\Phi=\alpha\circ\beta^{-1}$ (from \eqref{KhKi} and \eqref{combinatorics})
 and these isomorphisms, together with $\Phi'$, define bijections:\\
\begin{equation*}
\boxed{
 \begin{array}[htu]{ccccc}
\\
\parbox[t]{3.5cm}{twisted positive self-dual basis}&\leftrightarrow&\parbox[t]{3cm}{standard lifts
    of indec. tilting modules}&\leftrightarrow&\parbox[t]{3cm}{\text canonical
    basis}
\vspace{0.1cm}\\
(\underline M_{x}^i)^{Twist}&\mapsto&1\otimes[\tilde
    T(w_0x\cdot\omega_i)],
\vspace{0.1cm}\\
&&1\otimes[\tilde T({\bf a})]&\mapsto&
    v_{\bf a}^{\Diamond}\\
\\
\hline
\\
\parbox[t]{3.5cm}{standard basis}&\leftrightarrow&\parbox[t]{2.4cm}{standard lifts
    of Verma modules}&\leftrightarrow&\parbox[t]{3cm}{\text standard basis}
\vspace{0.1cm}\\
M_x^i&\mapsto&1\otimes[\tilde M(x\cdot\omega_i)],
\vspace{0.1cm}\\
&&1\otimes[\tilde M({\bf a})]&\mapsto&
    v_{\bf a}\\
\\
\hline
\\
\parbox[t]{3.5cm}{twisted standard
    basis}&\leftrightarrow&\parbox[t]{3.1cm}{standard lifts of dual Ver\-ma mo\-d\-u\-l\-es}&\leftrightarrow&\parbox[t]{3cm}{\text dual standard
    basis}
\vspace{0.1cm}\\
(M_{x}^i)^{Twist}&\mapsto&1\otimes [\tilde \nabla(w_0x\cdot\omega_i)],
\vspace{0.1cm}\\
&&1\otimes [\tilde\nabla({\bf a})]&\mapsto&
    v^{\bf a}\\
\\
\hline
\\
\parbox[t]{3.5cm}{\text twisted negative self-dual
    basis}&\leftrightarrow&\parbox[t]{2.4cm}{standard lifts of simple
    modules}&\leftrightarrow&\parbox[t]{3cm}{\text dual ca\-no\-ni\-cal basis}
\vspace{0.1cm}\\
\left(\tilde{\underline{M}}_x^i\right)^{Twist}&\mapsto&1\otimes
    [\tilde{L}(w_0x\cdot\omega_i)],
\vspace{0.1cm}\\
&&1\otimes [\tilde L({\bf a})]&\mapsto&
    v_{\bf a}^{\heartsuit}\\
\\ 
\end{array}
}
\end{equation*}
\item \label{bases2.5}
There is an isomorphism of $\mQ(q)$-modules 
\begin{eqnarray*}
  \gamma\quad:\bigoplus_{i=0}^n\cM^i&\cong&\op{G}(\bigoplus_{i=0}^n\op{gmod-}{}_i A\big),\\
 M_{x}^i&\mapsto&1\otimes [\tilde\nabla (w_0x\cdot\omega_i)].
\end{eqnarray*}
Under this isomorphism, the negative self-dual basis corresponds to standard
lifts of simple modules and the (twisted) positive self-dual basis corresponds to
standard lifts of tilting (resp. injective) modules; more precisely
\begin{eqnarray*}
  \begin{array}[th]{lcl}
\tilde{\underline M}_{x}^i&\mapsto&1\otimes[\tilde L(w_0x\cdot\omega_i)],
\vspace{0.1cm}\\
{\underline M}_{x}^i&\mapsto&1\otimes[\tilde T(w_0x\cdot\omega_i)],
\vspace{0.1cm}\\
({\underline M}_{x}^i)^{Twist}&\mapsto&1\otimes[\tilde I(x\cdot\omega_i)].
\end{array}
\end{eqnarray*}
\end{enumerate}
\end{theorem}

Before we proof the theorem we give an 
\begin{exs}{\rm
  Consider $\mathfrak{gl}_2$ (\ie $n=2$). Let $i=1$. Then $\op{gmod}-{_iA}$ is
  equivalent to the graded version of the principal block of
  $\cO(\mathfrak{gl}_2)$. Consider the Hecke module
  $\cM^1$. To avoid too many indices we will omit the sup-index $i$. The module
  $\cM=\cM^1$ has the standard basis given by the elements $M_e$ and $M_s$,
  where $s$ is the (only) simple reflection in the Weyl group. One can easily
  calculate the distinguished bases in the Hecke module $\cM$; together with \cite[Theorem 3.11.4 (ii)]{BGS} and \cite[Section 1.5]{FK} we get the following:
  \begin{itemize}
  \item the twisted standard basis is given by
  $(M_e)^{Twist}=M_s$ and $(M_s)^{Twist}=M_e+(q^{-1}-q)M_s$. The corresponding
  equations in the Grothendieck group are $\big[\op{d}\tilde\Delta(s\cdot
  0)\big]=\big[\tilde\Delta(s\cdot0)\big]$ and 
  $\big[\op{d}\tilde\Delta(0)\big]=\big[\tilde\Delta(0)\big]+\big[\tilde\Delta(s\cdot0)\langle
  -1\rangle\big]-\big[\tilde\Delta(s\cdot0)\langle 1\rangle\big]$.  
  \item the positive self-dual basis is given by $\underline M_e=M_e$ and
  $\underline M_s=M_s+qM_e$. These equations correspond to the equalities
  $\big[\tilde P(0)\big]=\big[\tilde\Delta(0)\big]$, and $\big[\tilde
  P(s)\big]=\big[\tilde\Delta(s\cdot0)\big]+\big[\tilde\Delta(e)\langle1\rangle\big]$ via the
  isomorphism $\beta$. 

  Under the isomorphism $\gamma$ the equations above correspond to the equalities
  $\big[\tilde T(s\cdot0)\big]=\big[\tilde\nabla(s\cdot0)\big]$ and $\big[\tilde T(0)\big]=\big[\tilde\nabla(0)\big]+q\big[\tilde\nabla(s\cdot0)\big]$.

  \item the twisted positive self-dual basis is given by $(\underline
  M_e)^{Twist}=M_s$ and $(\underline
  M_s)^{Twist}=(M_s)^{Twist}+q(M_e)^{Twist}=M_e+q^{-1}M_s$. These equations
  correspond to $\big[\tilde T(s\cdot0)\big]=\big[\tilde\Delta(s\cdot0)\big]$ and $\big[\tilde
  T(0)\big]=\big[\tilde\Delta(0)\big]+\big[\tilde\Delta(s\cdot0)\langle-1\rangle\big]$. Or,
  equivalently, to the equations  $\big[\tilde T((0,1))\big]=\big[\tilde\Delta((0,1))\big]$ and $\big[\tilde
  T((1,0))\big]=\big[\tilde\Delta((1,0))\big]+\big[\tilde\Delta((0,1))\langle-1\rangle\big]$. On the
  other hand we have $v_0\diamond v_1=v_0\otimes v_1$ and
  $v_1\diamond v_0=v_1\otimes v_0+q^{-1}v_0\otimes v_1$. Under the isomorphism
  $\gamma$ the twisted positive self-dual basis corresponds to the basis given
  by the standard lifts of the indecomposable injective modules, the
  corresponding relations are $\big[\tilde I(0)\big]=\big[\tilde\nabla(0)\big]$ and 
$\big[\tilde I(s\cdot0)\big]=\big[\tilde\nabla(s\cdot0)\big]+\big[\tilde\nabla(0)\langle-1\rangle\big]$.

  \item the negative self-dual basis is given by $\tilde{\underline M_e}=M_e$ and
  $\tilde{\underline{M}}_s=M_s-q^{-1}M_e$. Via the isomorphism $\gamma$, these
  equations become $\big[\tilde L(s\cdot0)\big]=\big[\op{d}\tilde\Delta(s\cdot0)\big]$ and $\big[\tilde L(0)\big]=\big[\op{d}\tilde\Delta(0)\big]-\big[
  \op{d}\tilde\Delta(s\cdot0)\langle -1\rangle\big]$, or to the equations  $\big[\tilde L((0,1))\big]=\big[\op{d}\tilde\Delta((0,1))\big]$ and $\big[\tilde L((1,0))\big]=\big[\op{d}\tilde\Delta((1,0))\big]-\big[
  \op{d}\tilde\Delta((0,1))\langle -1\rangle\big]$.
These formulas correspond to the
  following expressions of the dual canonical basis in terms of the dual
  basis: $v^0\heartsuit v^1=v^0\otimes v^1$ and $v^1\heartsuit v^0=v^1\otimes
  v^0-q^{-1}v^0\otimes v^1$.   
 \item the twisted negative self-dual basis is given by $(\tilde{\underline{
 M}}_e)^{Twist}=M_s$ and $(\tilde{\underline{M}}_s)^{Twist}=M_e-qM_s$. The
 corresponding equalities are $\big[\tilde L(s\cdot0)\big]=\big[\tilde\Delta(s\cdot0)\big]$
 and $\big[\tilde L(0)\big]=\big[\tilde\Delta(0)\big]-\big[\tilde\Delta(s\cdot0)\langle 1\rangle\big]$.  \end{itemize}
}
\end{exs}

We still have to proof the Theorem:

\begin{proof}[Proof of Theorem~\ref{bases2}]
 We have $\Phi=\alpha\circ\beta^{-1}$, since they agree by definition on
  standard modules. The existence of the isomorphism $\overline\Phi$ in
  \eqref{phi} implies that $\Phi$ is an isomorphism of $\mQ(q)$-modules. We have to show that it
  is a $U_q(\SL_2)$-morphism. Let $M=\tilde M({\bf a})\in\op{gmod-}{}_i A$.  From the definitions we get 
$\Phi({\bf K}[M])=\Phi([M\langle 2i-n\rangle])=q^{2i-n}\Phi([M])$. On the other
  hand $K v_{\bf a}=q^mv_{\bf a}$, where $m$ is the number of ones minus the
  numbers of zeros occurring in ${\bf a}$, hence $m=i-(n-i)=2i-n$. We get 
  $\Phi({\bf K}^{\pm 1}[M])=K^{\pm 1}\Phi([M])$ for any
  $M\in\oplus_{i=0}^n\op{gmod-}{}_i A$. \\
  Note that $\cE M(1^i0^{n-i})$ is
  projective. From \cite[Proposition 6]{BFK} we get that $P(1^i0^{n-i-1}1)$ is a direct
  summand and that $M(1^i0^{n-i})$ occurs with multiplicity one in any Verma
  flag. Hence we get in fact $\cE M(1^i0^{n-i})\cong P(1^i0^{n-i-1}1)$, since $M(1^i0^{n-i})$ would occur in a Verma flag of any other direct
  summand. Because of the indecomposability of $P(1^i0^{n-i-1}1)$ we get ${\bf E} \tilde M(1^i0^{n-i})\cong \tilde
  P(1^i0^{n-i-1}1)\langle k\rangle$ for some $k\in \mZ$ (\cite[Lemma 2.5.3]{BGS}). To determine $k$ we
  calculate $j$ such that
  $$\HOM_{\op{gmod-}{}_{i+1}A}(\tilde M(1^{i+1}0^{n-i-1})\langle
  j\rangle,{\bf E}\tilde M(1^i0^{n-i}))\not=0.$$
 (Since the homomorphism space
  in question is one-dimensional, the number $j$ is well-defined.) We set
  $M=M(1^i0^{n-i})$ and $N=M(1^{i+1}0^{n-i-1})$ From our
  definitions we get 
  \begin{eqnarray*}
&&   \HOM_{\op{gmod-}{}_{i+1}A}\big(\tilde M(1^{i+1}0^{n-i-1})\langle
  j\rangle,{\bf E}\tilde M(1^i0^{n-i})\big)\\
&=&
\HOM_{\op{gmod-}{}_{i+1}A}(X,Y),
\end{eqnarray*} 
where  
$Y=\HOM_C\big(\mV_{i+1} P_{i+1}, \RES_{i+1}^{i,i+1}
  C^{i,i+1}\otimes_{C^i}\mV_i M\langle -n+i+1\rangle\big)$ and $X=\HOM_{C^{i+1}}(\mV_{i+1}P_{i+1},\mV_{i+1} N)\langle
  j\rangle$. Then 
\begin{eqnarray*}
&&\HOM_{\op{gmod-}{}_{i+1}A}(X,Y)\\
&=&\HOM_{C^{i+1}}\big(\mV_{i+1}N\langle
  j\rangle, \RES_{i+1}^{i,i+1} C^{i,i+1}\otimes_{C^i}\mV_i M\big)\langle -n+i+1\rangle)\\
&=&\HOM_{C^{i+1}}\big(\mQ\langle
  j\rangle, \RES_{i+1}^{i,i+1} C^{i,i+1}\otimes_{C^i}\mQ\langle -n+i+1\rangle\big).\
  \end{eqnarray*}
Now Lemma~\ref{free} tells us that $j$ has to be the highest nonzero degree occurring
in $\oplus_{k=0}^{n-i-1}\mQ\langle 2k\rangle\langle -n+i+1\rangle$. 
That is $j=2n-2i-2-n+i+1=n-(i+1)$. On the other hand, the
formula~\cite[Theorem 3.11.4 (ii)]{BGS} gives 
$$[\tilde P(1^{i+1}0^{n-i})]= \sum_{k=0}^{n-i-1}[\tilde M(1^{i}0^{n-i-1-k}1
0^k\langle k\rangle].$$ Hence we finally get 
\begin{equation}
\label{EondomVerma}
{\bf E} \tilde M(1^i0^{n-i})\cong
\tilde P(1^i0^{n-i-1}1)  
\end{equation}
and  $\Psi([{\bf E} \tilde M(1^i0^{n-i})])=\sum_{k=0}^{n-i}\Psi([\tilde M(1^{i}0^{n-i-k}1
0^k\langle k\rangle])$. On the other hand we have to calculate $\blacktriangle(E)v_{\bf a}$,
where ${\bf a}=1^{i}0^{n-i}$. We get $\blacktriangle(E)v_{\bf a}=
\sum_{k=0}^{n-i}qv_{\bf a^k}$,  where ${\bf a^k}=1^i0^{n-i-1+k}10^k$. We get $\Psi([{\bf E}M])=E\Psi([M])$. The relation $\Psi([{\bf F}M])=F\Psi([M])$
follows from analogous calculations. The existence of the desired isomorphism
$\Phi$ follows. This proves part~\eqref{bases2.1}.\\

 Obviously, $\Phi'$ is an isomorphism of $\mC(q)$-modules. 
We first have to verify, that the functors ${\bf E}'$, ${\bf F}'$ and ${\bf
  K}$ satisfy the $U_q(\mathfrak{sl}_2)$ relations from
Definition~\ref{defquantumsl2}. Since we know that the functors ${\bf E}$,
${\bf F}$ and ${\bf K}$ satisfy these relations (see Theorem~\ref{funcrel}) it
is enough to verify the last equation. Since, however,  ${\bf
  F}_{i+1}'\bf{E}_i'={\bf F}_{i+1}\bf{E}_i$ and ${\bf E}_{i-1}'\bf{F}_i'={\bf
  E}_{i-1}\bf{F}_i$ this follows also directly from Theorem~\ref{funcrel}. It is left to show that $\Phi'$ is in fact a $U_q(\SL_2)$-morphism. We will
deduce this from part~\eqref{bases2.1}. Recall that if $\op{d}=\oplus_{i=0}^n \op{d}_{i}$ denotes the duality on
$\oplus_{i=0}^n\op{gmod-}{{}_i A}$ which fixes simple modules concentrated in
degree zero, then we put $\op{d}_i'=\langle 2i(n-i)\rangle\op{d}_{i}$
(see Section~\ref{psi}). Let for the moment be $\nabla=\op{d}\tilde\Delta({\bf a})$ for some
$\{0,1\}$-sequence containing exactly $i$ ones. Set $m_i=2i(n-i)$. We have the following equalities
\begin{eqnarray*}
  \Phi'\big({\bf E}'\nabla)&=& \Phi'(\langle -m_{i+1}\rangle{\bf E}_i\langle
  m_i\rangle\nabla\big)\\
&=&\Phi'\big(\langle -m_{i+1}\rangle{\bf E}_i\op{d}'_i\Delta({\bf a})\big)\\
&=&\Phi'\big(\langle -m_{i+1}\rangle\op{d}'_{i+1}{\bf E}_i\Delta({\bf a})\big)\\
&=&\Phi'\big(\op{d}_{i+1}{\bf E}_i\Delta({\bf a})\big).
\end{eqnarray*}
The first equation holds by definition of ${\bf E}'$, the second and the last
one by definition of $\op{d}'$. The remaining third equation is given by
Lemma~\ref{duality}. From the definitions of $\Phi$, $\Phi'$, and the duality
$\op{d}$ we get that
$\Phi'\big(\op{d}_{i+1}\rangle{\bf E}_i\Delta({\bf a})\big)$ is nothing else than
$\Phi({\bf E}_i\Delta({\bf a}))$, expressed in the standard basis, but the
involution $q\mapsto q^{-1}$ applies to all coefficients. From
part~\eqref{bases2.1} we know that this is the same as ${\bf E}\Phi(\Delta({\bf
  a}))$, expressed in the standard basis, but the
involution $q\mapsto q^{-1}$ applied to all coefficients. Since
$\Phi'(\nabla)=v^{\bf a}$, whereas $\Phi(\Delta({\bf a}))=v_{\bf a}$ we get
directly from the definition 
of the comultiplication $\triangle'$ (in comparison with $\triangle$) that 
$\Phi'({\bf E}'\nabla)=\Phi'(\op{d}_{i+1}{\bf E}_i\Delta({\bf
  a}))=E\Phi'(\nabla)$, where $E$ acts on the latter via the comultiplication
$\triangle'$. Analogous calculations show that $\Phi'({\bf
  F}'\nabla)=E\Phi'(\nabla)$. The equality $\Phi'({\bf
  K}\nabla)=K\Phi'(\nabla)$ is clear. The map $\Phi'$ is now in fact a
$U_q(\SL_2)$-morphism because the standard lifts of the dual Verma modules give rise to a $\mC(q)$-basis of
${\bf{G}}\big(\bigoplus_{i=0}^{n}{\op{gmod}-{}_i
  A}\big)$. Part~\eqref{bases2.1b} of the theorem follows.\\    

Part~\eqref{bases2.2} is well-known and follows for example from \cite[Theorem 3.11.4
(ii)]{BGS} and \cite[Remark 3.2 (2)]{SoKipp}.\\ 

To prove part~\eqref{bases2.3} we first verify that the form is bilinear. This
follows from the equalities
\begin{align*}
&  <1\otimes[M\langle k\rangle],
  1\otimes[N\langle l\rangle]>\pagebreak[0]\\
&=\sum_{i}\sum_{j}(-1)^j\DIM\HOM_\mathcal{D}\big(\op{d}(N\langle
  l\rangle), {\bf H} M\langle
  -i+k\rangle[j]\big)\; q^i\pagebreak[0]\\
&=\sum_{i}\sum_{j}(-1)^j\DIM\HOM_\mathcal{D}\big(\op{d}(N)\langle -l\rangle, {\bf H}M\langle  -i+k\rangle[j]\big)\; q^i
\end{align*}
\begin{align*}
&=\sum_{r}\sum_{j}(-1)^j\DIM\HOM_\mathcal{D}\big(\op{d}(N), {\bf H}M\langle
  -r\rangle[j]\big)\; q^{k+l}q^r\\
&=q^{k+l} <1\otimes[M],1\otimes[N]>.
\end{align*}  
On the other hand we also have
\begin{eqnarray*}
&&  <1\otimes [\tilde M(y\cdot\omega_i)],
  1\otimes [\op{d}\tilde M(x\cdot\omega_i)]>\\
&=&\sum_{i}\sum_{j}(-1)^j\DIM\HOM_\mathcal{D}
\big(\tilde M(x\cdot\omega_i), {\bf
  H} \tilde M(y\cdot\omega_i)\langle i\rangle[j]\big)\; q^i\\
&=&\sum_{i}\sum_{j}(-1)^j\DIM\HOM_\mathcal{D}\big(\tilde
  M(x\cdot\omega_i),\op{d}\tilde M(w_0y\cdot\omega_i)\langle i\rangle[j]\big)\; q^i\\
&=&\DIM\HOM_\mathcal{D}(\tilde
  M(x\cdot\omega_i),\op{d}\tilde M(w_0y\cdot\omega_i))=\delta_{x,w_oy}.
\end{eqnarray*}
From \eqref{form} the statement~\eqref{bases2.3} of the theorem follows.\\
 
Let us prove part~\eqref{bases2.4}.  By definition the standard basis of $\cM$ is mapped to
the standard lifts of Verma modules, and they are mapped to the standard basis
in $V_1^{\otimes n}$. By Proposition~\ref{SchurWeyl}~\eqref{twist2}, the twisted standard
basis is mapped to the standard lifts of the dual Verma modules, and they are
(by definition) mapped to the dual standard basis in $(V_1^{\otimes n})'$.  
From part~\eqref{bases2.2} and Proposition~\ref{SchurWeyl} we get that
  the twisted positive basis of $\cM$ corresponds to the standard lifts of
  tilting modules. The formula \cite[Theorem 2.6]{FKK} together with
  \cite[Proposition 3.4]{SoKipp} explicitly show that the twisted positive
  basis corresponds to the canonical basis. Finally we have
\begin{eqnarray*}
&& <[\tilde T(y\cdot\la)],
  [\op{d}\tilde L(x\cdot\la)]\rangle\\
&=&\sum_{i}\sum_{j}(-1)^j\DIM\HOM_\mathcal{D}\big(\tilde L(x\cdot\la), {\bf
  H} \tilde T(y\cdot\la)\langle i\rangle[j]\big)\; q^i\\
&=&\DIM\HOM_\mathcal{D}\big(\tilde L(x\cdot\la),\op{d}\tilde
  I(w_0y\cdot\la)\big)\; q^i\\
&=&\delta_{x,w_oy}.
\end{eqnarray*}
From formula~\eqref{form2} it follows that the standard lifts of the simple
modules correspond to the dual canonical basis. Finally, it is known that the negative
self-dual basis corresponds to the dual canonical basis (see \eg \cite[Theorem
2.5']{FKK}.\\ 

It is left to prove part~\eqref{bases2.5} of the theorem. Of course, $\gamma$
defines an isomorphism of $\mC(q)$-modules. That tilting modules correspond to
the positive self-dual bases as stated follows directly from
part~\eqref{bases2.2} of the theorem together with the
isomorphisms~\eqref{twist1} and~\eqref{twist2}. That simple modules correspond
to the negative self-dual basis elements follows directly from the second half
of part~\eqref{bases2.4}.
So we are done, the theorem follows.    
\end{proof}

\pagebreak
We finish with two additional remarks:
\begin{prop}
  \begin{enumerate}[(a)]
  \item \label{a} There are isomorphisms of functors 
    \begin{equation*}
{\bf H}_{w_0}\,{\bf E}\cong{\bf E}\,{\bf H}_{w_0},\quad 
{\bf H}_{w_0}\,{\bf F}\cong{\bf F}\,{\bf H}_{w_0},\quad {\text and}\quad  
{\bf H}_{w_0}\,{\bf K}\cong{\bf K}\,{\bf H}_{w_0}.
\end{equation*}
  \item \label{b} The ``categorical'' bilinear form $< ,>$ from
  Theorem~\ref{bases2}\eqref{bases2.3} satisfies 
\begin{equation*}
<\blacktriangle(x)(v_i\otimes
v_j), v^k\otimes v^l>=<v_i\otimes
v_j, \blacktriangle'(\omega(x)) (v^k\otimes v^l)>,  
\end{equation*}
  \end{enumerate}
\end{prop}
\begin{proof}
  The first statement is clear if we forget the grading. To prove the first
  isomorphism it is therefore enough to show
${\bf H}_{w_0}^2{\bf E}\cong{\bf
  E}{\bf H}_{w_0}^2$, even ${\bf H}_{w_0}^2{\bf E}_i\tilde M({\bf a})\cong{\bf
  E}_i{\bf H}_{w_0}^2\tilde M({\bf a})$, where ${\bf a}=(1,1,\ldots,1,
  0,0,\ldots 0)$ with exactly $i$ ones. From the formula~\eqref{EondomVerma}
  we know that ${\bf E}_i\tilde M({\bf a})\cong\tilde P({\bf b})$ for some
  $\{0,1\}$-sequence ${\bf b}$. Hence 
  \begin{align*}
    {\bf H}_{w_0}^2{\bf E}_i\tilde M({\bf a})\quad\cong\quad{\bf H}_{w_0}^2\tilde
    P({\bf b})\quad\cong\quad\tilde I({\bf b})\langle -2l(w_0^i)\rangle
  \end{align*}
by Proposition~\ref{SchurWeyl}. On the other hand 
\begin{alignat*}{10}
  {\bf
  E}_i{\bf H}_{w_0}^2\tilde M({\bf a})&\quad&\cong&\quad&{\bf
  E}_i(\op{d}\tilde M({\bf a}))\langle -2l(w_0^i)\rangle
&\quad&\cong&\quad&(\op{d}{\bf E}_i\tilde M({\bf a}))\langle -2l(w_0^i)\rangle\\
&\quad&\cong&\quad&(\op{d}P({\bf b}))\langle -2l(w_0^i)\rangle
&\quad&\cong&\quad&\tilde I({\bf b})\langle -2l(w_0^i)\rangle
\end{alignat*}
by Lemma~\ref{duality} and, again, Proposition~\ref{SchurWeyl}. The first
isomorphism of the proposition follows. The second is proved analogously, and the third is
clear. 

To prove statement~\eqref{b} set $m_i=i(n-i)$ and write just $\hom(N,M)$, short for
$\sum_{i}\sum_{j}(-1)^j\DIM\HOM^j_\mathcal{D}(N,M)q^i$, and calculate 

  \begin{alignat*}{6}
&  <{\bf E}_iM, N>\\
=&\hom(\op{d}N, {\bf H}{\bf E}_i M)\\
=&\hom\big(\op{d}N, {\bf H}_{w_0}{\bf E}_i M\langle l(w_0^{i+1})\rangle\big)&\text{ (by
definition of {\bf H})}\\
=&\hom\big(\op{d}N, {\bf E}_i{\bf H}_{w_0}M\langle l(w_0^{i+1})\rangle\big)&\text{ (by
  part~\eqref{a})}\\
=&\hom\big(\op{d}N, {\bf E}_i{\bf K}_i\langle 1\rangle {\bf H}_{w_0}{\bf
  K}^{-1}\langle -1\rangle M\langle l(w_0^{i+1})\rangle\big)\\
=&\hom\big( {\bf F}_{i+1}(\op{d}_{i+1}N), {\bf H}_{w_0}{\bf
  K}^{-1}\langle -1\rangle M\langle l(w_0^{i+1})\rangle\big)&\text{(Prop.~\ref{all-adjoint})}\\
=&\hom\big( {\bf F}_{i+1}\op{d}'_{i+1}\langle m_{i+1}\rangle N, {\bf H}_{w_0}{\bf
  K}^{-1}\langle -1\rangle M\langle l(w_0^{i+1})\rangle\big)&\text{(see
  Section~\ref{psi})}\displaybreak[0]\\
=&\hom\big( \op{d}'_{i}{\bf F}_{i+1}\langle m_{i+1}\rangle N, {\bf H}_{w_0}{\bf
  K}^{-1}\langle -1\rangle M\langle l(w_0^{i+1})\rangle\big)&\text{(Lemma~\ref{duality})}\displaybreak[0]\\
=&\hom\big(\langle m_i\rangle\op{d}_{i}{\bf F}_{i+1}\langle m_{i+1}\rangle N, {\bf H}_{w_0}{\bf
  K}^{-1}\langle -1\rangle M\langle l(w_0^{i+1})\rangle\big)&\text{(see Section~\ref{psi})}\displaybreak[0]\\
=&\hom\big(\op{d}_{i}\langle-m_i\rangle{\bf F}_{i+1}\langle m_{i+1}\rangle N, {\bf H}_{w_0}{\bf
  K}^{-1}\langle -1\rangle M\langle l(w_0^{i+1})\rangle\big)&\text{(by
  definition of $\op{d}$)}\displaybreak[0]\\
=&\hom\big(\op{d}_{i}{\bf F}_{i+1}'N, {\bf H}_{w_0}{\bf
  K}^{-1}\langle -1\rangle M\langle l(w_0^{i+1})\rangle\big)&\text{(by
  definition of ${\bf F}'$)}
\end{alignat*}
We claim that the latter is isomorphic to 
\begin{eqnarray*}
&&\hom\big( \op{d}{\bf F}_{i+1}'(N), {\bf H}_{w_0}M\langle
 (n-2i-1+l(w_0^{i+1})\rangle\big)\\
 &=&\hom\big( \op{d}{\bf F}_{i+1}'(N), {\bf H}_{w_0}M\big)\\
 &=&<M, {\bf F}_{i+1}' N>.
\end{eqnarray*}

To verify the claim we used the following: Note that
$l(w_0^i)=\frac{1}{2}(i(i-1)+(n-i)(n-i-1))$. Hence 
\begin{eqnarray*}
 &&l(w_0^{i+1})-l(w_0^i)\\
&=&\frac{1}{2}\big((i+1)i+(n-i-1)(n-i-2)-(i(i-1)+(n-i)(n-i-1))\big)\\
&=&\frac{1}{2}(2i-2n+2i+2)=2i-n+1. 
\end{eqnarray*}
Together with the definition of ${\bf H}$, the formula $(*)$ follows. Similarly
we get $<{\bf F}_iM, N>= <M,{\bf E}_{i-1}'N>$. The formula $<{\bf K}_iM, N>=
<M,{\bf K}_{i}N>$ is obvious from the bilinearity of the form. The proposition
follows. 
\end{proof}

\begin{remark}
\label{finalremark}
{\rm
  Theorem~\ref{bases2} generalises to arbitrary tensor products $V_{\bf d}$
  such that the standard basis corresponds to standard modules, the dual
  standard basis to dual standard modules, the canonical basis to tilting
  modules and the dual canonical bases to simple modules. Since the proofs
  involve deeply the properly stratified structure of the category of
  Harish-Chandra bimodules, the arguments will appear in another paper, where we
  moreover show that Proposition~\ref{JonesWenzl} is also true in the graded setup.}
\end{remark}

\subsection{Categorification dictionary} 
\small
\begin{tabular}{|c|c|}
\hline
\hspace{0.1in} & \hspace{0.1in} \\
{\bf Quantum} $\mf{sl}_2$ {\bf and its representations}  & 
{\bf Functors and categories} \\    
\hline
\hspace{0.1in}&  category $\Q$-Vect \\
ring $\Z[q,q^{-1}]$ & of graded vector spaces \\ 
\hline 
\hspace{0.1in}&  \hspace{0.1in} \\
multiplication by $q$ & grading shift up by $1$\\   
\hline   
\hspace{0.1in} &  \hspace{0.1in} \\
representation $V_n$ & category $\cC\cong\oplus_{i=0}^n C^i\gMOD$,\\
& where $C^i$ is the cohomology ring \\
&of a Grassmannian \\     
\hline 
\hspace{0.1in} &  \hspace{0.1in} \\
weight spaces of $V_n$ & the summands of $\cC\cong\oplus_{i=0}^n C^i\gMOD$,\\
\hline 
\hspace{0.1in} &  \hspace{0.1in} \\
semilinear form $<,>: \quad 
V_n\times V_n \to \Q(q)$ & bifunctor 
$\HOM_{\mc{C}}(\ast, \ast)$ \\ 
\hline 
\hspace{0.1in} & \\
canonical basis of $V_n$ & indecomposable projective modules in \\
&$\oplus_{i=0}^n C^i\gMOD$\\ 
\hline 
\hspace{0.1in} &  \hspace{0.1in} \\
dual canonical basis of $V_n$  &  simple modules $\oplus_{i=0}^n C^i\gMOD$\\
\hline 
\hspace{0.1in} &  \hspace{0.1in} \\
representation $V_1^{\otimes n}$ & graded version of
$\oplus_{i=0}^n\,{}_{i;n}\cO$,\\
&certain blocks of the category $\cO(\mathfrak{gl}_n)$.\\ 
\hline 
\hspace{0.1in} &  \hspace{0.1in} \\
weights spaces of $V_1^{\otimes n}$ & blocks in the graded version of
$\oplus_{i=0}^n\,{}_{i;n}\cO$,\\
\hline 
\hspace{0.1in}&  \hspace{0.1in} \\
standard basis of $V_1^{\otimes n}$ & standard (=Verma) modules\\ 
\hline 
\hspace{0.1in}&  \hspace{0.1in} \\
dual standard basis of $V_1^{\otimes n}$ & dual standard (=dual Verma) modules\\ 
\hline 
\hspace{0.1in}&  \hspace{0.1in} \\
canonical basis of $V_1^{\otimes n}$ & indecomposable tilting  modules\\ 
\hline 
\hspace{0.1in} &  \hspace{0.1in} \\
dual canonical basis of $V_1^{\otimes n}$  &  simple modules\\
\hline 
\hspace{0.1in} &  \hspace{0.1in} \\
representation $V_{\bf d}$ & graded version of
$\oplus_{i=0}^n\,{}_\mu\cH_{\omega_i}(\GL_n)$,\\ 
&certain blocks of the category\\ 
&of Harish-Chandra bimodules $\cH(\mathfrak{gl}_n)$.\\
\hline 
\hspace{0.1in} &  \hspace{0.1in} \\
weight spaces of $V_{\bf d}$ & blocks in the graded version of
$\oplus_{i=0}^n\,{}_\mu\cH_{\omega_i}(\GL_n)$,\\ 
\hline
\hspace{0.1in} &  \hspace{0.1in} \\
anti-automorphism $\tau$ & taking the right adjoint functor \\
\hline 
\hspace{0.1in} &  \hspace{0.1in} \\
involutions $\psi_n$, $\psi$  &  duality functor\\
\hline 
\hspace{0.1in} &  \hspace{0.1in} \\
Cartan involution $\sigma$ and $\sigma_n$  & 
the equivalences $\hat{\sigma}$
and $\hat{\sigma}_n$ \\
&arising from the equivalences of categories\\  
&$\cA_i^{\mu}\cong\cA_{n-i-1}^{\mu}$.\\
\hline 
\end{tabular}
\normalsize

\section{Geometric categorification}
\label{geometry}
Our treatment of the graded version of category $\cO$ (based on \cite{BGS})
uses substantially the cohomology rings of partial flag varieties. A more
systematic study of the structure of these rings naturally leads to an
alternative categorification which we call ``geometric''. We first consider
the geometric categorification of simple finite dimensional
$U_q(\mathfrak{sl}_2)$-modules using the cohomology rings of partial flag
varieties and relations between them. Then we proceed to the geometric
categorification of general tensor products using certain algebras of
functions which generalise the cohomology rings and point towards the Borel-Moore
homology of generalised Steinberg varieties. We conclude this section
formulating open problems related to the geometric categorification.   
 
\subsection{From algebraic to geometric categorification}
The categorification of simple $U_q(\SL_2)$-modules we propose gives rise to a
categorification of simple $\cU(\SL_2)$-modules by forgetting the grading. The
categorification for simple  $\cU(\SL_2)$-modules obtained in this way is
exactly the one appearing in \cite{CRsl2}. It doesn't seem to be obvious
from the approach of \cite{CRsl2} why exactly this categorification plays an
important role. As a motivation for choosing this categorification we first show that it naturally
emerges from our Theorem~\ref{cat} as follows: Theorem~\ref{cat} provides an isomorphism
$V_n \cong{\bf G}(\bigoplus_{i=0}^n\cA_i^{-\rho})$, since $W_{-\rho}=W$. Each
of the categories $\cA_i^{-\rho}$ contains (up to isomorphism and shift in the
grading) one single indecomposable
projective object (\cite[6.26]{Ja2}), hence also (up to isomorphism and
grading shift) one simple object. On the other hand the categories
$C^i\gMOD$ also have (up to isomorphism and shift in the
grading) one single simple object $S_i$. In particular, the Grothendieck
groups of the two categories coincide. However, we have the following stronger result:

\begin{prop}
\label{HCC}
  There is an equivalence of categories 
  \begin{eqnarray*}
    F:\bigoplus_{i=0}^n\cA_i^{-\rho}\rightarrow\bigoplus_{i=0}^n C^i\gMOD,
  \end{eqnarray*}
which intertwines the functors ${\bf E}_i$, ${\bf F}_i$, and ${\bf K}_i$ with
the functors 
 \begin{eqnarray*}
  E_i&=&\RES_{i+1}^{i,i+1} C^{i,i+1}\otimes_{C^i}\langle -n+i+1\rangle,\\
  F_i&=&\RES_{i-1}^{i,i-1} C^{i,i-1}\otimes_{C^i}\langle-i+1\rangle,\\
  K_i&=&\langle 2i-n\rangle.
\end{eqnarray*}
\end{prop}
\begin{proof}
  Let $P\in\cA_i^{-\rho}$ be the (up to isomorphism) unique  
  indecomposable projective module such that its head is concentrated in
  degree zero. Then we have an isomorphism of graded algebras 
  $\END_{\cA_i^{-\rho}}(P)\cong\END_{\cA_i}(P(w_0\cdot\omega_i))\cong
  C^i$. (The first isomorphism follows just from the definition of
  $\cA_i^{-\rho}$, the second is \cite[Endomorphismensatz]{Sperv} together
  with the definition of the grading on $\cA_i$.) Note that any module in
  $\cA_i^{-\rho}$ is a quotient of some $P\langle k\rangle$, $k\in\mZ$. In
  other words, $P$ is a $\mZ$-generator of $\cA_i^{-\rho}$ in the sense of
  \cite[E.3.]{AJS}. Hence (see \eg \cite[Proposition E.4.]{AJS}), the functor
  \begin{eqnarray*}
    \bigoplus_{k\in\mZ}\HOM_{\cA_i^{-\rho}}(P\langle
    k\rangle,\bullet):\quad\cA_i^{-\rho}\longrightarrow
    \op{gmod}-\END_{\cA_i^{-\rho}}(P)\cong C^i-\op{gmod}
  \end{eqnarray*}
defines an equivalence. Of course, we could fix any $l_i\in\mZ$ and still
have an equivalence of categories
\begin{eqnarray*}
    \bigoplus_{k\in\mZ}\HOM_{\cA_i^{-\rho}}(P\langle
    k+l_i\rangle,\bullet):\quad\cA_i^{-\rho}\longrightarrow
    \op{gmod}-\END_{\cA_i^{-\rho}}(P\langle l_i\rangle)\cong C^i-\op{gmod}.
  \end{eqnarray*}
Up to a grading shift, these equivalences obviously intertwine the functors in
question. (This follows directly from
the definitions of the functors and general arguments, see \eg
\cite[2.2]{Bass}.) It follows easily from the definitions of the graded lifts
that $l_i=-i(n-i)$ gives the required equivalences.   
\end{proof}
From Proposition~\ref{HCC} it follows in particular, that ${\bf
  G}(\bigoplus_{i=0}^n C^i\gMOD)$ becomes a $U_q(\SL_2)$-module via the
  functors $E=\oplus_{i=0}^n E_i$, $F=\oplus_{i=0}^n F_i$ and $K=\oplus_{i=0}^n K_i$. From
  Theorem~\ref{cat} we know that the resulting module is isomorphic to
  $V_n$. One could, of course, ask the question which $C^i$-modules correspond
  to the (dual) canonical basis elements. This will be answered in the next
  section. 

\subsection{Categorification of simple $U_q(\SL_2)$ modules via modules over
  the cohomology rings of Grassmannians}
\label{sectionVn}
We now want to describe the geometric categorification, motivated by
Proposition~\ref{HCC} separately, not as a consequence of the main
categorification theorem. The reason for this is that we
propose a generalisation of this construction to a ``geometric''
categorification of $V_{\bf d}$ using certain algebras of functions
(Section~\ref{between}). Later on we will discuss the relation of these
function algebras with the Borel-Moore
homology rings of Steinberg varieties and to the algebraic categorification
from Section~\ref{mainresult}.

Let $n\in\mZ_{>0}$ and let $I:=\{i_1, i_2,\ldots i_{k}\}\subseteq\{1,2,\ldots
n\}$. We consider the corresponding partial flag variety $\cG_I$ given by all flags
$\{0\}\subset F_1\subset F_2\subset\cdots\subset F_k\subset\mC^n$ where
$\{\DIM_\mC F_{j}\}=I$. Let $C^I=H^\bullet(\cG_I,\mC)$ denote its
cohomology. We will be only interested in the cases $k=1$, $k=2$ and $I={i}$,
and $I=\{i_1,i_2\}$.  
Let denote
$C^I\MOD$ the category of finitely generated $C^I$-modules. Each of
these categories has exactly one simple object, the trivial module
$\mQ$. If $I=\{i\}$ then the corresponding variety is a Grassmannian and we
denotes its cohomology ring by $C^i$. The Grothendieck group of $\oplus_{i=0}^n C^i\MOD$ is free of rank $n+1$.
The rings $C^I$ have a natural (positive, even) $\mZ$-grading and we may consider the
categories $C^I\gMOD$ of finitely generated {\it graded} $C^I$-modules and
the category $\oplus_{i=0}^n C^i\gMOD$. Its Grothendieck group is then a free
$\mZ[q,q^{-1}]$-module of rank $n+1$, where $q^i$ acts by shifting the grading
degree by $i$. Hence we have a candidate for a categorification of the simple $U_q(\SL_2)$-module $V_n$. 
We need also functors giving rise to the $U_q(\SL_2)$-action. If $J\subseteq I$ then there is an obvious surjection $\cG_I\rightarrow \cG_J$ inducing 
an inclusion $C^J\rightarrow C^I$ of rings. Let $\RES^I_J:C^I\gMOD\rightarrow
C^J\gMOD$ denote the restriction functor.    
For $0\leq i\leq n$ set  
\begin{eqnarray*}
  E_i&=&\RES_{i+1}^{i,i+1} C^{i,i+1}\otimes_{C^i}\langle -n+i+1\rangle,\\
  F_i&=&\RES_{i-1}^{i,i-1} C^{i,i-1}\otimes_{C^i}\langle-i+1\rangle,\\
  K_i&=&\langle 2i-n\rangle.
\end{eqnarray*}
Let $S_i\in C^i\gMOD$ be the
simple (trivial) module concentrated in degree $0$. Denote by $\mV_i T_i\in
C^i\gMOD$ the projective cover of $S_i\langle
-i(n-i)\rangle$. (The notation coincides with the one from Section~\ref{graded
  O}.)
We consider the following functors $E:=\oplus_{i=0}^n E_i$, $F:=\oplus_{i=0}^n
F_i$,  $K:=\oplus_{i=0}^n K_i$ from $\oplus_{i=0}^n C^i\gMOD$ to
itself. Since the functors are exact, and commute with shifts in the grading, they
induce $\mZ[q,q^{-1}]$-morphisms on $[\oplus_{i=0}^n
C^i\MOD]$. 
Note that $\frac{q^j-q^{-j}}{q-q^{-1}}=\sum_{k=0}^{j-1} q^{j-1-2k}$ for
$j\in\mZ_{>0}$. Therefore, if $2i-n>0$, we use the notation 
\begin{eqnarray}
  \label{eq:quotoffunc}
  \frac{K_i-K_i^{-1}}{q-q^{-1}}&=&-\frac{K_{n-i}-{K_{n-i}}^{-1}}{q-q^{-1}}\nonumber\\&=&
  \begin{cases}
    \displaystyle\bigoplus_{k=0}^{2i-n-1}\ID\langle 2i-n-1-2k\rangle&\text{ if $2i-n>0$},\\
    \hspace{1cm}0&\text{ if $2i-n=0$},
  \end{cases}  
\end{eqnarray}
(as endofunctors of $\oplus_{i=0}^n C^i\gMOD$). 

The following result categorifies simple $U_q(\mathfrak{sl}_2)$-modules: 

\begin{theorem}[The categorification of $V_n$]
\label{Vn}
Fix $n\in\mZ_{>0}$. 
  \begin{enumerate}[(a)]
\item \label{Vn1} The functors $E$, $F$ and $K$, $K^{-1}$ satisfy the
  relations   
\begin{eqnarray} 
  KE=\langle 2\rangle EK,\quad\quad KF=\langle -2\rangle FK,&&KK^{-1}=\ID=K^{-1}K,\label{easy}\\
  E_{i-1}F_i\cong F_{i+1}E_i\oplus\frac{K_i-K_i^{-1}}{q-q^{-1}},&&\text{   if
  $2i-n\geq 0$},\label{difficult1}\\
 F_{i+1}E_i\cong E_{i-1}F_i\oplus-\frac{K_i-K_i^{-1}}{q-q^{-1}},&&\text{   if
  $n-2i>0$}.
\label{difficult2}
  \end{eqnarray} 
In the Grothendieck group we have the equality 
\begin{eqnarray*}
(q-q^{-1})\big({E}_{i-1}^{\bf G}{F}_i^{\bf G}-{F}_{i+1}^{\bf G}{E}_i^{\bf G}\big)
&=&{K}_i^{\bf G}-({K}_i^{-1})^{\bf G}.
\end{eqnarray*}
Hence they induce a $U_q(\SL_2)$-structure on ${\bf G}(\oplus_{i=0}^n
C^i\MOD)$. 
\item  \label{Vn2} With respect to this structure, there is an isomorphism of $U_q(\SL_2)$-modules
  \begin{eqnarray*}
  \begin{array}[tb]{ccccc}
 V_n &\cong&{\bf
  G}(\bigoplus_{i=0}^n C^i\gMOD)\\
\\
   v^i&\mapsto&1\otimes [S_i]\\
   v_i&\mapsto&1\otimes [\mV_i T_i].
  \end{array}
\end{eqnarray*}
\item \label{Vn3} The involution $\sigma_n$ from \eqref{omegan} can be categorified in the
  following sense: There is an equivalence of categories
  \begin{eqnarray*}
    \hat{\sigma}_n:\quad \bigoplus_{i=0}^n C^i\gMOD&\cong&\bigoplus_{i=0}^n
    C^i\gMOD\\
S_i&\mapsto&S_{n-i}\\
\mV_i T_i&\mapsto&\mV_{n-i} T_{n-i}
  \end{eqnarray*}
such that 
\begin{eqnarray*}
  \begin{array}[th]{lcl}
  \hat{\sigma}_n\,\hat{\sigma}_n&\cong&\ID,\\
  \hat{\sigma}_n\,E\,\hat{\sigma}_n&\cong&F,\\
  \hat{\sigma}_n\,F\,\hat{\sigma}_n&\cong&E,\\
  \hat{\sigma}_n\,K\,\hat{\sigma}_n&\cong&K^{-1}\\
  \hat{\sigma}_n\langle k\rangle\,\hat{\sigma}_n&\cong&\langle k\rangle
  \end{array}
\end{eqnarray*}
for any $k\in\mZ$.
\item \label{Vn4} The involution $\psi_n$ from \eqref{psin} can be categorified in the following
  sense: The duality $M\mapsto \op{d}M:=\HOM_\mQ(M,\mQ)$ defines an
  involution  
 \begin{eqnarray*}
    \hat\psi_n:\quad \bigoplus_{i=0}^n C^i\gMOD&\cong&\bigoplus_{i=0}^n
    C^i\gMOD\\
    \mV_iP_i&\mapsto& \mV_iP_i,  
\end{eqnarray*}
and 
\begin{eqnarray*}
  \begin{array}[th]{lcl}
  \hat\psi_n\,E\,\hat\psi_n &\cong& F\\
  \hat\psi_n\,F\,\hat\psi_n&\cong& E\\
  \hat\psi_n\,K\,\hat\psi_n&\cong& K^{-1}\\
  \hat\psi_n\,\langle i\rangle\,\hat\psi_n&\cong&\langle -i\rangle.
  \end{array}
\end{eqnarray*}
for any $i\in\mZ$.   
\item \label{Vn5} The semilinear form becomes the form
  \begin{eqnarray*}
    <[M],[N]>=\sum_{i}\sum_{j}(-1)^j\DIM(\op{Ext}^j_{\cC}(N,M\langle -i\rangle)) q^i
  \end{eqnarray*}
\item \label{Vn6} The antilinear anti-automorphism $\tau$ can be viewed as the operation of
  taking right adjoint functors, namely: There are pairs of adjoint functors
\begin{eqnarray*}
({ E}, { F K}^{-1}\langle 1\rangle),& ({F}, {EK}\langle 1\rangle), \quad
({K},{K}^{-1}),\quad ({K}^{-1}
,K)
\end{eqnarray*}
and $(\langle -k\rangle, \langle k\rangle)$ for $k\in\mZ$.
 \end{enumerate}
\end{theorem}

We would like to stress again that essentially the same categorification, in the ungraded 
case (without $q$), was previously constructed by Chuang and Rouquier
\cite{CRsl2}. Before we prove the proposition we state the following 

\begin{lemma}
\label{dual}
  Consider the graded duality $\op{d}:
  \oplus_{i=0}^n C^i\gMOD\rightarrow\oplus_{i=0}^n C^i\gMOD$, $M\mapsto
  \op{d}(M)$, where $(\op{d}M)_j=\HOM_\mQ(M_{-j},\mQ)$. Then
  there are isomorphisms of functors $\op{d} E\cong E\op{d}$ and $\op{d} F\cong
  F\op{d}$.  
\end{lemma}

\begin{proof}
  To establish an isomorphism $\op{d} E\cong E\op{d}$ it is enough to show that
  $\op{d} E_i\mV_i\cong E_i\op{d}\mV_i$. Let us for a moment forget the
  grading. From Proposition~\ref{EFandV} we get $\op{d} E_i\mV_i\cong
  \op{d}\mV_i\cE_i$. On the other hand, $\op{d}\mV_i\cE_i\cong
  \mV_i\op{d}\cE_i\cong\mV_i\cE_i\op{d}$ by \cite[Lemma 8]{Sperv} and \cite[4.12
  (9)]{Ja2}. Finally we have
  $\mV_i\cE_i\op{d}\cong E_i\mV_i\op{d}\cong E_i\op{d}\mV_i$ again by
  Proposition~\ref{EFandV} and \cite[Lemma 8]{Sperv}. We get an isomorphism of
  functors $\op{d} E_i\cong E_i\op{d}:\oplus_{i=0}^n C^i\MOD\rightarrow
  \oplus_{i=0}^n C^i\MOD$. Since the functors $\cE_i$ are indecomposable (see
  Proposition~\ref{EFandV}), we can find (see \eg \cite[Lemma 2.5.3]{BGS}) some $k_i\in\mZ$ such that $\op{d} E_i\cong
  E_i\op{d}\langle k_i\rangle:\oplus_{i=0}^n C^i\gMOD\rightarrow
  \oplus_{i=0}^n C^i\gMOD$. On the other hand, using
  Lemma~\ref{free} we get
  isomorphisms of graded vector spaces 
  \begin{eqnarray*}
    E_i(S_i)&\cong& C^{i,i+1}\otimes_{C^i}\mQ\langle -n+i+1\rangle\\
&\cong&\oplus_{r=0}^{n-i+1}\mQ\langle 2r-n+i-1\rangle\\
&\cong&
\mQ\langle
n-i-1\rangle\oplus\mQ\langle n-i-3\rangle\oplus\cdots\oplus\mQ\langle
-n+i+3\rangle\oplus\mQ\langle-n+i+1\rangle\\
&\cong&\op{d} E_i(S_i). 
  \end{eqnarray*} 
Hence $k_i=0$ and so $\op{d} E\cong E\op{d}$. The arguments establishing an isomorphism  $\op{d} F\cong F\op{d}$ are analogous. 
\end{proof}
  
\begin{proof}[Proof of Proposition~\ref{Vn}]
The relations~\eqref{easy} follow directly from the definitions. The
verification of the relations~\eqref{difficult1} and \eqref{difficult2} is
much more involved. We prove it here only on the level of the
Grothendieck group. For the full statement we refer to
Theorem~\ref{funcrel} and Proposition~\ref{HCC}. 
By Lemma~\ref{free} we get
isomorphisms of graded $C^i$-modules
\begin{eqnarray*}
  F_{i+1}E_i(\mQ)&\cong&
  C^{i,i+1}\otimes_{C^{i+1}}C^{i,i+1}\otimes_{C^i}\mQ\langle -n+1\rangle\\
&\cong&
  \oplus_{l=0}^{i}\oplus_{k=0}^{n-i-1}\mQ\langle2k+2l\rangle\langle -n+1\rangle,\nonumber\\
E_{i-1}F_i(\mQ)&\cong&C^{i,i-1}\otimes_{C^{i-1}}C^{i,i-1}\otimes_{C^i}\mQ\langle
  -n+1\rangle\\
&\cong&
  \oplus_{l=0}^{n-i}\oplus_{k=0}^{i-1}\mQ\langle2l+2k\rangle\langle -n+1\rangle.
\end{eqnarray*}
If $n-2i\geq 0$, then  
\begin{eqnarray*}
{F}_{i+1}^{\bf G}{E}_i^{\bf G}-{E}_{i-1}^{\bf G}{F}_i^{\bf G}
&=&
\big(\oplus_{k=0}^{n-i-1}\langle(2k+2i)\rangle\langle(1-n)\rangle\big)^{\bf G}\\
&&-\;(\oplus_{k=0}^{i-1}\langle(2(n-i)+2k)\rangle\langle 1-n\rangle)^{\bf G}\\
&=&
\begin{cases}
0&\text{ if $n-2i=0$,}\\
\big(\oplus_{k=0}^{n-2i-1}\langle (2i+2k)\rangle\langle 1-n\rangle\big)^{\bf G}&\text{ if $n-2i>0$,}
\end{cases}\\
&=& [n-2i]\op{id}.
\end{eqnarray*}
If $n-2i<0$, then 
\begin{eqnarray*}
{E}_{i-1}^{\bf G}{F}_i^{\bf G}-{F}_{i+1}^{\bf G}{E}_i^{\bf G}
&\cong&
\big(\oplus_{k=0}^{i-1}\langle (2(n-i)+2k)\rangle\langle(-n+1)\rangle\big)^{\bf G}\\
&-&\big(\oplus_{k=0}^{n-i-1}\langle (2k+2i)\rangle\langle 1-n\rangle\big)^{\bf G}\\
&=&\big(\oplus_{k=0}^{2i-n-1}\langle (n-2i+1+2k)\rangle\big)^{\bf G}\\
&=& [2i-n]\op{id}.
\end{eqnarray*}
In particular, 
\begin{eqnarray*}
(q-q^{-1})\big({E}^{\bf G}{F}^{\bf G}-F^{\bf G}E^{\bf G}\big)
&=&K^{\bf G}-(K^{-1})^{\bf G}.
\end{eqnarray*}
This proves the part~\eqref{Vn1} of the proposition on the level of the
Grothendieck group.

Obviously, the map  $v^i\mapsto[S_i]$ from part~\eqref{Vn2} defines an
  isomorphism of vector spaces. We first have to verify that this is in fact a
  morphisms of $U_q(\SL_2)$-modules. 
Using
  Lemma~\ref{free} we get
  isomorphisms of $\mZ[q,q^{-1}]$-modules
  \begin{eqnarray*}
    &&[E_i(S_i)]=[C^{i,i+1}\otimes_{C^i}\mQ\langle -n+i+1\rangle]\\
    &=&[\oplus_{r=0}^{n-i-1}\mQ\langle 2r-n+i-1\rangle]\\
&=&
[\mQ\langle
n-i-1\rangle\oplus\mQ\langle n-i-3\rangle\oplus\cdots\oplus\mQ\langle
    -n+i+3\rangle\oplus\mQ\langle-n+i+1\rangle]\\
&=&
[S_{i+1}\langle
n-i-1\rangle\oplus S_{i+1}\langle n-i-3\rangle\oplus\cdots\oplus S_{i+1}\langle
    -n+i+3\rangle\oplus\mQ\langle-n+i+1\rangle]\\
&=&
\sum_{k=0}^{n-i-1}[S_{i+1}\langle-n+i+1+2k\rangle].
  \end{eqnarray*}
Similarly we have
\begin{eqnarray*}
   [F_i(S_i)]&=&[\mQ\langle
-i+1\rangle\oplus\mQ\langle -i+3\rangle\oplus\cdots\oplus\mQ\langle
i-3\rangle\oplus\mQ\langle i-1\rangle]\\
&=&[S_{i-1}\langle
-i+1\rangle\oplus S_{i-1}\langle -i+3\rangle\oplus\cdots\oplus S_{i-1}\langle
i-3\rangle\oplus S_{i-1}\langle i-1\rangle]\\
&=&
\sum_{k=0}^{i-1}[S_{i-1}\langle 1-i+2k\rangle]. 
\end{eqnarray*}
This fits with the formula~\eqref{dualcan}. Hence, the assignment $v^i\mapsto
[S_i]$ defines a $U_q(\SL_2)$-morphism, where the simple modules
concentrated in degree zero correspond
to the dual canonical basis elements. It is left to show that $v_i$ is mapped
to $[\mV_i T_i]$.  
Recall that the algebra $C^i$ has a basis ${b_y}$ naturally indexed by
elements $y\in W^i$. With our convention on the grading $b_y$ is homogeneous
of degree $2l(y)$. Hence, from the definition of $\mV_i T_i$ we get 
\begin{eqnarray}
\label{candual}
[\mV_i T_i]=[n,i][S_i], 
\end{eqnarray}
and the part~\eqref{Vn2} of the proposition follows directly from the formula
$v_i=[n,i]v^i$. 

To prove statement~\eqref{Vn3}
we fix the standard basis $\{b_i\}_{1\leq i\leq n}$ of $\mC^n$. Then there is an
isomorphism of vector spaces sending $b_i$ to $b_{n-i+1}$. This induces an
isomorphism $C^i\cong C^{n-i}$ of graded algebras, hence an equivalence of
categories 
\begin{eqnarray*}
    \hat{\sigma}_n:\quad \bigoplus_{i=0}^n C^i\gMOD&\cong&\bigoplus_{i=0}^n
    C^i\gMOD\\
S_i&\mapsto&S_{n-i},\\
\mV_i T_i&\mapsto&\mV_{n-i} T_{n-i}.
  \end{eqnarray*}   
The remaining isomorphisms of functors of part~\eqref{Vn3} follow then directly from the definitions. 

To prove statement~\eqref{Vn4}, recall first that $\mV_iP_i$ is self-dual,
hence we have $\hat\psi_n\mV P_i\cong\mV P_i$. The
isomorphisms  $\hat\psi_n\circ K^i\cong K^{-i}\circ\hat\psi_n$ and $\hat\psi_n\circ
\langle i\rangle\cong\langle -i\rangle\circ\hat{\sigma}_n$ follow directly from the
definitions. For the remaining isomorphisms of functors of part~\eqref{Vn4} we
refer to Lemma~\ref{dual}.
The semi-linearity of the form in part~\eqref{Vn5} follows directly from the formula $\HOM_\cC(M\langle
 k\rangle,N)=\HOM_\cC(M,N\langle -k\rangle)=\HOM_\cC(M,N)_k$ of graded vector
 spaces. On the other hand we have 
\begin{eqnarray*}
\sum_{i}\sum_{j}(-1)^j\DIM\op{Ext}^j_{\cC}(S_k,\mV_l T_l\langle
-i\rangle)q^i
&=&\sum_{i}\DIM\HOM_\cC(S_k,
\mV_l T_l\langle -i\rangle)q^i\\
&=&\delta_{k,l}q^{l(n-l)}=\delta_{k,l}q^{k(n-k)},
\end{eqnarray*}
since $\mV_l T_l$ is injective and its socle is concentrated in degree $l(n-l)$. From
formula \eqref{candual} it follows
$<[\mV_lT_l],[\mV_kT_k]>=\delta_{l,k}q^{l(n-l)}[n,l]$. Looking at formula~\eqref{semilin} completes the proof.         
It is left to prove statement~\eqref{Vn6}. From Lemma~\ref{adjunctions} and Lemma~\ref{free} we get the followings pairs
  of adjoint functors
  \begin{eqnarray*}
    &&\big(\RES^{i,i+1}_{i+1} C^{i,i+1}\otimes_{C^i}\langle -n+i+1\rangle,
    \RES^{i,i+1}_{i} C^{i,i+1}\otimes_{C^{i+1}}\langle -2i\rangle\langle
    n-i-1\rangle\big)\\
&=&\big(E_i, F_{i+1}\langle n-2(i+1)\rangle\langle 1\rangle\big).     
  \end{eqnarray*}
  and 
   \begin{eqnarray*}
    &&\big(\RES^{i,i-1}_{i-1} C^{i,i-1}\otimes_{C^i}\langle -i+1\rangle,
    \RES^{i,i-1}_{i} C^{i,i-1}\otimes_{C^{i-1}}\langle -2(n-i)\rangle\langle
    i-1\rangle\big)\\
&=&\big(F_i, E_{i-1}\langle 2(i-1)-n\rangle\langle 1\rangle\big).     
  \end{eqnarray*}
This gives the first two pairs of adjoint functors. The remaining ones are
obvious from the definitions. 
\end{proof}

The answer to the question raised at the end of the previous section follows
now directly: Theorem~\ref{cat} together with Theorem~\ref{Vn}~\eqref{Vn2} and Proposition~\ref{HCC} provide isomorphisms of $U_q(\SL_2)$-modules
  \begin{eqnarray*}
  \begin{array}[tb]{ccccc}
 V_n &\cong&{\bf G}(\bigoplus_{i=0}^n\cA_i^{-\rho})&\cong&{\bf
  G}(\bigoplus_{i=0}^n C^i\gMOD)\\
\\
   v^k&\mapsto&1\otimes [L_i]&\mapsto&1\otimes [S_i]\\
   v_k&\mapsto&1\otimes [P_i]&\mapsto&1\otimes [\mV_i T_i].
  \end{array}
\end{eqnarray*}
where $L_i$ denote a graded lift of $\cL(M(-\rho),L(w_o\cdot\omega_i))$ with
head concentrated in degree zero and $P_i$ denotes its projective cover.
 
\subsection{An elementary categorification of $\ov{V}_1^{\otimes n}$ using
  algebras of functions}
Next we would like to give a geometric categorification of the tensor products
$\ov{V}_{\bf d}$. Since we do not have a general version of Soergel's theory
which on the one hand side naturally leads to the algebras $C^i$ and on the
other hand generalises directly Proposition~\ref{HCC}, we will start from the
opposite end and propose a class of finite dimensional algebras whose graded
modules yield the desired geometric categorification.  

In this subsection we will do the first step
in this direction by giving a rather elementary categorification of $\ov
V_1^{\otimes n}$ using the algebras $B^i$ of functions on the finite set of
cosets $W/W_i$ (Proposition~\ref{naive}). We believe that this is a necessary
ingredient of a more substantial and general construction which will be
considered in the next subsection. 

Let $W=S_n$ be the symmetric
group of order $n!$ with subgroup $W_i$ as above. 
For any $0\leq i\leq n$ Let $B^i=\op{Func}(W/W_i)$
be the algebra of complex valued functions on the (finite) set
$W/W_i$. Similarly, for $0\leq i,i+1\leq n$ let
$B^{i,i+1}=\op{Func}(W/W_{i,i+1})$ be the algebra of functions on
$W/W_{i,i+1}$. For any $w\in W/W_{i}$ we have an idempotent $e^{(i)}_w\in
B^i$, namely the
characteristic function on $w$, \ie $e^{(i)}_w(x)=\delta_{w,x}$. In fact, the
$e^{(i)}_w$, $w\in W/W_i$ form a complete set of primitive, pairwise orthogonal,
idempotents. The algebra $B^i$ is semisimple with simple (projective) modules
$S_w^i=B^i e_w^{(i)}$. 
On the other hand $B^{i,i+1}$ is both, a $B^i$-module and a $B^{i+1}$-module
as follows: Because $W_{i,i+1}$ is a subgroup of $W_i$ and $W_{i+1}$ we have
surjections $\pi_i: W/W_{i,i+1}\rightarrow W/W_i$ and $\pi_{i+1}:W/W_{i,i+1}\rightarrow W/W_{i+1}$. If $g\in B^{j}$ for $j\in\{i,i+1\}$ and
$f\in B^{i,i+1}$ we put $g.f(x)=g(\pi_j(x))f(x)$ for $x\in W/W_{i,i+1}$. The $B^i$'s are commutative, hence we get a left and a
right module structure. Clearly, $B^{i,i+1}$ becomes a free $B^j$-module of rank
equal to the order of the group 
$({W/W_{i,i+1}})/({W/W_i})$, hence equal to the order of $W_i/W_{i,i+1}$. Let
$$\cC_{func}:=\bigoplus_{i=0}^{n} B^i\MOD.$$  For technical reasons, if $i>n$
or $i<0$, let $B^i\MOD$ denote the
category consisting of the zero $\mC$-module. We define the following endofunctors of $\cB$: 
\begin{itemize}
\item $E_{func}=\bigoplus_{i=0}^n E_i$, \\
where $E_i:B^i\MOD\rightarrow B^{i+1}\MOD$
  is the functor $B^{i,i+1}\otimes_{B^i}\bullet$ if $i<n$ and the zero functor otherwise. 
\item $F_{func}=\bigoplus_{i=0}^n F_i$, \\
where $F_i:B^i\MOD\rightarrow B^{i-1}\MOD$
  is the functor $B^{i,i-1}\otimes_{B^i}\bullet$ if $i>0$ and the zero functor
  otherwise. 
\end{itemize}

For any $w\in W$ we denote by ${\bf a}_{w,n,i}$ the $\{0,1\}$-sequence $w(1,
1,\ldots 1, 0,0,\ldots 0)$ of length $n$, where we used exactly $i$ ones. We get the
following elementary categorification: 
\begin{prop}
\label{naive}
  The category $\cC_{func}$ together with the isomorphism 
  \begin{eqnarray*}
    \eta:\quad\quad{\bf G}(\cC_{func})&\longrightarrow&\ov V_1^{\otimes n}\\
    S_w^i=B^i e_w^{(i)}&\longmapsto& v_{{\bf a}_{w,n,i}}, 
  \end{eqnarray*}
and the functors $E_{func}$ and $F_{func}$ is a categorification (in the sense of
Subsection~\ref{section1}) of the module $\overline{V}_1^{\otimes n}$.  
\end{prop}

\begin{proof}
  Clearly, the map $\eta$ is an isomorphism of vector spaces and the functors
  $E_{func}$ and $F_{func}$ are exact. As $B^{i,i+1}$ is a free $B^i$-module of
  rank equal to $n-i$, the order of $W_i/W_{i,i+1}$, it follows that the
  $B^{i+1}$-module $B^{i,i+1}\otimes_{B^{i}}S_w^i$ is of dimension $n-i$. Hence it is a direct sum of
  $n-i$ simple $B^{i+1}$-modules. A basis of $B^{i,i+1}\otimes_{B^{i}} S_w^i$ is given by
  elements of the form $f_x\otimes 1$, where $x\in W_i/W_{i,i+1}$ and $f_x$ is
  the characteristic function for $xw\in W_iw$. If $x\in W_i$ then
  $x{\bf a}(w,n,i+1)$ is equal to $x{\bf a}(w,n,i)$, but exactly one zero occurring in the
  sequence replaced by a one. If we allow only $x\in W_i/W_{i,i+1}$ then the 
  $x{\bf a}(w,n,i)$ provide each sequence exactly once. Therefore,
  $\eta(E_{func}S_w^i)=E\eta(S_w^i)$. Similarly, $\eta(F_{func}S_w^i)=F\eta(S_w^i)$. The
  statement follows.      
\end{proof}

Given any finite dimensional algebra, say $A$, we could equip $A$ with a
trivial $\mZ$-grading by putting $A=A_0$. Then a graded $A$-module is nothing
else than an $A$-module $M$ which carries the structure of a $\mZ$-graded
vector space. In particular, we could consider the function algebras $B^{i}$
as trivially graded. Then $B^{i,i+1}$ becomes a $\mZ$-graded
 $(B^{i+1},B^{i})$-bimodule by putting the characteristic function
 corresponding to the shortest coset representative $w\in W/W_{i,i+1}$ in degree
 $l(w)$. Similarly $B^{i,i-1}$ becomes a $\mZ$-graded
 $(B^{i-1},B^{i})$-bimodule. In this way, we get graded lifts 
 \begin{eqnarray*}
 {\bf E}_{func}&:&\quad\oplus_{i=0}^{n} B^i\gMOD\longrightarrow \oplus_{i=0}^{n}
 B^i\gMOD\\
{\bf F}_{func}&:&\quad\oplus_{i=0}^{n} B^i\gMOD\longrightarrow \oplus_{i=0}^{n}
 B^i\gMOD
\end{eqnarray*}
of our functors $E_{func}$ and $F_{func}$. We define ${\bf
  \cC}_{func}=\oplus_{i=0}^{n} B^i\gMOD$ with the endofunctor  
$\mathbf{K}=\oplus_{i=0}^{n}\langle
 2i-n\rangle$. The following statement follows directly from
 Proposition~\ref{naive} and the definition of the grading and provides a
  categorification of the $U_q(\mathfrak{sl}_2)$-module $V_1^{\otimes n}$ in
  terms of graded modules over our function algebras:

 \begin{corollary}
The isomorphism 
  \begin{eqnarray*}
    \eta:\quad\quad{\bf G}({\bf \cC}_{func})&\longrightarrow&V_1^{\otimes n}\\
    S_w^i=B^i e_w^{(i)}&\longmapsto& v_{{\bf a}_{w,n,i}}, 
  \end{eqnarray*}
defines an isomorphism of $U_q(\mathfrak{sl}_2)$-modules, where the module
structure on the left hand side is induced by the functors ${\bf E}_{func}$, ${\bf F}_{func}$ and ${\bf K}_{func}$. 
 \end{corollary}
  
\subsection{A categorification of $\ov V_{\bf d}$ using finite dimensional algebras}
\label{between}

Now we combine the categorification of irreducible modules from
Section~\ref{sectionVn} and the elementary categorification from
Proposition~\ref{naive} to a (partly conjectural) categorification of an arbitrary tensor product
$\ov V_{\bf d}$. The present construction is parallel to the categorification
of irreducible modules as described in Section~\ref{sectionVn} and coincides
with this categorification in the special case when $\ov V_{\bf d}$ has a
single factor. 

Let again $W=S_n$ with the Young subgroup $S_{\bf d}$
corresponding to the composition ${\bf d}$. Let $B^{\bf d}=\op{Func}(W/S_{\bf
  d})$ be the algebra of (complex values) functions on
$W/S_{\bf d}$. Recall (from Section~\ref{section1}) the subalgebras $C^i$,
$C^{i,i+1}$, $C^{i,i-1}$ in the coinvariant algebra corresponding to
$W$. The Weyl group $W$ is acting on both, $B^{\bf d}$ and $C$.

We set 
\begin{eqnarray}
  \label{eq:BC}
  H_{\bf d}^i=(B^{\bf d}\otimes C)^{W_i},\quad\quad\quad H_{\bf d}^{i,i+1}=(B^{\bf d}\otimes C)^{W_{i,i+1}},
\end{eqnarray}
where we take the $W_i$-invariants with respect to the diagonal action. For
any $w\in W_i\backslash W/ S_{\bf d}$, there is an idempotent 
$f_w=e_w\otimes 1$, where $e_w(x)=e_w(yx)=\delta_{w,x}$ for any $x\in W_i\backslash W/S_{\bf d}$ and $y\in W_i$. These $f_w$ form in fact a complete set
of primitive pairwise orthogonal idempotents. In particular, the simple
modules of $H_{\bf d}^i$ are naturally indexed by (longest coset
representatives of) the double cosets
$W_i\backslash W/ S_{\bf
  d}$. Let $P_{{\bf d},x}^i=H_{\bf d}^if_w$ be the corresponding
indecomposable projective module with simple head $S_{{\bf d},x}^i$. We define
$$\cC_{geom}:=\bigoplus_{i=0}^{n} H_{\bf d}^i\MOD.$$ 
Obviously, $H_{\bf d}^j$ is a subset of  $H_{\bf d}^{i,i+1}$ for $j=i, i+1$
(if they are defined). For technical reasons we denote by $H_{\bf d}^i\MOD$
the category containing only the zero $\mC$-module if $i>n$ or $i<0$.

Analogous to our elementary construction we define for $0\leq i\leq n$ the functors 
\begin{eqnarray*}
 E_i:&& H_{\bf d}^i\MOD\rightarrow H_{\bf d}^{i+1}\MOD, \quad\text{   as    }\quad
 H_{\bf d}^{i+1}\otimes_{ H_{\bf d}^i}\bullet\quad\text{if $0\leq i<n$},\\
 F_i:&& H_{\bf d}^i\MOD \rightarrow H_{\bf d}^{i-1}\MOD, \quad\text{   as
 }\quad H_{\bf d}^{i-1}\otimes_{ H_{\bf d}^i}\bullet\quad\text{if $0< i\leq
 n$},
\end{eqnarray*}
otherwise it should be just the zero functor. We set
\begin{eqnarray*}
E_{geom}=\bigoplus_{i}^n E_i,&\quad&F_{geom}=\bigoplus_{i}^n F_i,
\end{eqnarray*}

\begin{prop}
  \begin{enumerate}
There are isomorphisms of vector spaces
  \begin{eqnarray*}
\Phi_1:\quad{\bf G}(\cC_{geom})&\longrightarrow&\ov V_{\bf d}\\
    S_{{\bf d},w}^i&\longmapsto& v^{{\bf a}(\mu)}. \\
\Phi_2:\quad {\bf G}(\cC_{geom})&\longrightarrow&\ov V_{\bf d}\\
    P_{{\bf d},w}^i&\longmapsto& v_{{\bf a}(\mu)}.
  \end{eqnarray*}
where ${{\bf a}(\mu)}=w(1,\ldots 1,0,\ldots 0)$. 
\end{enumerate}
\end{prop}

\begin{proof}
  This follows directly from the definitions.
\end{proof}
Note that if the tensor product $\ov V_{\bf d}$ has only one factor, that is
$\mu=(n)$ and hence $S_{\bf d}=W$ we get $H_{\bf d}^i\cong C^i$ and the functors
$E^{\cG}$ and $F^{\cG}$ become the functors $E$ and $F$ from
Section~\ref{sectionVn}. For the general case we would like to formulate the following  
\begin{conjecture}
  The isomorphisms $\Phi_1$ and $\Phi_2$ agree and are isomorphisms of $\mathfrak{sl}_2$-modules,
  where the action on the left hand side is induced by the functors $E_{geom}$ and $F_{geom}$.
\end{conjecture}
The quantum version of this 
  conjectural categorification should again arise from the corresponding
  graded version.

\begin{remark}{\rm 
  The Conjecture implies in particular, that the functors $E_{geom}$ and
  $F_{geom}$ preserve the additive category of projective modules. By direct
  calculations it can be shown that the conjecture is true for all cases where
  $n=2, 3$. In these cases we also know that $(B^{\bf d}\otimes C)^{W'}$ is a
  free $(B^{\bf d}\otimes C)^{W''}$-module of rank $|W''/W'|$ for any
  subgroups $W\supseteq W''\supseteq W'$.} 
\end{remark}

\subsection{Open problems related to a geometric categorification}
\label{BorelMoore}
The categorification of $\ov{V}_{\bf d}$ via the modules over the finite
dimensional algebras $H_{\bf d}^i$ from the previous section strongly suggest
that the geometry of Grassmannians and partial flag varieties used in the case
of a single factor $\overline{V}_n$ should be replaced by the geometry of generalised Steinberg
varieties (as defined in \cite{RD1}) for the general linear group $GL_n(\mC)$
to obtain a geometric categorification of arbitrary tensor products $\ov{V}_{\bf
  d}$. Note first that the dimension of the algebra $B\otimes C$ coincides
with the dimension of the Borel-Moore homology $H_*(Z)$ of the (full)
Steinberg variety for $GL_n(\mC)$ and is equal to $|W|^2$ (see
e.g. \cite[Proposition 8.1.5, Lemma 7.2.11]{CG}). Moreover, the algebras
$H_{\bf d}^i$ can be viewed as the subalgebra of $W_{\bf d}\otimes
W_i$-invariants in $B\otimes C$, and it was proven in \cite[(1.1'')]{RD2} that the
Borel-Moore homology of the generalised Steinberg variety $Z_{\bf d}^i$
associated to the pair $(W_{\bf {d}}, W_i)$ is isomorphic to the subspace of
$H_*(Z)$ given by $W_{\bf d}\otimes
W_i$-invariants. We expect that using intersection theory, (see \cite{Fulton}) one can define a
commutative algebra structure on $Z_{\bf d}^i$ which yields the algebras
$H_{\bf d}^i$ introduced in Section~\ref{between}.

We also note that the generalised Steinberg varieties $Z_{\bf
  d}^i$ are precisely the tensor product varieties of Malkin (\cite{Malkin})
and Nakajima (\cite{Nakajima}) in the special case of tensor products of
finite dimensional irreducible $\cU_q(\mathfrak{sl_2})$-modules. Thus the
geometric categorification we propose is also a natural next step in the
geometric description of the tensor products $V_{\bf d}$ from \cite{Savage}. 

Looking from the algebraic categorification side of the picture one notices
that the generalised Steinberg varieties also appear as characteristic
varieties of Harish-Chandra bimodules (\cite{Borho}). Thus one expects that the geometry of
characteristic varieties that underlies our Categorification Theorem~\ref{cat}
should provide a conceptual relation between the algebraic and the geometric
categorifications giving rise to isomorphisms of
$\cU(\mathfrak{sl}_2)$ modules:
\begin{eqnarray*}
  \mathbf{G}\big(\bigoplus_{i=0}^n\cA_i^\mu(\mathfrak{gl}_n)\big)\cong\mathbf{G}(\bigoplus_{i=0}^n
  H_{\bf d}^i\MOD)\cong V_{\bf d}.
\end{eqnarray*}

One can also relate the algebraic categorification with a geometric
categorification by studying the projective functors acting on the category
$\cO(\mathfrak{gl}_n)$, extending our constructions in
Section~\ref{Ograd}. In fact the complexified Grothendieck ring of
projective endofunctors of the principal block of $\cO(\mathfrak{gl}_n)$  is
isomorphic to the group algebra of the Weyl group $W$ by the classification
theorem from \cite{BG}  and the
Kazhdan-Lusztig theory. On the other hand, the group algebra of $W$ is
canonically isomorphic to the top degree $H_{top}(Z)$ of the Borel-Moore
homology $H_*(Z)$ with respect to the convolution product (see e.g. \cite[Theorem 3.4.1]{CG}). One can show that the whole Borel Moore homology ring
$H_*(Z)$ encodes the Grothendieck ring of these projective functors together with the
action of the centre of the category. This picture can be generalised to
$H_{top}(Z_{\bf d}^i)$ and $H_*(Z_{\bf d}^i)$ by looking at projective functors
between different singular blocks of the category
$\cO(\mathfrak{gl}_n)$. Details of this alternative approach will appear in a
subsequent paper.    
\bibliography{ref}

\newpage
Igor Frenkel, \\
Yale University, Department of Mathematics,
10 Hillhouse Avenue, PO Box 208283, New Haven, Connecticut 06520-8283.\\
e-mail: {\tt frenkel-igor@yale.edu} \\

Mikhail Khovanov,\\ 
Department of Mathematics, Columbia University, New York, NY
10027.\\
 email: {\tt khovanov\symbol{64}math.columbia.edu}\\

Catharina Stroppel, \\
Department of Mathematics,
University of Glasgow, University Gardens,
Glasgow G12 8QW, United Kingdom.\\
e-mail: {\tt cs\symbol{64}maths.gla.ac.uk}

\end{document}